\documentclass[12pt,reqno]{amsproc}
\usepackage{amsmath,exscale,amssymb}
\usepackage{float}
\usepackage{graphicx}
\usepackage{amscd}
\usepackage{amsthm,amsxtra}
\usepackage[size=2em]{diagrams}%,PostScript=ghostscript]{diagrams}

\theoremstyle{plain}
\newtheorem{X}{X}[section]
\newtheorem{theorem}[X]{Theorem}
\newtheorem{proposition}[X]{Proposition}
\newtheorem{lemma}[X]{Lemma}
\newtheorem{corollary}[X]{Corollary}

\newtheorem{E}[X]{}

\theoremstyle{definition}
\newtheorem{definition}[X]{Definition}
\newtheorem{example}[X]{Example}
\newtheorem{plain}[X]{}

\newtheorem{remark}[X]{Remark}

\newenvironment{mylist}
{\begin{list} {--}
{\setlength{\leftmargin}{.5in}\setlength{\rightmargin}{.5in}}}
{\end{list}}

\renewenvironment{itemize}
{\begin{mylist}} {\end{mylist}}

\def\1{{1\mkern-7mu1}}  %% since alas there's no "\Bbb{1}";

\newcommand\Aut{\operatorname{Aut}}

\newcommand\Coker{\operatorname{Coker}}
\newcommand\Corr{\operatorname{Corr}}

\newcommand\End{\operatorname{End}}
\newcommand\Ext{\operatorname{Ext}}

\newcommand\Gal{\operatorname{Gal}}

\newcommand\Hom{\operatorname{Hom}}

\newcommand\id{\operatorname{id}}
\newcommand\Ind{\operatorname{Ind}}

\newcommand\im{\operatorname{Im}}  %\Im is already used by TeX.

\newcommand\Ker{\operatorname{Ker}}

\newcommand\Map{\operatorname{Map}}

\newcommand\ob{\operatorname{ob}}
\newcommand\ord{\operatorname{ord}}

\newcommand\rank{\operatorname{rank}}

\newcommand\Rep{\operatorname{Rep}}

\newcommand\Spec{\operatorname{Spec}}

\newcommand\Tor{\operatorname{Tor}}
\newcommand\Tr{\operatorname{Tr}}

\def\Ab{{\mathbf{Ab}}}%abelian groups
\def\CM{{\bold{CM}}}%CM-motives
\def\Crys{{\bold{Crys}}}%crystals
\def\Et{{\bold{Et}}}%Etale algebras.

\def\Isoc{{\bold{Isoc}}}%Isocrystals
\def\Mod{{\bold{Mod}}}%modules
\def\Modf{{\bold{Modf}}}%modules of finite type.

\def\Rep{{\bold{Rep}}}%representations
\def\Vct{{\bold{Vec}}}
\def\ord{{\mathrm{ord}}}

\let\xr=\xrightarrow

\def\tsum{\mathop{\textstyle \sum }}%
\def\tprod{\mathop{\textstyle \prod }}%
\def\tcoprod{\mathop{\textstyle \coprod }}%
\begin{document}
\title{Integral Motives and Special Values of Zeta Functions}
\author{James S. Milne}
\address{2679 Bedford Rd., Ann Arbor, MI 48104, USA.}
\email{math@jmilne.org}\urladdr{www.jmilne.org/math/}
\author{Niranjan Ramachandran}
\address{Mathematics Dept, University of Maryland, College Park, MD 20742,
USA.}
\email{atma@math.umd.edu}
\thanks{May 22, 2002; submitted version.\\
The first author received support from the National Science Foundation and
the second author from MPIM (Bonn) and a GRB Summer Grant (UMD)}

\begin{abstract}
For each field $k$, we define a category of rationally decomposed mixed
motives with $\mathbb{Z}{}$-coefficients. When $k$ is finite, we show that
the category is Tannakian, and we prove formulas relating the behaviour of
zeta functions near integers to certain $\Ext$ groups.
\end{abstract}

\maketitle
\tableofcontents

\subsection{Introduction}

For each field $k$, there is conjecturally a category of mixed isomotives
whose full subcategory of semisimple objects is the category $\mathcal{M}%
{}(k;\mathbb{Q}{})$ of pure isomotives. In this paper, we define a $\mathbb{Z%
}$-category $\mathcal{M}{}(k;\mathbb{\mathbb{Z}{}}$) of motives whose
quotient by its subcategory of torsion objects is $\mathcal{M}{}(k;\mathbb{Q}%
{})$. Thus, $\mathcal{M}{}(k;\mathbb{Z}{})$ is the abelian category of mixed
motives over $k$ whose weight filtrations split modulo torsion. As
Grothendieck observed, when $k$ is finite, the weight filtrations always
split modulo torsion, and so, in this case, $\mathcal{M}{}(k;\mathbb{Z}{})$
is the full category of mixed motives over $k$.

The category $\mathcal{M}{}(k;\mathbb{Z}{})$ depends on our choice of $%
\mathcal{M}{}(k;\mathbb{Q})$. For definiteness, take $\mathcal{M}{}(k;%
\mathbb{Q}{})$ to be the category of isomotives based on the smooth
projective varieties over $k$ whose K\"{u}nneth projectors are algebraic
with the correspondences taken to be the numerical equivalence classes of
algebraic classes. This category is abelian and semisimple (Jannsen 1992),
hence Tannakian (Deligne 1990), and when we assume\footnote{%
See (\ref{cm3}) for an alternative choice of $\mathcal{M}{}(k;\mathbb{Q}{})$
that does not require any unproven conjectures.} numerical equivalence
coincides with homological equivalence, it admits canonical $l$-adic fibre
functors for each $l$. The category $\mathcal{M}{}(k;\mathbb{Z}{})$ we
construct is noetherian, abelian, and, at least when $k$ is finite,
Tannakian.

The $\Ext$ groups in the full category $\mathcal{M}{}(k;\mathbb{Z}{})$ seem
to be pathological, but those in the subcategory $\mathcal{M}{}^{+}(k;%
\mathbb{Z}{})$ of effective motives are of interest. When $k$ is finite, we
show that the Tate conjecture implies\footnote{%
When $\mathcal{M}{}(k;\mathbb{Q})$ is as in (\ref{cm3}), these relations are
hold unconditionally.} relations between these $\Ext$ groups and the
behaviour of zeta functions near integers (Theorems \ref{gc.a}, \ref{gc.i}).
Moreover, to a smooth projective variety $V$ satisfying a certain
degeneration condition on its crystalline cohomology, we are able to attach
\textquotedblleft Weil\textquotedblright\ motivic cohomology groups $H_{%
\text{mot}}^{i}(V,\mathbb{Z}{}(r))$ satisfying the spectral sequence%
\begin{equation}
H^{i}(\Gamma _{0},H_{\text{mot}}^{j}(\bar{V},\mathbb{Z}{}(r))\Longrightarrow
H_{\text{mot}}^{i+j}(V,\mathbb{Z}{}(r))  \label{e29}
\end{equation}%
posited by Lichtenbaum (2002), and deduce from Theorem \ref{gc.i}\ that%
\begin{equation}
\zeta (V,s)\sim \pm \chi ^{\times }(V,\mathbb{Z}{}(r))\cdot q^{\chi (V,%
\mathcal{O}{}_{V},r)}\cdot (1-q^{r-s})^{\rho }\text{ as }s\rightarrow r
\label{e21}
\end{equation}%
(Theorem \ref{gc.n}). Here $\Gamma _{0}$ is the subgroup of $\Gal(\mathbb{F}%
{}/\mathbb{F}_{q})$ generated by the Frobenius element (so $\Gamma _{0}\cong
\mathbb{Z}{}$), $\chi ^{\times }(V,\mathbb{Z}{}(r))$ is the alternating
product of the orders of the cohomology groups of a complex%
\begin{equation*}
\cdots \rightarrow H_{\text{mot}}^{i}(V,\mathbb{Z}{}(r))\rightarrow H_{\text{%
mot}}^{i+1}(V,\mathbb{Z}{}(r))\rightarrow \cdots
\end{equation*}%
arising from (\ref{e29}), $\chi (V,\mathcal{O}_{{}V},r)$ is as in Milne
1986a, and $\rho =\tsum (-1)^{i}i\cdot \rank H_{\text{mot}}^{i}(V,\mathbb{Z}%
{}(r))$.

The category $\mathcal{M}{}(k;\mathbb{Z}{})$ is defined in \S 5 as a full
subcategory of the fibre product of the $\mathbb{Q}{}$-category $\mathcal{M}%
{}(k;\mathbb{Q}{})$ with a certain $\mathbb{\hat{Z}}{}$-category of
realizations. The construction requires that we show that the various $l$%
-adic fibre functors on $\mathcal{M}{}(k;\mathbb{Q}{})$ define an ad\`{e}lic
fibre functor --- the proof of this depends crucially on the theorem of
Gabber (1983). In \S 4, we show that the usual Tate conjecture implies an ad%
\`{e}lic version. For example, for a variety over a finite field, if the $l$%
-adic Tate conjecture holds for a single $l$, then it holds $l$-integrally
for almost all${}$ $l$.

The category $\mathcal{M}{}(\mathbb{F}{}_{q};\mathbb{Z}{})$ is shown to be
Tannakian in \S 6, and the main results on Exts${}$ and motivic cohomology
groups are obtained in \S 10 after preliminaries in \S \S 7--9.

This article is part of a project whose ultimate goal is to define a
triangulated category of motives with $t$-structure whose heart is the
category $\mathcal{M}{}(k;\mathbb{Z}{})$ defined above, and to show that the
Tate conjecture for smooth projective varieties implies the formula (\ref%
{e21}) for \emph{arbitrary }varieties over $\mathbb{F}{}_{q}$.

We note that our approach to motives is opposite to that of other authors
(Hanamura 1995, 1999; Levine 1998; Voevodsky 2000) who define a triangulated
category of motives with the goal of putting a $t$-structure on it whose
heart will be the abelian category of motives. While our approach can be
expected to give the full category of mixed motives only when the ground
field is $\mathbb{F}{}_{q}$ or $\mathbb{F}{}$, we do expect to obtain the
correct $p$-structure, as illustrated by (\ref{e21}). Most of the
difficulties in the present work involve the $p$-part in characteristic $p$.

\subsection{Notations\label{aaa}}

Throughout, $l$ is a prime number, $p$ is the characteristic (or
characteristic exponent) of the ground field, and $\ell $ is a prime number $%
\neq p$. The $l$-adic valuations $|\quad |_{l}$ and $\ord_{l}$ are
normalized so that $|l|_{l}=l^{-1}$ and $\ord_{l}(l)=1$. The expression
\textquotedblleft for almost all $l$\textquotedblright\ means
\textquotedblleft for all but (possibly) a finite number of $l$%
\textquotedblright .

The profinite completion of $\mathbb{Z}{}$ is denoted $\mathbb{\hat{Z}}{}$
and $\mathbb{A}_{f}=\mathbb{\mathbb{Q}{}}\otimes \mathbb{\hat{Z}}$. Thus, $%
\mathbb{\hat{Z}}=\prod \mathbb{Z}{}_{l}$, and we let $\mathbb{Z}%
{}^{p}=\prod_{l\neq p}\mathbb{Z}{}_{l}$ and $\mathbb{A}_{f}^{p}=\mathbb{Q}%
{}\otimes \mathbb{Z}{}^{p}$. When $p=1$, $\mathbb{Z}{}^{p}=\mathbb{\hat{Z}}%
{} $ and $\mathbb{A}_{f}^{p}=\mathbb{A}_{f}$.

We use $M_{m}$ and $M^{(m)}$ respectively to denote the kernel and cokernel
of multiplication by $m$ on $M$. For a prime $l$, $M(l)=\bigcup M_{l^{n}}$
is the $l$-primary component of $M$, and $M_{\text{tors}}=$ $\bigcup M_{m}$
is the torsion subobject.

For a module $M$ endowed with an action of a group $\Gamma $, $M^{\Gamma }$
(resp. $M_{\Gamma }$) is the largest submodule (resp. quotient module) of a $%
\Gamma $-module $M$ on which $\Gamma $ acts trivially.

Lattice always means full lattice: when $R$ is a ring and $V$ is a free $%
\mathbb{Q}{}\otimes R$-module of finite rank, an $R$-lattice in $V$ is a
finitely generated $R$-submodule $\Lambda $ of $V$ such that the inclusion $%
\Lambda \hookrightarrow V$ induces an isomorphism $\mathbb{Q}{}\otimes
\Lambda \rightarrow V$.

For a field $k$, $k^{\text{al}}$ denotes an algebraic closure of $k$, $k^{%
\text{sep}}$ a separable (algebraic) closure of $k$, and $k^{\text{pf}}$ a
perfect closure of $k$ (for example, $k^{\text{pf}}=k^{p^{-\infty }}\subset
k^{\text{al}}$).

For a perfect field $k$ of characteristic $p$, $W(k)$ denotes the ring of
Witt vectors with coefficients in $k$, and $B(k)$ denotes the field of
fractions of $W(k)$. The unique lift of the automorphism $x\mapsto
x^{p}\colon k\rightarrow k$ to $W(k)$ is denoted $\sigma $. We let $A$
denote the skew polynomial ring $W[F,\sigma ]$ (polynomials in $F$ with the
relation $F\cdot a=\sigma (a)\cdot F$, $a\in W$).

Generally, rings are assumed to be commutative. For a ring $R$, $\Mod(R)$ is
the category of $R$-modules and $\Modf(R)$ the category of finitely
presented $R$-modules, and for a field $k$, $\Vct(k)$ is the category of
finite dimensional $k$-vector spaces. We usually abbreviate $M\otimes _{R}S$
to $M_{S}$. For an $R$-linear category (alias, $R$-category) $\mathcal{M}{}$%
, $\mathcal{M}{}_{S}$ denotes the category with the same objects as $%
\mathcal{M}{}$, but with $\Hom_{\mathcal{M}{}_{S}}(X,Y)=\Hom_{\mathcal{M}%
{}}(X,Y)_{S}$.

For an affine group scheme $G$ over a noetherian ring $R$, $\Rep(G;R)$
denotes the category of linear representations of $G$ on finitely generated $%
R$-modules.

Given functors $F\colon \mathcal{A}{}\rightarrow \mathcal{C}{}$ and $G\colon
\mathcal{B}{}\rightarrow \mathcal{C}{}$, the \emph{fibre product category }$%
\mathcal{A}{}\times _{\mathcal{C}{}}\mathcal{B}{}$ has as objects the
triples $(A,B,\gamma )$ with $A$ an object of $\mathcal{A}{}$, $B$ an object
of $\mathcal{B}$, and $\gamma \colon F(A)\rightarrow G(B)$ an isomorphism;
the morphisms are the pairs $a\colon A\rightarrow A^{\prime }$ and $b\colon
B\rightarrow B^{\prime }$ giving rise to a commutative diagram in $\mathcal{C%
}{}$.

A variety $V$ is a geometrically reduced separated scheme of finite type
over a field. The group of algebraic cycles of codimension $r$ on $V$ is
denoted $Z^{r}(V)$, and $Z_{\sim }^{r}(V)$ is the quotient of $Z^{r}(V)$ by
some adequate relation $\sim $. We let $H^{\ast }(V,F)=\oplus _{i}H^{i}(V,F)$%
.

An equivalence class containing $x$ is denoted $[x]$. We also use $[S]$ to
denote the cardinality of a finite set $S$. For a homomorphism $f\colon
M\rightarrow N$ of abelian groups whose kernel and cokernel are finite, $%
z(f) $ is defined to be $[\Ker(f)]/[\Coker(f)]$.

Isomorphisms are denoted $\approx $ and canonical isomorphisms $\cong $.

\section{The Realization Categories}

\subsection{The categories $\mathcal{R}$ away from $p$}

Throughout this subsection, $k$ is a field of characteristic exponent $p$,
and $R$ is a topological ring, for example, $R=\mathbb{Z}{}$, $\mathbb{Z}%
_{\ell }$, $\mathbb{Q}_{\ell }$, $\mathbb{Z}{}^{p}$, or $\mathbb{A}_{f}^{p}$
with its natural topology (discrete for $\mathbb{Z}{}$).

\begin{definition}
\label{pl5}(a) For a profinite group $\Gamma $, $\mathcal{R}{}(\Gamma ;R)$
is the category of continuous linear representations of $\Gamma $ on
finitely presented $R$-modules.

(b) When $k$ is finitely generated over the prime field,
\begin{equation*}
\mathcal{R}{}(k;R)=\mathcal{R}{}(\Gamma ;R),\quad \Gamma =\Gal(k^{\text{sep}%
}/k).
\end{equation*}%
Otherwise, $\mathcal{R}{}(k;R)=\varinjlim\mathcal{R}{}(k^{\prime };R)$ where $%
k^{\prime }$ runs over the subfields of $k$ finitely generated over the
prime field.
\end{definition}

\begin{remark}
\label{pl7}Let $\Gamma $ be a profinite group. Two continuous linear
representations $\rho _{1}$ and $\rho _{2}$ of open subgroups $U_{1}$ and $%
U_{2}$ of $\Gamma $ on an $R$-module $M$ are said to be \emph{related}\/ if
they agree on an open subgroup of $U_{1}\cap U_{2}$. This is an equivalence
relation, and an equivalence class will be called a \emph{germ} of a
continuous linear representation\/ of $\Gamma $. A \emph{morphism} $(M,[\rho
_{U}])\rightarrow (M^{\prime },[\rho _{U^{\prime }}^{\prime }])$ \emph{of
germs} is an $R$-linear map $\alpha \colon M\rightarrow M^{\prime }$ that is
equivariant for some open subgroup of $U\cap U^{\prime }$. When $k$ is a
separable closure of a field $k_{0}$ finitely generated over the prime
field, an object of $\mathcal{R}{}(k;R)$ can be identified with a germ of a
continuous linear representation of $\Gal(k/k_{0})$.
\end{remark}

\begin{lemma}
\label{pl11}For each object $M$ of $\mathcal{R}{}(\Gamma ;\mathbb{Z}_{\ell
}) $, there is an open subgroup $U\subset \Gamma $ such that $M$ decomposes
into the sum of a torsion object and a torsion-free object in $\mathcal{R}%
{}(U;\mathbb{Z}_{\ell })$. Consequently, when $k$ is separably closed, for
each $M$ in $\mathcal{R}{}(k;\mathbb{Z}{}_{\ell })$, $M_{\text{tors}}$ is a
direct summand of $M$.
\end{lemma}

\begin{proof}
Let $M=M_{\text{tors}}\oplus M_{1}$ as $\mathbb{Z}_{\ell }$-modules, and let
$\gamma \in \Gamma $ act on $M$ as the matrix $\left(
\begin{smallmatrix}
a(\gamma ) & b(\gamma ) \\
0 & c(\gamma )%
\end{smallmatrix}%
\right) $. The map $\gamma \mapsto a(\gamma )$ is the homomorphism
describing the action of $\Gamma $ on $M_{\text{tors}}$. After replacing $%
\Gamma $ with an open subgroup $U$, we may suppose this action to be
trivial. The map $\gamma \mapsto c(\gamma )$ defines an action of $U$ on $%
M_{1}$, and hence an action of $U$ on $\Hom(M_{1},M_{\text{tors}})$, which,
again, we may suppose to be trivial. Then $\gamma \mapsto b(\gamma )$ is a
homomorphism $U\rightarrow \Hom(M_{1},M_{\text{tors}})$, and after we
replace $U$ with the kernel of this map, $M=M_{\text{tors}}\oplus M_{1}$
will be a decomposition of $U$-modules.
\end{proof}

For $M$ in $\mathcal{R}{}(\Gamma ;\mathbb{Z}{}_{\ell })$, set $M^{\vee }=\Hom%
_{\mathbb{Z}{}_{\ell }}(M,\mathbb{Z}{}_{\ell })$ with $\gamma \in \Gamma $
acting according to the rule $(\gamma f)(m)=f(\gamma ^{-1}m)$. We say $M$ is
\emph{reflexive }if the canonical map $M\rightarrow M^{\vee \vee }$ is an
isomorphism.

\begin{proposition}
\label{pl10}\emph{(a)} An object of $\mathcal{R}{}(\Gamma ;\mathbb{Z}%
{}_{\ell })$ is reflexive if and only if it is torsion-free.

\emph{(b)} Every object of $\mathcal{R}{}(\Gamma ;\mathbb{Z}_{\ell })$ is a
quotient of a reflexive object.
\end{proposition}

\begin{proof}
(a) As $M^{\vee }$ is torsion-free, the condition is necessary, and it is
obviously sufficient.

(b) If $M$ decomposes in $\mathcal{R}{}(\Gamma ;\mathbb{Z}{}_{\ell })$ into
the direct sum of a torsion module $M_{t}$ with trivial $\Gamma $-action and
a torsion-free module $M_{1}$, then $M$ will be the quotient of $\mathbb{Z}%
{}_{\ell }^{r}\oplus M_{1}$ for some $r$. According the lemma, there is an
open subgroup $U$ of $\Gamma $ such that $M$ has such a decomposition in $%
\mathcal{R}{}(U;\mathbb{Z}{}_{\ell })$, and hence there will be a surjection
$\mathbb{Z}{}_{\ell }^{r}\oplus M_{1}\rightarrow M$ in $\mathcal{R}{}(U;%
\mathbb{Z}{}_{\ell })$ with $M_{1}$ free. Now%
\begin{equation*}
\Ind_{U}^{\Gamma }(\mathbb{Z}{}_{\ell }^{r}\oplus M_{1})\twoheadrightarrow %
\Ind_{U}^{\Gamma }(M)\xr{\varphi\mapsto\sum_{s\in\Gamma/U} s\varphi(s^{-1})}M
\end{equation*}%
realizes $M$ as the quotient of a torsion-free $\Gamma $-module. Here $\Ind%
_{U}^{\Gamma }(N)$ denotes the induced module.
\end{proof}

\begin{remark}
\label{pl10m}(a) The category $\mathcal{R}{}(k;\mathbb{Q}{}_{\ell })$ is
denoted \underline{$\mathrm{Tate}$}$(k)$ in Saavedra 1972, VI A4.2, A4.3,
and its objects are called Tate modules. It is a neutral Tannakian category
over $\mathbb{Q}{}_{\ell }$.

(b) The functor $M\mapsto (M_{\ell })_{\ell \neq p}\colon \mathcal{R}%
{}(\Gamma ;\mathbb{Z}{}^{p})\rightarrow \prod_{\ell \neq p}\mathcal{R}%
{}(\Gamma ;\mathbb{Z}_{l})$ is exact, full, and faithful, with essential
image the objects $(M_{\ell })_{\ell \neq p}$ such that $\dim _{\mathbb{F}%
{}_{\ell }}M_{\ell }/\ell M_{\ell }$ is bounded.
\end{remark}

\subsection{The categories $\mathcal{R}{}$ at $p$}

An $F$\emph{-isocrystal} \emph{over} \emph{a perfect field} $k$ is a finite
dimensional $B(k)$-vector space $M$ endowed with a $\sigma $-linear
bijection $F_{M}\colon M\rightarrow M$. With the obvious structures, the $F$%
-isocrystals over $k$ form a Tannakian category $\Isoc(k)$ over $\mathbb{Q}%
{}_{p}$, and the forgetful functor is a $B(k)$-valued fibre functor
(Saavedra 1972, VI 3.2). The identity object is $\1=(B(k),\sigma )$, and the
Tate object is $\mathbb{T}=(B(k),p^{-1}\cdot \sigma )$.

An $F$-isocrystal $(M,F_{M})$ is \emph{effective }if $F_{M}$ stabilizes a $%
W(k)$-lattice in $M$. Let $\Isoc^{+}(k)$ denote the full subcategory of $%
\Isoc(k)$ whose objects are the effective $F$-isocrystals. Every object of $%
\Isoc(k)$ is of the form $X\otimes \mathbb{T}^{\otimes m}$ with $X$
effective. In other words, for each $(M,F_{M})$ in $\Isoc(k)$, there exists
an $m\in \mathbb{Z}{}$ such that $p^{m}F_{M}$ stabilizes a lattice (Saavedra
1972, VI 3.1.3, 3.2.1).

An $F$\emph{-crystal\footnote{%
Some authors use \textquotedblleft crystal\textquotedblright\ to mean free
crystal, i.e., an $M$ that is torsion-free (as a $W(k)$-module).} over a
perfect field }$k$ is a finitely generated $W(k)$-module $\Lambda $ endowed
with a $\sigma $-linear map $F_{\Lambda }\colon \Lambda \rightarrow \Lambda $
such that the kernel of $F_{\Lambda }$ is torsion. The $F$-crystals over $k$
form an abelian tensor category $\Crys^{+}(k)$ over $\mathbb{Z}_{p}$. For an
$F$-crystal $(\Lambda ,F_{\Lambda })$, $(\Lambda ,F_{\Lambda })_{\mathbb{Q}%
{}}=_{\text{df}}(\Lambda _{\mathbb{Q}{}},F_{\Lambda }\otimes 1_{\mathbb{Q}%
{}})$ is an $F$-isocrystal, and the functor
\begin{equation*}
(\Lambda ,F_{\Lambda })\mapsto (\Lambda ,F_{\Lambda })_{\mathbb{Q}{}}\colon %
\Crys^{+}(k)_{\mathbb{Q}{}}\rightarrow \Isoc(k)
\end{equation*}%
is fully faithful with essential image $\Isoc^{+}(k)$.

Let $\mathbb{L}=(W(k),p\cdot \sigma )$ be the Lefschetz $F$-crystal. The
functor%
\begin{equation*}
-\otimes \mathbb{L}\colon \Crys^{+}(k)\rightarrow \Crys^{+}(k)
\end{equation*}%
is faithful (full and faithful on torsion-free objects). We define $\Crys%
\left( k\right) $ to be the tensor category obtained from $\Crys^{+}(k)$ by
inverting $\mathbb{L}$. Thus, the objects of $\Crys\left( k\right) $ are
pairs $(\Lambda ,m)$ with $\Lambda $ an $F$-crystal and $m\in \mathbb{Z}{}$,
and%
\begin{equation*}
\Hom((\Lambda ,m),(\Lambda ^{\prime },m^{\prime
}))=\varinjlim_{N\geq m,m^{\prime }}\Hom(\Lambda \otimes
\mathbb{L}^{\otimes N-m},\Lambda ^{\prime }\otimes
\mathbb{L}^{\otimes N-m^{\prime }}).
\end{equation*}%
The tensor product is defined by%
\begin{equation*}
(\Lambda ,m)\otimes (\Lambda ^{\prime },m^{\prime })=(\Lambda \otimes
_{W(k)}\Lambda ^{\prime },m+m^{\prime })\text{.}
\end{equation*}%
The \emph{Tate object }$\mathbb{T}$ of $\Crys(k)$ is $(\1,1)$.

The functor%
\begin{equation*}
\Lambda \mapsto (\Lambda ,0)\colon \Crys^{+}(k)\rightarrow \Crys(k)
\end{equation*}%
is faithful, and it is full on torsion-free objects. The forgetful functor
defines an equivalence of the full subcategory of $\Crys(k)$ of torsion
objects with the category of $W(k)$-modules of finite length.

\begin{definition}
\label{pl13}Let $k$ be a field of characteristic $p\neq 0$. If $k$ is
finitely generated over $\mathbb{F}{}_{p}$,
\begin{equation*}
\mathcal{R}{}(k;\mathbb{Q}{}_{p})=\Isoc(k^{\text{pf}}),\quad \mathcal{R}{}(k;%
\mathbb{Z}{}_{p})=\Crys(k^{\text{pf }}).
\end{equation*}%
Otherwise,
\begin{equation*}
\mathcal{R}{}(k;\mathbb{Z}_{p})=\varinjlim_{k^{\prime }\subset k}\mathcal{R}%
{}(k^{\prime };\mathbb{Z}{}_{p}),\quad \mathcal{R}{}(k;\mathbb{Q}{}_{p})=%
\varinjlim_{k^{\prime }\subset k}\mathcal{R}{}(k^{\prime
};\mathbb{Q}_{p}),
\end{equation*}%
where the limits are over the subfields of $k$ finitely generated over $%
\mathbb{F}{}_{p}$. We define $\mathcal{R}{}^{+}(k;\mathbb{Q}{}_{p})$ and $%
\mathcal{R}{}^{+}(k;\mathbb{Z}{}_{p})=\varinjlim_{k^{\prime }\subset k}$ $\Crys%
^{+}(k^{\prime \text{pf}})$ similarly.
\end{definition}

\begin{lemma}
\label{pl13m}For each object $M$ of $\mathcal{R}{}^{+}(k;\mathbb{Z}{}_{p})$,
there is an $n$ such that $M\otimes \mathbb{L}^{\otimes n}$ decomposes into
a direct sum of a torsion object and a torsion-free object. Consequently,
for each $M$ in $\mathcal{R}{}(k;\mathbb{Z}{}_{p})$, $M_{\text{tors}}$ is a
direct summand of $M$.
\end{lemma}

\begin{proof}
Let $M=M_{\text{tors}}\oplus M_{1}$ as $W(k^{\prime \text{pf}})$-modules,
and let $F_{M}$ act as $\left(
\begin{smallmatrix}
a & b \\
0 & c%
\end{smallmatrix}%
\right) $. Then $F_{M\otimes \mathbb{L}^{\otimes n}}$ acts as $\left(
\begin{smallmatrix}
p^{n}a & p^{n}b \\
0 & p^{n}c%
\end{smallmatrix}%
\right) $. Choose $n$ so that $p^{n}M_{\text{tors}}=0$, and then $p^{n}b=0$.
This proves the first statement, and the second follows because $-\otimes
\mathbb{L}$ is an equivalence of categories on $\mathcal{R}{}(k;\mathbb{Z}%
{}_{p})$.
\end{proof}

Let $(\Lambda ,F)$ be an $F$-crystal. For some $N$ there exists a $\sigma
^{-1}$-linear map $V\colon \Lambda \rightarrow \Lambda $ such that $%
FV=p^{N}=VF$. Let $\Lambda ^{\prime }=\Hom_{W(k)}(\Lambda ,W(k))$ and define
$F^{\prime }\colon \Lambda ^{\prime }\rightarrow \Lambda ^{\prime }$ by $%
(F^{\prime }f)(\lambda )=f(V\lambda )$. For an object $(\Lambda ,F,m)$ of $%
\mathcal{R}{}(k;\mathbb{Z}{}_{p})$, define $(\Lambda ,F,m)^{\vee }=(\Lambda
^{\prime },F^{\prime },N-m)$. We say that $(\Lambda ,F,m)$ is \emph{%
reflexive }if the canonical map $(\Lambda ,F,m)\rightarrow (\Lambda
,F,m)^{\vee \vee }$ is an isomorphism.

\begin{proposition}
\label{pl13p} \emph{(a)} An object of $\mathcal{R}{}(k;\mathbb{Z}_{p})$ is
reflexive if and only if it is torsion-free.

\emph{(b)} Every object of $\mathcal{R}{}(k;\mathbb{Z}_{p})$ is a quotient
of a reflexive object.
\end{proposition}

\begin{proof}
(a) is a obvious. After (\ref{pl13m}) it suffices to prove (b) for a torsion
$F$-crystal $N$. Such an $N$ is a quotient $M\rightarrow N$ of a free
finitely generated $W$-module $M$. Endow $M$ with the $F$-crystal structure
such that $F_{M}=\sigma $. After $M\rightarrow N$ has been tensored with a
high power of $\mathbb{L}$, it will be a morphism of $F$-crystals.
\end{proof}

\begin{remark}
\label{pl14}(a) There is a canonical functor $\varinjlim\Isoc(\mathbb{F}%
{}_{q})\rightarrow \Isoc(\mathbb{F}{})$, which is faithful and essentially
surjective (Demazure 1972, p.~ 85) but \textbf{not }full ($W(\mathbb{F}%
{})\neq \varinjlim W(\mathbb{F}{}_{q})$).

(b) Ekedahl 1986, p36, defines a \emph{virtual }$F$\emph{-crystal} to be a
triple $(M,F,\Lambda )$ with $(M,F)$ an isocrystal and $\Lambda $ a $W$%
-submodule of $M$ such that $B(k)\cdot \Lambda =M$; a virtual $F$-crystal is
of \emph{finite type} if $\Lambda $ is finitely generated. The functor $%
((\Lambda ,F_{\Lambda }),m)\mapsto (\Lambda _{\mathbb{Q}{}},F_{\Lambda
}\otimes p^{-m},\Lambda )$ defines an equivalence of the full subcategory of
$\Crys(k)$ of torsion-free objects with the category of virtual $F$-crystals
of finite type.
\end{remark}

\subsection{The categories $\mathcal{R}{}$ in general}

\begin{definition}
\label{pl15}For a field $k$ of characteristic $p\neq 0$,
\begin{eqnarray*}
\mathcal{R}{}(k;\mathbb{\hat{Z}}) &=&\mathcal{R}{}(k;\mathbb{Z}{}^{p})\times
\mathcal{R}{}(k;\mathbb{Z}{}_{p}),\quad \mathcal{R}{}(k;\mathbb{A}_{f})=%
\mathcal{R}{}(k;\mathbb{A}_{f}^{p})\times \mathcal{R}{}(k;\mathbb{Q}{}_{p})
\\
\mathcal{R}{}^{+}(k;\mathbb{\hat{Z}}) &=&\mathcal{R}{}(k;\mathbb{Z}%
{}^{p})\times \mathcal{R}{}^{+}(k;\mathbb{Z}{}_{p}),\quad \mathcal{R}%
{}^{+}(k;\mathbb{A}_{f})=\mathcal{R}{}(k;\mathbb{A}_{f}^{p})\times \mathcal{R%
}{}^{+}(k;\mathbb{Q}{}_{p})
\end{eqnarray*}%
For a field of characteristic zero,%
\begin{equation*}
\mathcal{R}{}(k;\mathbb{\hat{Z}})=\mathcal{R}{}(k;\mathbb{Z}{}^{p}),\quad {}%
\mathcal{R}{}(k;\mathbb{A}_{f})=\mathcal{R}{}(k;\mathbb{A}_{f}^{p})\text{.}
\end{equation*}
\end{definition}

Note that $\mathcal{R}{}(k;\mathbb{A}_{f})\cong \mathcal{R}{}(k;\mathbb{\hat{%
Z}})_{\mathbb{Q}{}}{}$. For notational convenience, we sometimes use $%
\mathcal{R}{}^{+}(k;\mathbb{Z}{}_{l})$ to denote $\mathcal{R}{}(k;\mathbb{Z}%
{}_{l})$ when $l\neq p$.

\subsection{The motivic subcategory of $\mathcal{R}{}(k;\mathbb{Q}_{\ell })$.%
}

Let $\mathcal{M}{}(k;\mathbb{Q}_{\ell })$ be the strictly full subcategory
of $\mathcal{R}(k;\mathbb{Q}{}_{\ell }){}$ whose objects are subquotients of
a Tate twist of a tensor power of an $\ell $-adic \'{e}tale cohomology group
of a smooth projective variety over $k$. When $k=\mathbb{F}{}_{q}$, $%
\mathcal{M}(k;\mathbb{Q}_{\ell })$ consists of the semisimple
representations of $\Gamma =\Gal(k^{\text{al}}/k)$ for which the eigenvalues
of the Frobenius generator of $\Gamma $ are Weil $q$-numbers: the Weil
conjectures (proved by Grothendieck and Deligne) imply that every object of $%
\mathcal{M}{}(k;\mathbb{Q}{}_{\ell })$ has this property, but all such
representations arise already in the cohomology of abelian varieties (cf.
Milne 1994, 3.7). When $k$ is a number field, Conjecture 1 of Fontaine and
Mazur (1995) suggests a description of $\mathcal{M}{}(k;\mathbb{Q}{}_{\ell
}) $.

\section{Ad\`{e}lic Cohomology}

\subsection{Modules over $\mathbb{A}_{f}^{p}$ and lattices}

Recall that $\mathbb{Z}{}^{p}=\prod_{\ell \neq p}\mathbb{Z}{}_{\ell }$ and $%
\mathbb{A}_{f}^{p}$ is the restricted direct product $\mathbb{A}%
_{f}^{p}=\prod_{\ell \neq p}(\mathbb{Q}_{\ell },\mathbb{Z}_{\ell })$.

\begin{plain}
\label{ac1}Endow $(\mathbb{A}_{f}^{p})^{m}$ with its topology as the
restricted direct product of the spaces $\mathbb{Q}{}_{\ell }^{m}$ relative
to the subspaces $\mathbb{Z}{}_{\ell }^{m}$. Every $\mathbb{A}_{f}^{p}$%
-linear bijection $(\mathbb{A}_{f}^{p})^{m}\rightarrow (\mathbb{A}%
_{f}^{p})^{m}$ is bicontinuous. This allows us to endow any free $\mathbb{A}%
_{f}^{p}$-module $M$ of finite rank with a topology by choosing an
isomorphism $\alpha \colon M\rightarrow (\mathbb{A}_{f}^{p})^{m}$.
\end{plain}

\begin{plain}
\label{ac2}Let $M$ be a free $\mathbb{A}_{f}^{p}$-module of finite rank $m$,
so that $M_{\ell }=_{\text{df}}M\otimes _{\mathbb{A}_{f}^{p}}\mathbb{Q}%
{}_{\ell }$ is a $\mathbb{Q}_{\ell }$-vector space of dimension $m$ for all $%
\ell \neq p$. A \emph{lattice} in $M$ is a finitely generated $\mathbb{Z}%
{}^{p}$-submodule $\Lambda $ such that the inclusion $\Lambda
\hookrightarrow M$ induces an isomorphism $\Lambda _{\mathbb{A}%
_{f}^{p}}\rightarrow M$; then $\Lambda _{\ell }=_{\text{df}}\Lambda \otimes
_{\mathbb{Z}{}^{p}}\mathbb{Z}{}_{\ell }$ is a $\mathbb{Z}{}_{\ell }$-lattice
in $M_{\ell }$ for all $\ell \neq p$, and $\Lambda $ is the free $\mathbb{Z}%
^{p}$-submodule spanned by a basis for $M$. Clearly, every lattice in $M$ is
compact and open, and so any two lattices $\Lambda $ and $\Lambda ^{\prime }$
in $M$ are \emph{commensurable}, that is, both quotients $\Lambda /\Lambda
\cap \Lambda ^{\prime }$ and $\Lambda ^{\prime }/\Lambda \cap \Lambda
^{\prime }$ are finite (hence, $\Lambda _{\ell }=\Lambda _{\ell }^{\prime }$
for almost all${}$ $\ell $). Let $M_{0}$ be a $\mathbb{Q}{}$-structure on $M$%
. There is a canonical one-to-one correspondence between the lattices $%
\Lambda $ in $M$ and the lattices $\Lambda _{0}$ in $M_{0}$, namely,
\begin{equation*}
\Lambda \mapsto \Lambda \cap M_{0}\text{, }\Lambda _{0}\mapsto \text{its
closure in }M\text{.}
\end{equation*}
\end{plain}

\begin{plain}
\label{ac3}Let $(M_{\ell })_{\ell \neq p}$ be a family of vector spaces over
the fields $\mathbb{Q}{}_{\ell }$, all of the same finite dimension. Two
families $(\Lambda _{\ell })_{\ell \neq p}$ and $(\Lambda _{\ell }^{\prime
})_{\ell \neq p}$ of lattices in the $M_{\ell }$ will be said to be \emph{%
equivalent} if $\Lambda _{\ell }=\Lambda _{\ell }^{\prime }$ for almost all $%
\ell $. To give a free $\mathbb{A}_{f}^{p}$-module of finite rank is the
same as to give a family $(M_{\ell })_{\ell \neq p}$ together with an
equivalence class of families $(\Lambda _{\ell })_{\ell \neq p}$: given a
family $(M_{\ell },\Lambda _{\ell })_{\ell \neq p}$, let $M=\prod (M_{\ell
},\Lambda _{\ell })$; given $M$, choose an isomorphism $M\rightarrow (%
\mathbb{A}_{f}^{p})^{m}$, and let $\Lambda _{\ell }$ be the inverse image of
$\mathbb{Z}_{\ell }^{m}$. Let $M$ and $M^{\prime }$ be the free $\mathbb{A}%
_{f}^{p}$-modules of finite rank defined by systems $(M_{\ell },\Lambda
_{\ell })_{\ell \neq p}$ and $(M_{\ell }^{\prime },\Lambda _{\ell }^{\prime
})_{\ell \neq p}$ respectively. A family of $\mathbb{Q}{}_{\ell }$-linear
maps $\alpha _{\ell }\colon M_{\ell }\rightarrow M_{\ell }^{\prime }$
defines an $\mathbb{A}_{f}^{p}$-linear map $\alpha \colon M\rightarrow
M^{\prime }$ if and only if $\alpha _{\ell }(\Lambda _{\ell })\subset
\Lambda _{\ell }^{\prime }$ for almost all${}$ $\ell $, in which case $%
\alpha $ is injective if and only if $\alpha _{\ell }$ is injective for all $%
\ell $ and it is surjective if and only if $\alpha _{\ell }$ is surjective
for all $\ell $ and maps $\Lambda _{\ell }$ onto $\Lambda _{\ell }^{\prime }$
for almost all${}$ $\ell $.
\end{plain}

\begin{plain}
\label{ac4}An $\mathbb{A}_{f}^{p}$-linear map $\alpha \colon M\rightarrow
M^{\prime }$ of free $\mathbb{A}_{f}^{p}$-modules of finite rank will be
said to be of \emph{constant rank} if the rank of $\alpha _{\ell }=_{\text{df%
}}\alpha \otimes _{\mathbb{A}_{f}^{p}}\mathbb{Q}{}_{\ell }$ is independent
of $\ell $. Then the kernel, cokernel, image, and co-image of $\alpha $ are
again free $\mathbb{A}_{f}^{p}$-modules of finite rank.
\end{plain}

\subsection{Ad\`{e}lic cohomology (characteristic $p\neq 0$)}

Throughout this subsection, all varieties are smooth and projective over a
separably closed field $k$ of characteristic $p\neq 0$. For such a variety $%
V $, $H^{r}(V,-)$ denotes an \'{e}tale cohomology group and $\Lambda
{}_{\ell }^{i}(V,r)$ denotes the image of $H^{i}(V,\mathbb{Z}_{\ell }(r))$
in $H^{i}(V,\mathbb{Q}{}_{\ell }(r))$. Thus,%
\begin{equation*}
\Lambda {}_{\ell }^{i}(V,r)\cong H^{i}(V,\mathbb{Z}_{\ell }(r))/\{\text{%
torsion}\}.
\end{equation*}

\begin{proposition}
\label{ac5}

\begin{enumerate}
\item The torsion subgroup of $H^{i}(V,\mathbb{Z}{}_{\ell })$ is finite for
all $\ell $, and is zero for almost all${}$ $\ell $.

\item Assume $V$ is connected of dimension $d_{V}$ and fix a prime $\ell
\neq p$. If the cohomology groups $H^{j}(V,\mathbb{Z}{}_{\ell })$ of $V$ are
torsion-free for all $j$, then the pairing
\begin{equation*}
H^{i}(V,\mathbb{Q}{}_{\ell })\times H^{2d_{V}-i}(V,\mathbb{Q}{}_{\ell
}(d_{V}))\rightarrow H^{2d_{V}}(V,\mathbb{Q}{}_{\ell }(d_{V}))\cong \mathbb{Q%
}{}_{\ell }
\end{equation*}%
induces a perfect\footnote{%
A bilinear form $M\times N\rightarrow R$ of $R$-modules is \emph{%
nondegenerate }if its left and right kernels are zero, and \emph{perfect }if
the maps $M\rightarrow \Hom_{R}(N,R)$ and $N\rightarrow \Hom_{R}(M,R)$ it
defines are isomorphisms.} pairing on the cohomology groups with
coefficients in $\mathbb{Z}_{\ell }$.

\item Let $W$ be a second variety over $k$ and fix a prime $\ell \neq p$. If
the groups $H^{i}(V,\mathbb{Z}{}_{\ell })$ and $H^{j}(W,\mathbb{Z}_{\ell })$
are torsion-free for all $i$ and $j$, then the K\"{u}nneth isomorphism
\begin{equation*}
H^{r}(V\times W,\mathbb{Q}{}_{\ell })\cong \bigoplus_{i+j=r}H^{i}(V,\mathbb{Q%
}{}_{\ell })\otimes H^{j}(W,\mathbb{Q}{}_{\ell })
\end{equation*}%
induces an isomorphism on the groups with coefficients in $\mathbb{Z}_{\ell
} $.
\end{enumerate}
\end{proposition}

\begin{proof}
(a) The finiteness of the torsion subgroup of $H^{i}(V,\mathbb{Z}_{\ell })$
is a standard result, and the rest of statement (a) is proved in Gabber 1983.

(b) The pairing on the cohomology groups with coefficients in $\mathbb{Z}%
_{\ell }$ is the unique pairing for which the diagrams
\begin{equation*}
\begin{CD} H^{i}(V,\mathbb{Z}_{\ell
})@.\times@.H^{2d_V-i}(V,\mathbb{Z}_{\ell}(d_V)) @>>>
H^{2d_V}(V,\mathbb{Z}{}_{\ell }(d_V))\cong \mathbb{Z}_{\ell }\\
@VVV@.@VVV@VVV\\ H^{i}(V,\mathbb{Z}/{\ell}^{m} \mathbb{Z}) @.\times@.
H^{2d_V-i}(V,(\mathbb{Z}/{\ell}^{m} \mathbb{Z})(d_V)) @>>>
H^{2d_V}(V,(\mathbb{Z}/{\ell}^{m} \mathbb{Z})(d_V))\cong
\mathbb{Z}{}/{\ell}^m \mathbb{Z} \end{CD}
\end{equation*}%
commute for all $m$. The lower pairing is perfect and the vertical maps are
the cokernels of multiplication by $\ell ^{m}$ on $H^{i}(V,\mathbb{Z}%
{}_{\ell })$ and $H^{2d_{V}-i}(V,\mathbb{Z}{}_{\ell }(d_{V}))$ because we
are assuming there is no torsion. This implies (even with $m=1$) that the
discriminant of the top pairing is an $\ell $-adic unit, and so the top
pairing is also perfect.\footnote{%
We thank M. Nori for this proof, which is simpler than our original.}

(c) This is a standard result.
\end{proof}

\begin{definition}
\label{ac5m}For $V$ as above, $H^{i}(V,\mathbb{A}_{f}^{p}(r))$ is the
restricted direct product of the spaces $H^{i}(V,\mathbb{Q}{}_{\ell }(r))$
with respect to the lattices $\Lambda {}^{i}(V,r)$. Equivalently,
\begin{equation*}
H^{i}(V,\mathbb{A}_{f}^{p}(r))=H^{i}(V,\mathbb{Z}{}^{p}(r))\otimes \mathbb{Q}%
{}
\end{equation*}%
where $H^{i}(V,\mathbb{Z}{}^{p}(r))=\varprojlim_{(n,p)=1}H^{i}(V,\mathbb{\mu }%
_{n}^{\otimes r})$.
\end{definition}

\begin{remark}
\label{ac6} (a) Because the Betti numbers of a smooth projective variety are
independent of $\ell $, $H^{i}(V,\mathbb{A}_{f}^{p})$ is a free $\mathbb{A}%
_{f}^{p}$-module of finite rank. A regular map $\alpha \colon V\rightarrow W$
defines an $\mathbb{A}_{f}^{p}$-linear map $\alpha ^{\ast }\colon H^{i}(W,%
\mathbb{A}_{f}^{p})\rightarrow H^{i}(V,\mathbb{A}_{f}^{p})$.

(b) For a connected $V$, the proposition implies the following:
%TCIMACRO{\TeXButton{BeginQuote}{\begin{quotation}}}%
%BeginExpansion
\begin{quotation}%
%EndExpansion
\noindent the pairings in (b) of the proposition define a perfect pairing
\begin{equation*}
\quad \quad H^{i}(V,\mathbb{A}_{f}^{p})\times H^{2d_{V}-i}(V,\mathbb{A}%
_{f}^{p}(d_{V}))\rightarrow H^{2d_{V}}(V,\mathbb{A}_{f}^{p}(d_{V}))\cong
\mathbb{A}_{f}^{p}
\end{equation*}%
of free $\mathbb{A}_{f}^{p}$-modules of finite rank.%
%TCIMACRO{\TeXButton{EndQuote}{\end{quotation}}}%
%BeginExpansion
\end{quotation}%
%EndExpansion
In other words, Poincar\'{e} duality holds for $\mathbb{A}_{f}^{p}$%
-cohomology. This implies that the map $\alpha ^{\ast }$ in (a) has an $%
\mathbb{A}_{f}^{p}$-linear adjoint
\begin{equation*}
\alpha _{\ast }\colon H^{j}(V,\mathbb{A}_{f}^{p}(d_{V}))\rightarrow
H^{j+2d_{W}-2d_{V}}(W,\mathbb{A}_{f}^{p}(d_{W})),\quad j=2d_{V}-i.
\end{equation*}

(c) Part (c) of the theorem implies the following: there is a canonical
isomorphism of free $\mathbb{A}_{f}^{p}$-modules of finite rank
\begin{equation*}
H^{r}(V\times W,\mathbb{A}_{f}^{p})\cong \bigoplus_{i+j=r}H^{i}(V,\mathbb{A}%
_{f}^{p})\otimes H^{j}(W,\mathbb{A}_{f}^{p}).
\end{equation*}

(d) Assume (for simplicity) that $V$ is connected. Let $f$ be an algebraic
correspondence of degree $r$ from $V$ to $W$ with rational coefficients,
i.e., an element of $Z^{d_{V}+r}(V\times W)_{\mathbb{Q}{}}$. For all but
finitely many $\ell $, $f$ will be $\ell $-integral, and so its $\ell $-adic
cohomology classes define an element $cl(f)\in H^{2d_{V}+2r}(V\times W,%
\mathbb{A}_{f}^{p}(d_{V}+r))$. Combining this remark with (b) and (c), we
obtain a map (also denoted $f$)%
\begin{equation*}
x\mapsto q_{W\ast }(cl(f)\cup q_{V}^{\ast }x)\colon H^{i}(V,\mathbb{A}%
_{f}^{p}(m))\rightarrow H^{i+2r}(W,\mathbb{A}_{f}^{p}(m+r))
\end{equation*}%
for all $m$. Here $q_{W}$ and $q_{V}$ are the projection morphisms.
\end{remark}

\begin{remark}
\label{ac7}The main theorem of de Jong 1996 shows that Gabber's theorem
(Gabber 1983) holds also for the cohomology with compact support of
arbitrary varieties over separably closed fields. Thus, appropriately
interpreted, so do the remaining results of this subsection.
\end{remark}

\subsection{Ad\`{e}lic cohomology (characteristic zero)}

Let $V$ be a variety (not necessarily smooth or projective) over an
algebraically closed field $k$ of characteristic zero. For any model $%
V_{0}/k_{0}$ of $V$ over an algebraically closed subfield $k_{0}$ of $k$ and
embedding $k_{0}\hookrightarrow \mathbb{C}{}$,%
\begin{equation*}
H^{i}(V,\mathbb{Z}{}_{l})\cong H^{i}(V_{0},\mathbb{Z}{}_{l})\cong
H^{i}(V_{0}(\mathbb{C}{}),\mathbb{Z}{})_{\mathbb{Z}{}_{l}}.
\end{equation*}%
The group $H^{i}(V_{0}(\mathbb{C}{}),\mathbb{Z}{})$ is finitely generated,
and so $H^{i}(V,\mathbb{Z}{}_{l})$ is torsion-free for almost all${}$ $l$.
The results in the last subsection hold \emph{mutatis mutandis} for the
group $H^{i}(V,\mathbb{A}_{f}(r))$ defined to be the restricted topological
product of the $H^{i}(V,\mathbb{Q}{}_{l})$ with respect to the $\Lambda
_{l}^{i}(V,r)$. For any choice of a pair $V_{0}/k_{0}$, $k_{0}%
\hookrightarrow \mathbb{C}{}$ as above,%
\begin{equation*}
H^{i}(V,\mathbb{A}_{f})\cong H^{i}(V_{0},\mathbb{A}_{f})\cong H^{i}(V_{0}(%
\mathbb{C}{}),\mathbb{Q}{})_{\mathbb{A}_{f}}.
\end{equation*}

\section{An Ad\`{e}lic Tate Conjecture}

Throughout this section, $V$ is a smooth projective variety of dimension $d$
over a field $k$. For $\ell \neq \mathrm{char}(k)$, $H^{i}(V,\mathbb{Z}%
_{\ell }(r))$ denotes the \'{e}tale cohomology group. We let $\bar{k}$ $=k^{%
\text{sep}}$ and $\bar{V}=V_{/\bar{k}}$. Finally, $\Gamma =\Gal(\bar{k}/k)$.

\subsection{Conjectures on algebraic cycles in $\ell $-adic cohomology}

The group $Z_{\ell \text{-hom}}^{r}(V){}$ is the co-image of the cycle class
map $Z^{r}(V)\rightarrow H^{2r}(\bar{V},\mathbb{Z}_{\ell }(r))$, and $Z_{%
\text{num}}^{r}(V)$ is the quotient of $Z^{r}(V)$ by numerical equivalence.

Recall (Tate 1994) that, when $k$ is finitely generated over the prime
field, there are the following conjectures.

\begin{description}
\item[$\mathbf{T}^{r}\mathbf{(V/k,\ell )}$] The cycle map $Z_{\ell \text{-hom%
}}^{r}(V)_{\mathbb{Q}{}_{\ell }}\rightarrow H^{2r}(\bar{V},\mathbb{Q}%
{}_{\ell }(r))^{\Gamma }$ is surjective.

\item[$\mathbf{E}^{r}\mathbf{(V/k,\ell )}$] The quotient map $Z_{\ell \text{%
-hom}}^{r}(V)_{\mathbb{Q}}\rightarrow Z_{\text{num}}^{r}(V)_{\mathbb{Q}{}}$
is an isomorphism.

\item[$\mathbf{I}^{r}\mathbf{(V/k,\ell )}$] The cycle map $Z_{\ell \text{-hom%
}}^{r}(V)_{\mathbb{Q}{}_{\ell }}\rightarrow H^{2r}(\bar{V},\mathbb{Q}%
{}_{\ell }(r))$ is injective.

\item[$\mathbf{S}^{r}\mathbf{(V/k,\ell )}$] The map $H^{2r}(\bar{V},\mathbb{Q%
}_{\ell }(r))^{\Gamma }\rightarrow H^{2r}(\bar{V},\mathbb{Q}_{\ell
}(r))_{\Gamma }$ induced by the identity map is bijective.
\end{description}

\subsection{An ad\`{e}lic Tate conjecture}

\begin{lemma}
\label{ac8}Assume $k$ is finitely generated over the prime field. If $%
T^{r}(V/k,\ell )$, $T^{d-r}(V/k,\ell )$, and $S^{r}(V/k,\ell )$ hold for
almost all${}$ $\ell $, then the cycle class map%
\begin{equation*}
Z_{\ell \text{-hom}}^{r}(V)_{\mathbb{Z}{}_{\ell }}\rightarrow H^{2r}(\bar{V},%
\mathbb{Z}_{\ell }(r))^{\Gamma }
\end{equation*}%
is surjective for almost all${}$ $\ell $.
\end{lemma}

\begin{proof}
According to Tate 1994, 2.9 (c)$\Rightarrow $(d), the conditions imply that $%
E^{r}(V/k,\ell )$ and $E^{d-r}(V/k,\ell )$ hold for almost all${}$ $\ell $.

Note that $E^{r}(V/k,\ell )$ implies $Z_{\ell \text{-hom}}^{r}(V)/\{$torsion$%
\}\cong Z_{\text{num}}^{r}(V)$, which is finitely generated (see, for
example, Milne 1980, VI 11.7). Therefore, for those $\ell $ for which $%
E^{r}(V/k,\ell )$ and $E^{d-r}(V/k,\ell )$ hold,
\begin{equation*}
Z_{\ell \text{-hom}}^{r}(V)/\{\text{torsion}\}\times Z_{\ell \text{-hom}%
}^{d-r}(V)/\{\text{torsion}\}\rightarrow \mathbb{Z}{}
\end{equation*}%
is a nondegenerate pairing of finitely generated abelian groups whose
discriminant $D$ is independent of $\ell $.

Consider the diagram%
\begin{equation*}
\begin{CD} Z_{{\ell}\text{-hom}}^{r}(V)_{\mathbb{Z}_{\ell}} @.\times @.
Z_{{\ell}\text{-hom}}^{d-r}(V)_{\mathbb{Z}_{\ell}} @>\bullet>>
\mathbb{Z}_{\ell}\\ @VVV @.@VVV \Vert \\
H^{2r}(\bar{V},\mathbb{Z}_{\ell}(r))^{\Gamma }@.\times @. H^{2d-2r}(\bar{V},
\mathbb{Z}_{\ell}(d-r))^{\Gamma } @>\cup>> \mathbb{Z}_{\ell}.\\ \end{CD}
\end{equation*}%
For those $\ell $ for which $T^{r}(V/k,\ell )$ and $T^{d-r}(V/k,\ell )$
hold, the vertical arrows have finite cokernel.

For those $\ell $ for which

\begin{itemize}
\item the groups $H^{2r}(\bar{V},\mathbb{Z}{}_{\ell }(r))$ and $H^{2d-2r}(%
\bar{V},\mathbb{Z}{}_{\ell }(d-r))$ are torsion-free,

\item $E^{r}(V/k,\ell )$ and $E^{d-r}(V/k,\ell )$ hold,

\item $T^{r}(V/k,\ell )$ and $T^{d-r}(V/k,\ell )$ hold, and

\item $\ell $ does not divide $D$,
\end{itemize}

\noindent the vertical maps are surjective. Since these conditions hold for
almost all${}$ primes $\ell $, this proves the lemma.
\end{proof}

\begin{theorem}
\label{ac9}Let $k$ be as in the lemma. If $T^{r}(V/k,\ell )$, $%
T^{d-r}(V/k,\ell )$, and $S^{r}(V/k,\ell )$ hold for all $\ell \neq p$, then
the cycle class map defines a bijection%
\begin{equation*}
Z_{\text{num}}^{r}(V)_{\mathbb{A}_{f}^{p}}\rightarrow H^{2r}(\bar{V},\mathbb{%
A}_{f}^{p}(r))^{\Gamma }.
\end{equation*}
\end{theorem}

\begin{proof}
We have to show that the cycle class map $Z^{r}(V)\rightarrow H^{2r}(\bar{V},%
\mathbb{Z}_{\ell }(r))$ defines

\begin{enumerate}
\item a bijection $Z_{\text{num}}^{r}(V)_{\mathbb{Q}{}_{\ell }}\rightarrow
H^{2r}(\bar{V},\mathbb{Q}{}_{\ell }(r))^{\Gamma }$ for all $\ell \neq p$, and

\item a surjection $Z_{\text{num}}^{r}(V)_{\mathbb{Z}{}_{\ell }}\rightarrow
\Lambda {}_{\ell }^{2r}(\bar{V},r)^{\Gamma }$ for almost all${}$ $\ell $.
\end{enumerate}

\noindent Applying Tate 1994, 2.9, we see that $E^{r}(V/k,\ell )$, $%
I^{r}(V/k,\ell )$, and $T^{r}(V/k,\ell )$ hold for all $\ell \neq p$. This
implies (a), and (b) follows from the lemma.
\end{proof}

\begin{corollary}
\label{ac9a}Let $V$ be a (smooth projective) variety over $\mathbb{F}_{q}$,
and let $r$ be an integer. The following statements are equivalent:

\begin{enumerate}
\item for some prime $\ell \neq p$, $\dim _{\mathbb{Q}{}}Z_{\text{num}%
}^{r}(V)=\dim _{\mathbb{Q}{}_{\ell }}H^{2r}(V,\mathbb{Q}{}_{\ell
}(r))^{\Gamma }$;

\item the order of the pole of the zeta function $Z(V,t)$ of $V$ at $%
t=q^{-r} $ is equal to the rank of $Z_{\text{num}}^{r}(V)$;

\item the cycle class map defines a bijection $Z_{\text{num}}^{r}(V)_{%
\mathbb{A}_{f}^{p}}\rightarrow H^{2r}(\bar{V},\mathbb{A}_{f}^{p}(r))^{\Gamma
}$.
\end{enumerate}
\end{corollary}

\begin{proof}
Apply Tate 1994, 2.9, and the theorem.
\end{proof}

\subsection{Counterexamples to an integral Tate conjecture\label{ac10}}

Assume $k$ is finitely generated over the prime field, and consider the
cycle class map%
\begin{equation}
Z^{r}(V)_{\mathbb{Z}{}_{\ell }}\rightarrow H^{2r}(\bar{V},\mathbb{Z}{}_{\ell
}(r))^{\Gamma }\text{.}  \label{e9}
\end{equation}%
Over infinite fields $k$ there exist $V$'s without a rational zero cycle of
degree $1$ (for example, any curve of genus zero without a rational point or
any curve\footnote{%
For such a curve, even $CH^{1}(\bar{V})^{\Gamma }\xr{\text{degree}}\mathbb{Z}%
{}$ fails to be surjective (Lang and Tate 1958).} of genus $1$ with index $%
\neq 1$) and for such a variety (\ref{e9}) will fail to be surjective for $%
r=d$ and suitable $\ell $. The argument that Atiyah and Hirzebruch (1962)
used to prove that not all torsion cohomology classes on a complex variety
are algebraic (contradicting Hodge's original conjecture) also applies over
finite fields\footnote{%
For this, one needs to use that the odd-dimensional Steenrod operations with
values in $H^{\ast }(\bar{V},\mathbb{Z}{}/\ell \mathbb{Z}{})$ vanish on
algebraic cycles (Brosnan 1999).} and proves that (\ref{e9}) is not
surjective on torsion classes.

Griffiths and Harris (1985) conjecture that if $V$ is a very general
hypersurface\footnote{%
A statement is said to hold for a very general hypersurface of degree $d$ if
it holds for all hypersurfaces in the complement of countably many proper
subvarieties of the space of hypersurfaces of degree $d$.} of degree $d\geq
6 $ in $\mathbb{P}{}^{4}$, then $d$ divides the degree of every curve on $V$%
. Koll\'{a}r has obtained some results in the direction of this conjecture
(see Ballico et al. 1992, p.~ 135). Thus, it is known that, if $k$ is an
uncountable algebraically closed field and $\ell >3$ is distinct from the
characteristic of $k$, then there exist smooth hypersurfaces $V\subset
\mathbb{P}{}^{4}$ such that the cycle class map
\begin{equation*}
Z^{2}(V)_{\mathbb{Z}{}_{\ell }}\rightarrow H^{4}(V,\mathbb{Z}{}_{\ell
}(2))\cong \mathbb{Z}{}_{\ell }
\end{equation*}%
is not surjective (cf. Schoen 1998, 0.6).

We do not know of an example with $k$ finite for which (\ref{e9}) is not
surjective modulo torsion (however, Schoen 1998, 2.1, gives an example with $%
V$ singular). For $r=1$ and $k$ finite, one sees from the Kummer sequence
that $T^{1}(V/k,\ell )$ implies that the map (\ref{e9}) is surjective.

\section{The Ad\`{e}lic Fibre Functor}

\subsection{Isomotives}

Fix a field $k$ of characteristic exponent $p$, and let $\mathcal{S}{}$ be a
class of smooth projective varieties over $k$ satisfying the following
condition:

\begin{plain}
\label{af1}$\mathcal{S}{}$ is closed under passage to a connected component
and under the formation of products and disjoint unions, and it contains the
zero-dimensional varieties and the projective spaces.
\end{plain}

\begin{plain}
\label{af1a}Fix an adequate equivalence relation $\sim $. Recall that, for
smooth projective varieties $V$ and $W$ over $k$, the group of
correspondences of degree\emph{\ }$i$ from $V$ to $W$ is
\begin{equation*}
\Corr^{i}(V,W)=\oplus _{j}Z_{\sim }^{\dim (V_{j})+i}(V_{j}\times W)_{\mathbb{%
Q}{}}
\end{equation*}%
where the $V_{j}$ are the equidimensional components of $V$.
\end{plain}

\begin{plain}
\label{af1b}The category $\mathcal{CV}_{\sim }^{0}(k){}$ of correspondences
based on $\mathcal{S}{}$, has one object $hV$ for each $V\in \mathcal{%
\mathcal{S}{}}$, and%
\begin{equation*}
\Hom_{\mathcal{CV}_{\sim }^{0}(k){}}(hV,hW)=\Corr^{0}(V,W)_{\mathbb{Q}{}}%
\text{.}
\end{equation*}%
With the definition
\begin{equation}
hV\otimes hW=h(V\times W)  \label{e11}
\end{equation}%
and the obvious constraints, $\mathcal{CV}_{\sim }^{0}(k){}$ becomes a
tensor category over $\mathbb{Q}{}$ (Saavedra 1972, VI 4.1.1).
\end{plain}

\begin{plain}
\label{af1c}The category of effective isomotives $\mathcal{M}{}_{\sim
}^{+}(k;\mathbb{Q}{};\mathcal{S}{})$ is the pseudo-abelian hull of $\mathcal{%
CV}_{\sim }^{0}(k)$. Its objects are pairs $(hV,e)$ with $e$ an idempotent
in $\End(hV)$. It is a pseudo-abelian tensor category over $\mathbb{Q}$.
\end{plain}

\begin{plain}
\label{af1d}In $\mathcal{M}{}_{\sim }^{+}(k;\mathbb{Q}{};\mathcal{S}{})$
there is a canonical decomposition%
\begin{equation*}
h\mathbb{P}{}^{1}=\1\oplus \mathbb{L}\text{,}
\end{equation*}%
and the category of isomotives $\mathcal{M}_{\sim }{}(k;\mathbb{Q}{};%
\mathcal{S}{})$ is obtained from $\mathcal{M}{}_{\sim }^{+}(k;\mathbb{Q}{};%
\mathcal{S}{})$ by inverting $\mathbb{L}$ (Saavedra 1972, VI 4.1.3.1). Its
objects are pairs $((hV,e),m)$ with $(hV,e)$ an effective motive and $m\in
\mathbb{Z}{}$. It is a pseudo-abelian rigid tensor category over $\mathbb{Q}%
{}$ (Saavedra 1972, VI 4.1.3.5). We sometimes write $h(V,e,m)$ for $%
((hV,e),m)$.
\end{plain}

\begin{plain}
\label{af1e}We say that \emph{the K\"{u}nneth projectors are algebraic} for
a smooth projective variety $V$ if this is so for the $l$-adic \'{e}tale
cohomologies for all $l\neq p$ and the crystalline cohomology
(characteristic $p\neq 0$). There are then well-defined elements $\pi
^{i}\in Z_{\text{num}}^{i}(V)_{\mathbb{Q}{}}$ (independent of the cohomology
theory) that give these projectors for each cohomology theory (when $k$ is
finite, this follows from the Riemann hypothesis, and when $k$ is arbitrary
it can be proved by a specialization argument).

When the K\"{u}nneth projectors are algebraic for all $V$ in $\mathcal{S}{}$%
, we use them to modify the commutativity constraints in $\mathcal{M}%
{}_{\sim }^{+}(k;\mathbb{Q}{};\mathcal{S}{})$ and $\mathcal{M}{}_{\sim }(k;%
\mathbb{Q}{};\mathcal{S}{}){}$ (Saavedra 1972, VI 4.2.1.5).
\end{plain}

For the remainder of this section, we assume that \ $\sim $ is not coarser
than $\ell $-adic homological equivalence on $Z^{r}(V)_{\mathbb{Q}{}}$ for
any $\ell \neq p$, so that the formula
\begin{equation*}
\omega _{\ell }(h(V,e,m))=e\left( H^{\ast }(\bar{V},\mathbb{Q}{}_{\ell
}(m))\right) ,\quad \bar{V}=V_{k^{\text{sep}}},
\end{equation*}%
defines a tensor functor $\omega _{\ell }\colon \mathcal{M}{}_{\sim }(k;%
\mathbb{Q}{};\mathcal{S}{})\rightarrow \Vct(\mathbb{Q}{}_{\ell })$, and we
drop $\sim $ from the notation.

\subsection{The ad\`{e}lic fibre functor}

With the notation of (\ref{ac5m}), define
\begin{equation*}
\omega ^{p}(hV)=H^{\ast }(\bar{V},\mathbb{A}_{f}^{p}).
\end{equation*}%
Because of (\ref{ac6}d) and our assumption on $\sim $, a correspondence $f$
of degree $0$ from $V$ to $W$ defines a homomorphism $\omega ^{p}(f)\colon
\omega ^{p}(hV)\rightarrow \omega ^{p}(hW)$, and one sees immediately from
the definitions that $\omega ^{p}(f\circ g)=\omega ^{p}(f)\circ \omega
^{p}(g)$. The isomorphisms
\begin{equation*}
\omega ^{p}(hV\otimes hW)\overset{(\ref{e11})}{=}\omega ^{p}(h(V\times W))%
\overset{\ref{ac6}c}{\cong }\omega ^{p}(hV)\otimes \omega ^{p}(hW)
\end{equation*}%
provide $\omega ^{p}$ with the structure of a tensor functor $\mathcal{CV}%
^{0}(k){}\rightarrow \Mod(\mathbb{A}_{f}^{p})$. It is $\mathbb{Q}{}$-linear
(in particular, additive).

\begin{theorem}
\label{af2}The functor $\omega ^{p}\colon \mathcal{CV}^{0}(k){}\rightarrow %
\Mod(\mathbb{A}_{f}^{p})$ extends (essentially uniquely) to a $\mathbb{Q}{}$%
-linear tensor functor $\omega ^{p}\colon \mathcal{M}{}(k;\mathbb{Q}{};%
\mathcal{S}{})\rightarrow \Mod(\mathbb{A}_{f}^{p})$. For each $\ell \neq p$,
$\omega ^{p}\otimes _{\mathbb{A}_{f}^{p}}\mathbb{Q}{}_{\ell }\cong \omega
_{\ell }$.
\end{theorem}

\begin{proof}
As $\Mod(\mathbb{A}_{f}^{p})$ is abelian, $\omega ^{p}$ extends (essentially
uniquely) to $\mathcal{M}{}^{+}(k;\mathbb{Q}{};\mathcal{S}{})$ because of
the universality of the pseudo-abelian hull. Because $\omega ^{p}(\mathbb{L}%
)=\mathbb{A}_{f}^{p}(-1)$ is invertible in $\Mod(\mathbb{A}_{f}^{p})$, $%
\omega ^{p}$ has an essentially unique extension from $\mathcal{M}{}^{+}(k;%
\mathbb{Q}{};\mathcal{S}{})$ to $\mathcal{M}{}(k;\mathbb{Q}{};\mathcal{S}{})$%
. For $V\in \mathcal{S}{}$, $\omega ^{p}(hV)\otimes _{\mathbb{A}_{f}^{p}}%
\mathbb{Q}{}_{\ell }\cong \omega _{\ell }(hV)$, and this property is
preserved under each extension.
\end{proof}

\begin{remark}
\label{af3}There is the following formula,%
\begin{equation*}
\omega ^{p}(h(V,e,m))=e(H^{\ast }(\bar{V},\mathbb{A}_{f}^{p}(m)))
\end{equation*}%
where, on the right, $e$ denotes the map defined in (\ref{ac6}d). For a
motive $X=h(V,e,m)$, define $\lambda {}_{\ell }(X)$ to be the image of
\begin{equation*}
H^{\ast }(\bar{V},\mathbb{Z}_{\ell }(m))\rightarrow H^{\ast }(\bar{V},%
\mathbb{\mathbb{Q}{}}_{\ell }(m))\overset{e}{\rightarrow }H^{\ast }(\bar{V},%
\mathbb{Q}{}_{\ell }(m))\text{.}
\end{equation*}%
Then, $\omega ^{p}(X)$ is the restricted direct product
\begin{equation*}
\omega ^{p}(X)=\prod_{\ell \neq p}(\omega _{\ell }(X),\lambda _{\ell }(X)).
\end{equation*}
\end{remark}

\begin{proposition}
\label{af4}If $\sim $ is numerical equivalence and the K\"{u}nneth
projectors of all varieties in $\mathcal{S}{}$ are algebraic, then for all
isomotives $X$,

\begin{enumerate}
\item $\omega ^{p}(X)$ is a free $\mathbb{A}_{f}^{p}$-module of finite rank;

\item for all morphisms $\alpha $ of isomotives, $\omega ^{p}(\alpha )$ has
constant rank.
\end{enumerate}
\end{proposition}

\begin{proof}
The hypotheses imply that $\mathcal{M}{}(k;\mathbb{Q};\mathcal{S}{}{})$ is a
Tannakian category over $\mathbb{Q}{}$ (Jannsen 1992). For any fibre functor
$\omega $ on a Tannakian category with values in the vector spaces over a
field, and any morphism $\alpha \colon X\rightarrow Y$ in the category
\begin{eqnarray*}
\dim (\omega (X)) &=&\rank(X) \\
\rank(\omega (\alpha )) &=&\rank X-\rank\Ker(\alpha ).
\end{eqnarray*}
In particular, $\dim (\omega (X))$ and $\rank(\omega (\alpha ))$ are
independent of $\omega $. When we apply this statement to the functors $%
\omega _{\ell }$, we obtain the theorem.
\end{proof}

\begin{remark}
\label{af5}Because we are assuming $\sim $ is not coarser than $\ell $-adic
homological equivalence, the hypotheses of the proposition are that, for all
$\ell \neq p$, numerical equivalence coincides with $\ell $-adic homological
equivalence on algebraic cycles with $\mathbb{Q}{}$-coefficients and that
the K\"{u}nneth projectors are algebraic. Without these hypotheses, the
conclusion is not known (see the discussion Katz 1994, p.~28, (3)).
\end{remark}

\subsection{The functor $\protect\omega ^{p}\colon \mathcal{M}{}(k;\mathbb{Q}%
{};\mathcal{S}{})\rightarrow \mathcal{R}{}(k;\mathbb{A}_{f}^{p})$.}

\begin{theorem}
\label{af6}The functor $\omega _{l}$ has a canonical factorization through
the forgetful functor $\mathcal{R}{}(k;\mathbb{Q}{}_{\ell })\rightarrow \Vct(%
\mathbb{Q}{}_{\ell })$. Under the hypotheses of \ref{af4}, the functor $%
\omega ^{p}$ has a canonical factorization through the forgetful functor $%
\mathcal{R}{}(k;\mathbb{A}_{f}^{p})\rightarrow \Mod(\mathbb{A}_{f}^{p})$.
\end{theorem}

\begin{proof}
Let $V\in \mathcal{S}{}$. When $k$ is finitely generated over the prime
field, there is a canonical continuous action of $\Gal(k^{\text{sep}}/k)$ on
$H^{i}(\bar{V},\mathbb{Q}{}_{\ell }(m))$, which therefore lies in $\mathcal{R%
}{}(k;\mathbb{Q}{}_{\ell })$. For any correspondence $f\in \Corr^{r}(V,W)$,
the map%
\begin{equation*}
f\colon H^{i}(\bar{V},\mathbb{Q}_{\ell }(m))\rightarrow H^{i+r}(\bar{W},%
\mathbb{Q}{}_{\ell }(m+r))
\end{equation*}%
respects the actions of $\Gal(k^{\text{sep}}/k)$. Therefore, the functor
\begin{equation*}
\omega _{\ell }\colon \mathcal{M}{}(k;\mathbb{Q}{})\rightarrow \Vct(\mathbb{Q%
}_{\ell })
\end{equation*}%
has a canonical factorization through the forgetful functor $\mathcal{R}{}(k;%
\mathbb{Q}{}_{\ell })\rightarrow \Vct(\mathbb{Q}{}_{\ell })$. The argument
for $\omega ^{p}$ is similar.

For an arbitrary field $k$, any varieties $V,W$ and correspondence $f$ will
have models $V_{1},W_{1},f_{1}$ over some subfield $k_{1}$ of $k$ finitely
generated over the prime field. The image of $\omega _{l}(f_{1})\colon
\omega _{l}(V_{1})\rightarrow \omega _{l}(W_{l})$ under the functor $%
\mathcal{R}{}(k_{1};\mathbb{Q}_{l})\rightarrow \mathcal{R}{}(k;\mathbb{Q}%
{}_{l})$ is independent of the choice of $k_{1},V_{1},W_{1},f_{1}$. Thus,
the statements hold also for $k$.
\end{proof}

\section{The Category of Motives\label{cm}}

We fix a field $k$ of characteristic exponent $p$, and a class $\mathcal{S}%
{} $ of smooth projective varieties over $k$ satisfying the condition (\ref%
{af1}) and such that

\begin{enumerate}
\item for all $V$ in $\mathcal{S}{}$, the K\"{u}nneth projectors are
algebraic, and

\item for all $V$ in $\mathcal{S}{}$ and all $l$ (including $p$), numerical
equivalence coincides with $l$-adic homological equivalence for algebraic
cycles with $\mathbb{Q}{}$-coefficients.
\end{enumerate}

\noindent Let $\mathcal{M}{}^{+}(k;\mathbb{Q}{})$ and $\mathcal{M}{}(k;%
\mathbb{Q}{})$ be the categories $\mathcal{M}{}_{\text{num}}^{+}(k;\mathbb{Q}%
{};\mathcal{S}{})$ and $\mathcal{M}{}_{\text{num}}(k;\mathbb{Q};\mathbb{%
\mathcal{S}{})}$ defined in \ref{af1e}. Recall that they are semisimple
abelian $\mathbb{Q}{}$-categories (Jannsen 1992) and that $\mathcal{M}{}(k;%
\mathbb{Q}{})$ is Tannakian (Deligne 1990).

Because $\omega _{l}$ is a fibre functor on $\mathcal{M}{}(k;\mathbb{Q}{})$,
the map%
\begin{equation}
\Hom_{\mathcal{M}{}(k;\mathbb{Q}{})}(X,Y)_{\mathbb{Q}{}_{l}}\rightarrow \Hom%
_{\mathcal{R}{}(k;\mathbb{Q}{}_{l})}(\omega _{l}(X),\omega _{l}(Y))
\label{e27}
\end{equation}%
defined by $\omega _{l}$ is injective (Deligne 1990, 2.13). It follows that
it is also injective when $\mathcal{M}{}$ and $\mathcal{R}{}$ are replaced
by $\mathcal{M}{}^{+}$ and $\mathcal{R}{}^{+}$. We say that the $l$\emph{%
-adic Tate conjecture holds} for $\mathcal{M}{}(k;\mathbb{Q}{})$ or $%
\mathcal{M}{}^{+}(k;\mathbb{Q})$ if the corresponding map is an isomorphism
for $l$.

While conditions (a) and (b) are conjectured to hold for all smooth
projective varieties, they have been proved for very few varieties (see
Kleiman 1994 for a list). In the final subsection of this section (Variants,
p\pageref{variants}), we list some alternative categories $\mathcal{M}%
{}^{+}(k;\mathbb{Q}{})$ and $\mathcal{M}{}(k;\mathbb{Q}{})$ to which our
constructions apply, but which require no unproven conjectures.

Our construction of $\mathcal{M}{}(k;\mathbb{Z}{})$ is suggested by the
following simple observation: for a finitely generated $\mathbb{Z}{}$-module
$M$, $M\cong M_{\mathbb{\hat{Z}}{}}\times _{M_{\mathbb{A}_{f}}}M_{\mathbb{Q}%
} $; therefore, any category whose Hom-sets are finitely generated $\mathbb{Z%
}{}$-modules admits a fully faithful functor into the fibre product of a $%
\mathbb{\hat{Z}}{}$-category with a $\mathbb{Q}{}$-category over a $\mathbb{A%
}{}_{f}$-category.

\subsection{The category of effective motives${}$ $\mathcal{M}^{+}(k;\mathbb{%
\mathbb{Z}{}})$}

\begin{definition}
\label{cm7}The \emph{category of effective motives }$\mathcal{M}^{+}{}(k;%
\mathbb{Z}{})$ over $k$ is the full subcategory of the fibre product
category
\begin{equation*}
\mathcal{R}^{+}(k;\mathbb{\hat{Z}})\times _{\mathcal{R}{}^{+}(k;\mathbb{A}%
_{f})}\mathcal{M}{}^{+}(k;\mathbb{Q}{})
\end{equation*}%
whose objects $(X_{f},X_{0},x_{f})$ are those for which the torsion subgroup
of $X_{f}$ is finite.
\end{definition}

Thus, an effective motive\emph{\ }$X$ over $k$ is a triple $%
(X_{f},X_{0},x_{f})$ consisting of

\begin{enumerate}
\item an object $X_{f}=(X_{l})_{l}$ of $\mathcal{R}{}^{+}(k;\mathbb{\hat{Z})}
$ such that $X_{l}$ is torsion-free for almost all${}$ $l$,

\item an effective isomotive $X_{0}$, and

\item an isomorphism $x_{f}\colon (X_{f})_{\mathbb{Q}{}}\rightarrow \omega
_{f}(X_{0})$ in $\mathcal{R}^{+}(k;\mathbb{A}_{f})$.
\end{enumerate}

\noindent A morphism of effective motives\emph{\ }$\alpha \colon
X\rightarrow Y$ is a pair of morphisms
\begin{equation*}
(\alpha _{f}\colon X_{f}\rightarrow Y_{f},\;\alpha _{0}\colon
X_{0}\rightarrow Y_{0})
\end{equation*}%
such that%
\begin{equation*}
\begin{CD} X_{f} @>{\alpha_f }>> Y_{f} \\ @VV{x_{f}}V@VV{y_{f}}V \\ \omega
_{f}(X_{0}) @>{\omega _{f}(\alpha_0)}>>\omega_f(Y_{0})\end{CD}
\end{equation*}%
commutes. We say that a motive $X=(X_{f},X_{0},x_{f})$ is \emph{finite }if $%
X_{0}=0$.

Note that to give $(X_{f},x_{f})$ amounts to giving an $X_{l}$ in $\mathcal{R%
}{}(k;\mathbb{Z}{}_{l})$ for each $l$ together with homomorphisms $%
x_{l}\colon X_{l}\rightarrow \omega _{l}(X_{0})$ that induce isomorphisms $%
X_{l}\rightarrow \lambda _{l}(X_{0})$ for almost all${}$ $l$ and
isomorphisms $(X_{l})_{\mathbb{Q}{}}\rightarrow \omega _{l}(X_{0})$ for all $%
l$. We usually regard $x_{f}$ as a map $X_{f}\rightarrow \omega _{f}(X_{0})$
with finite kernel. We use the notations: $X_{f}=X^{p}\times X_{p}$, $%
X^{p}=\prod_{l\neq p}X_{l}$, $x_{f}=x^{p}\times x_{p}$, $\alpha _{f}=\alpha
^{p}\times \alpha _{p}$.

\begin{proposition}
\label{cm8}The category $\mathcal{M}^{+}(k;\mathbb{Z})$ is abelian. A
sequence
\begin{equation*}
X\rightarrow Y\rightarrow Z
\end{equation*}%
in $\mathcal{M}^{+}(k;\mathbb{\mathbb{Z}})$ is exact if and only if the
sequences%
\begin{equation*}
X_{0}\rightarrow Y_{0}\rightarrow Z_{0}\text{ and }X_{f}\rightarrow
Y_{f}\rightarrow Z_{f}
\end{equation*}%
are exact in $\mathcal{M}^{+}(k;\mathbb{Q})$ and $\mathcal{R}{}^{+}(k;%
\mathbb{\hat{Z}})$ respectively.
\end{proposition}

\begin{proof}
Both statements are true for the fibre product category, and so it remains
to show that $\mathcal{M}{}^{+}(k;\mathbb{Z}{})$ is closed under the
formation of kernels and cokernels in the larger category. Any subobject of
an object in $\mathcal{M}{}^{+}(k;\mathbb{Z}{})$ is also in $\mathcal{M}%
{}^{+}(k;\mathbb{Z}{})$. Therefore, $\mathcal{M}{}^{+}(k;\mathbb{Z}{})$ is
closed under the formation of kernels, and, in proving that the cokernel of $%
\alpha $ is in $\mathcal{M}{}^{+}(k;\mathbb{Z}{})$, we may assume that $%
\alpha $ is injective. Thus, given a short exact sequence
\begin{equation*}
0\rightarrow X\overset{\alpha }{\rightarrow }Y\rightarrow Z\rightarrow 0
\end{equation*}%
in the fibre product category with $X,Y$ in $\mathcal{M}{}^{+}(k;\mathbb{Z}%
{})$, we have to show that $Z$ is in $\mathcal{M}{}^{+}(k;\mathbb{Z}{})$. As
$\mathcal{M}{}^{+}(k;\mathbb{Q}{})$ is semisimple, the sequence of $\mathcal{%
M}{}^{+}(k;\mathbb{Q}{})$-components splits. But $\omega ^{p}$ is additive,
and so%
\begin{equation*}
0\rightarrow \omega ^{p}(X_{0})\xr{\omega^{p}(\alpha_0)}\omega
^{p}(Y_{0})\rightarrow \omega ^{p}(Z_{0})\rightarrow 0
\end{equation*}%
is split-exact, which implies that%
\begin{equation*}
0\rightarrow \lambda {}{}_{\ell }(X_{0})\rightarrow \lambda {}{}_{\ell
}(Y_{0})\rightarrow \lambda {}{}_{\ell }(Z_{0})\rightarrow 0
\end{equation*}%
is split-exact for almost all${}$ $\ell $. In turn, this implies that
\begin{equation*}
0\rightarrow X_{\ell }\overset{\alpha _{\ell }}{\rightarrow }Y_{\ell
}\rightarrow Z_{\ell }\rightarrow 0
\end{equation*}%
is split-exact for almost all${}$ $\ell $. Thus, the torsion subgroup of $%
Z_{\ell }$ is zero for almost all${}$ $\ell $ and finite for all $\ell $,
which shows that $Z$ is in $\mathcal{M}{}^{+}(k;\mathbb{Z}{})$.
\end{proof}

\begin{corollary}
\label{cm8m}Let $X$ and $X^{\prime }$ be submotives of an effective motive $%
Y $, and let $X\cap X^{\prime }=X\times _{Y}X^{\prime }$. If $X_{0}=Y_{0}$,
then the cokernel $Z$ of $X\cap X^{\prime }\rightarrow X^{\prime }$ is
finite.
\end{corollary}

\begin{proof}
Because $X\mapsto X_{0}$ is exact, $(X\times _{Y}X^{\prime }\rightarrow
X^{\prime })_{0}\cong (X_{0}\times _{Y_{0}}X_{0}^{\prime }\rightarrow
X_{0}^{\prime })$, which is an isomorphism if $X_{0}=Y_{0}$.
\end{proof}

\begin{proposition}
\label{cm7m}The category $\mathcal{M}{}^{+}(k;\mathbb{Z}{})$ is noetherian
(i.e., its objects are noetherian).
\end{proposition}

\begin{proof}
Let%
\begin{equation*}
X(1)\subset X(2)\subset X(3)\subset \cdots
\end{equation*}%
be a sequence of a submotives of an effective motive $Y$. Eventually, the
sequence $X(1)_{0}\subset X(2)_{0}\subset \cdots $ becomes stationary,
equal, say, to $X(\infty )_{0}\subset Y_{0}$. Let $X(\infty
)_{f}=Y_{f}\times _{\omega _{f}(Y_{0})}\omega _{f}(X(\infty )_{0})$. Then $%
X(\infty )=(X(\infty )_{f},X(\infty )_{0},\cdot )$ is an effective motive,
and $X(1)\subset X(2)\subset \cdots $ is a sequence of submotives of $%
X(\infty )$ such that $X(r)_{0}=X(\infty )_{0}$ for $r\geq r_{1}$, some $%
r_{1}$. Now $X(\infty )/X(r_{1})$ is finite (\ref{cm8m}), and so the
sequence $X(1)\subset X(2)\subset \cdots $ becomes stationary.
\end{proof}

\begin{lemma}
\label{cm9}Given a morphism $\beta \colon X_{0}\rightarrow Y_{0}$ of
effective isomotives, there exists an integer $m$ such that $m\beta $ is in
the image of $\Hom(X,Y)\rightarrow \Hom(X_{0},Y_{0})$.
\end{lemma}

\begin{proof}
We have to find an integer $m$ and morphisms $\alpha _{l}$ such that the
diagrams
\begin{equation*}
\begin{CD} X_{l} @>{\alpha_{l}}>>Y_{l} \\ @VV{x_{l }}V @VV{y_{l }}V \\
\omega _{l}(X_{0})@>{\omega _{l}(m\beta )}>> \omega _{l}(Y_{0}) \end{CD}
\end{equation*}%
commute. Because $\omega _{l}(\beta )$ is the $l$-component of a
homomorphism $\omega ^{p}(X_{0})\rightarrow \omega ^{p}(Y_{0})$, it maps $%
x_{l}(X_{l})$ into $y_{l}(Y_{l})$ for almost all $l$, and so there exists an
integer $m$ such that $m\omega _{l}(\beta )$ maps $x_{l}(X_{l})$ into $%
y_{l}(Y_{l})$ for all $l$. Choose an integer $n$ that kills the torsion in $%
Y_{l}$ for all $l$. Then the map
\begin{equation*}
\begin{CD} X_{l }@>x_{l }>>x_{l}X_{l} @>\omega _{l}(mn\beta )>>
y_{l}Y_{l}=Y_{l}/{\text{{torsion}}} .\end{CD}
\end{equation*}%
lifts to a map $X_{l}\rightarrow Y_{l}$.
\end{proof}

\begin{proposition}
\label{cm10}For all $X,Y$ in $\mathcal{M}^{+}(k;\mathbb{\mathbb{Z}})$, $\Hom%
(X,Y)$ is a finitely generated $\mathbb{Z}$-module (modulo $p$-torsion when $%
k$ is infinite and $p\neq 1$).
\end{proposition}

\noindent \textsc{Proof. }As $\Hom(X,Y)\subset \End(X\oplus Y)$, it suffices
to show that all endomorphism rings are finitely generated $\mathbb{Z}{}$%
-modules (modulo \ldots ). Because $\omega ^{p}$ is faithful, $X\mapsto
X_{f} $ is faithful. In particular, the map%
\begin{equation*}
\alpha \mapsto \alpha _{f}\colon \End(X)\rightarrow \End(X_{f})
\end{equation*}%
is injective, and so the torsion subgroup of $\End(X)$ is contained in that
of $\End(X_{f})$, which is finite except possibly for the $p$-torsion when $%
k $ is infinite and $p\neq 1$. The kernel of the map
\begin{equation*}
\alpha \mapsto \alpha _{0}\colon \End(X)\rightarrow \End(X_{0})
\end{equation*}%
is obviously torsion, and so it remains to show that its image $R$ is
finitely generated over $\mathbb{Z}$. If $\beta \in \End(X_{0})$ is in $R$,
then $\omega _{l}(\beta )$ stabilizes the lattice $x_{l}(X_{l})$ for all $l$%
. Thus, the characteristic polynomial of $\beta $ has coefficients in $%
\mathbb{Z}\subset \mathbb{Q}{}$. Let $\alpha _{1},\ldots ,\alpha _{n}$ be a
basis for $\End(X_{0})$ as a $\mathbb{Q}$-vector space. According to Lemma %
\ref{cm9}, after possibly replacing each $\alpha _{i}$ with $m\alpha _{i}$
for some integer $m$, we may suppose that each $\alpha _{i}\in R$. Because $%
\mathcal{M}{}^{+}(k;\mathbb{Q}{})$ is semisimple, $\End(X_{0})$ is a
semisimple $\mathbb{Q}{}$-algebra, and so the trace pairing $\alpha ,\beta
\mapsto \Tr_{\End(X_{0})/\mathbb{Q}{}}(\alpha \beta )$ is nondegenerate.
Therefore, there is a $\mathbb{Q}$-basis $\beta _{1},\ldots ,\beta _{n}$ for
$\End(X_{0})$ dual to $\alpha _{1},\ldots ,\alpha _{n}$. Now the usual
argument shows that
\begin{equation*}
R\subset \mathbb{Z}\beta _{1}+\cdots +\mathbb{Z}\beta _{n}.
\end{equation*}%
Namely, let $\beta \in R$. Then $\beta $ can be written uniquely as a linear
combination $\beta =\sum b_{j}\beta _{j}$ of the $\beta _{j}$ with
coefficients $b_{j}\in \mathbb{Q}{}$, and we have to show that each $%
b_{j}\in \mathbb{Z}$. As $\alpha _{i}$ and $\beta $ are in $R$, so also is $%
\beta \cdot \alpha _{i}$, and so $\Tr(\beta \cdot \alpha _{i})\in \mathbb{Z}$%
. But
\begin{equation*}
\Tr(\beta \cdot \alpha _{i})=\Tr(\sum_{j}b_{j}\beta _{j}\alpha
_{i})=\sum_{j}b_{j}\Tr(\beta _{j}\alpha _{i})=b_{i}.\quad \quad \square
\end{equation*}

\begin{remark}
\label{cm10m}(a) If $Y$ has no $p$-torsion (i.e., $Y_{p}$ is torsion-free),
then $\Hom(X,Y)$ has no $p$-torsion and hence is finitely generated over $%
\mathbb{Z}{}$.

(b) Let $p\neq 1$. The effective motive $X$ with $X_{0}=0$ and $%
X_{f}=X_{p}=k $ ($F$ acting as zero) has $\End(X)=k$. The same phenomenon
occurs with subgroup schemes of abelian varieties: a supersingular elliptic
curve has $\alpha _{p}=\Spec k[T]/(T^{p})$ as a subgroup, and $\End(\alpha
_{p})=k$.
\end{remark}

\begin{proposition}
\label{cm11}The functor $X\mapsto X_{0}\colon \mathcal{M}^{+}(\mathbb{F};%
\mathbb{Z}{})\rightarrow \mathcal{M}^{+}(\mathbb{F}{};\mathbb{Q}{})$ defines
an equivalence of categories
\begin{equation*}
\mathcal{M}^{+}(\mathbb{F};\mathbb{Z})_{\mathbb{Q}{}}\rightarrow \mathcal{M}%
^{+}(\mathbb{F};\mathbb{Q}).
\end{equation*}
\end{proposition}

\begin{proof}
Because the kernel of
\begin{equation*}
\alpha \mapsto \alpha _{0}\colon \Hom(X,Y)\rightarrow \Hom(X_{0},Y_{0})
\end{equation*}%
is torsion, the map
\begin{equation*}
r\otimes \alpha \mapsto r\alpha _{0}\colon \mathbb{Q}{}\otimes \Hom%
(X,Y)\rightarrow \Hom(X_{0},Y_{0})
\end{equation*}%
is injective, and Lemma \ref{cm9} shows that it is surjective. Thus, the
functor is faithful and full. Given an $X$ in $\mathcal{M}{}^{+}(k;\mathbb{Q}%
{})$, clearly there exists an $X^{p}$ in $\mathcal{R}{}^{+}(k;\mathbb{Z}%
^{p}) $ such that $(X^{p})_{\mathbb{Q}{}}\approx \omega ^{p}(X)$. Using
this, one sees that $X\mapsto X_{0}$ is essentially surjective.
\end{proof}

\begin{remark}
\label{cm12}An object $X$ of $\mathcal{M}^{+}(k;\mathbb{\mathbb{Z}})$ is
said to be \emph{torsion} if $mX=0$ for some $m$ (equivalently, $X$ is
finite). The torsion objects form a thick subcategory of $\mathcal{M}%
{}^{+}(k;\mathbb{Z)}$ and Proposition \ref{cm11} shows that the functor
\begin{equation*}
X\mapsto X_{0}\colon \mathcal{M}^{+}(k;\mathbb{\mathbb{Z}})\rightarrow
\mathcal{M}^{+}(k;\mathbb{Q})
\end{equation*}%
realizes $\mathcal{M}^{+}(k;\mathbb{Q})$ as the quotient of $\mathcal{M}%
^{+}(k;\mathbb{\mathbb{Z}})$ by its subcategory of torsion objects.
\end{remark}

\begin{proposition}
\label{cm13}The functor
\begin{equation*}
\mathcal{M}^{+}(k;\mathbb{Z}{})_{\mathbb{Z}{}_{l}}\rightarrow \mathcal{R}%
{}^{+}(k;\mathbb{Z}_{l})
\end{equation*}%
defined by $X\mapsto X_{l}$ is faithful; it is full if and only if the $l$%
-adic Tate conjecture holds for $\mathcal{M}{}^{+}(k;\mathbb{Q}_{l})$.
\end{proposition}

\begin{proof}
We first show that%
\begin{equation*}
z\otimes \alpha \mapsto z\alpha _{l}\colon \mathbb{\mathbb{Z}}_{l}\otimes _{%
\mathbb{Z}{}}\Hom_{\mathcal{M}^{+}(k;\mathbb{Z})}(X,Y)\rightarrow \Hom_{%
\mathcal{R}^{+}(k;\mathbb{\mathbb{Z}}_{l})}(X_{l},Y_{l})
\end{equation*}%
is injective. The torsion subgroup of $\mathbb{\mathbb{Z}}_{l}\otimes \Hom_{%
\mathcal{M}^{+}(k;\mathbb{Z})}(X,Y)$ is equal to the $l$-primary component
of the torsion subgroup of $\Hom_{\mathcal{M}^{+}(k;\mathbb{Z})}(X,Y)$, and
this maps injectively to $\Hom_{\mathcal{R}^{+}(k;\mathbb{\mathbb{Z}}%
_{l})}(X_{l},Y_{l})$ (because the functor $X\mapsto X_{f}$ is faithful).
Thus, the kernel of $z\otimes \alpha \mapsto z\alpha _{l}$ is torsion-free.
Consider the diagram%
\begin{equation*}
\begin{CD} \mathbb{Z}_{l }\otimes _{\mathbb{Z}}\Hom(X,Y)@>>>
\Hom_{\mathcal{R}^+(k;\mathbb{Z}_{l })}(X_{l },Y_{l }) \\ @VVV @VVV\\
\mathbb{Q}_{l }\otimes _{\mathbb{Q}}\Hom(X_{0},Y_{0}) @>>>
\Hom_{\mathcal{R}^+(k;\mathbb{Q}_{l })}(\omega _{l}(X_{0}),\omega _{l
}(Y_{0})). \end{CD}
\end{equation*}%
The left vertical arrow is obtained from
\begin{equation*}
\Hom(X,Y)\rightarrow \Hom(X_{0},Y_{0})
\end{equation*}%
by tensoring with $\mathbb{Z}{}_{l}$, and hence has torsion kernel (\ref%
{cm11}). The bottom arrow is injective (\ref{e27}). Hence, the kernel of the
top horizontal arrow is torsion, but we have already shown it to be
torsion-free; it is therefore zero.

We next show that the cokernel of $z\otimes \alpha \mapsto z\alpha _{l}$ is
torsion-free. Suppose that $\beta \in \Hom_{\mathcal{R}^{+}(k;\mathbb{%
\mathbb{Z}}_{l})}(X_{l},Y_{l})$ is such that $l^{m}\beta =z\alpha _{l}$ for
some $m\in \mathbb{N}{}$, $\alpha \in \Hom(X,Y)$, and $z\in \mathbb{Z}{}_{l}$%
. Because $z\otimes \alpha \mapsto z\alpha _{l}$ is (obviously) surjective
on torsion, we may suppose that $z$ is not divisible by $l$. Then $\alpha
_{l}$ is divisible by $l^{m}$ in $\Hom(X_{l},Y_{l})$. For $l^{\prime }\neq l$%
, $l^{m}$ is a unit in $\mathbb{Z}{}_{l^{\prime }}$, and so $\alpha
_{l^{\prime }}$ is also divisible by $l^{m}.$ Now $\alpha ^{\prime }$, with
\begin{equation*}
\alpha _{0}^{\prime }=l^{-m}\alpha _{0}\text{,}\quad \alpha _{l}^{\prime
}=z^{-1}\beta ,\quad \alpha _{l^{\prime }}^{\prime }=l^{-m}\alpha
_{l^{\prime }}\text{ for }l^{\prime }\neq l
\end{equation*}%
is a morphism $X\rightarrow Y$ such that $z\alpha _{l}^{\prime }=\beta $.

Finally, we assume the $l$-adic Tate conjecture for $\mathcal{M}{}^{+}(k;%
\mathbb{Q}{})$, and we prove that the cokernel of $z\otimes \alpha \mapsto
z\alpha _{l}$ is torsion. For this, it suffices to prove that
\begin{equation*}
\mathbb{Q}{}_{l}\otimes _{\mathbb{Z}{}}\Hom(X,Y)\rightarrow \mathbb{Q}%
{}\otimes _{\mathbb{Z}{}}\Hom_{\mathcal{R}^{+}(k;\mathbb{Z}%
_{l})}(X_{l},Y_{l})
\end{equation*}%
is surjective. But this follows from the isomorphisms
\begin{equation*}
\mathbb{Q}\otimes _{\mathbb{Z}{}}\Hom(X,Y){}\overset{\text{(\ref{cm11})}}{%
\cong }\Hom(X_{0},Y_{0}),
\end{equation*}%
and
\begin{equation*}
\mathbb{Q}_{l}\otimes _{\mathbb{Q}{}}\Hom(X_{0},Y_{0})\overset{\text{Tate}}{%
\rightarrow }\Hom_{\mathcal{R}^{+}(k;\mathbb{Q}_{l})}(\omega
_{l}(X_{0}),\omega _{l}(Y_{0}))\cong \mathbb{Q}\otimes _{\mathbb{Z}}\Hom_{%
\mathcal{R}^{+}(k;\mathbb{\mathbb{Z}}_{l})}(X_{l},Y_{l})\text{.}
\end{equation*}
\end{proof}

\begin{proposition}
\label{cm13m}The category%
\begin{equation*}
\mathcal{M}{}^{+}(k;\mathbb{Z}{})=\varinjlim_{k^{\prime }}\mathcal{M}%
{}^{+}(k^{\prime };\mathbb{Z}{})
\end{equation*}%
where $k^{\prime }$ runs over the subfields of $k$ finitely generated over
the prime field.
\end{proposition}

\begin{proof}
The statement is true for each of the categories $\mathcal{M}{}^{+}(k;%
\mathbb{Q}{})$, $\mathcal{R}{}^{+}(k;\mathbb{\hat{Z}})$, and $\mathcal{R}%
{}^{+}(k;\mathbb{A}_{f})$, and it follows easily for their fibre product and
its subcategory $\mathcal{M}{}^{+}(k;\mathbb{Z}{})$.
\end{proof}

\subsection{The category of motives}

\begin{definition}
The \emph{category of motives }$\mathcal{M}{}(k;\mathbb{Z}{})$ is the full
subcategory of
\begin{equation*}
\mathcal{R}(k;\mathbb{\hat{Z}})\times _{\mathcal{R}{}(k;\mathbb{A}_{f})}%
\mathcal{M}{}(k;\mathbb{Q}{})
\end{equation*}%
whose objects $(X_{f},X_{0},x_{f})$ are those for which the torsion subgroup
of $X_{f}$ is finite.
\end{definition}

The results in the previous subsection hold \emph{mutandis mutatis }for $%
\mathcal{M}{}(k;\mathbb{Z)}$ --- in particular, $\mathcal{M}{}(k;\mathbb{Z}%
{})$ is a noetherian abelian category whose quotient by its torsion
subcategory is $\mathcal{M}{}(k;\mathbb{Q})$.

\begin{proposition}
The evident functor $\mathcal{M}{}^{+}(k;\mathbb{Z}{})\rightarrow \mathcal{M}%
{}(k;\mathbb{Z}{})$ is faithful, and it is full on torsion-free objects (on
all objects when $k$ has characteristic zero). It realizes $\mathcal{M}{}(k;%
\mathbb{Z}{})$ as the category obtained from $\mathcal{M}{}^{+}(k;\mathbb{Z}%
{})$ by inverting $\mathbb{L}$.
\end{proposition}

\begin{proof}
Straightforward.
\end{proof}

\begin{remark}
\label{cm14}(a)\emph{\ }When $k=\mathbb{Q}{}$ there is a canonical functor
from $\mathcal{M}{}(k;\mathbb{Z}{})$ to the category of systems of
realizations over $\mathbb{Q}{}$ with integer coefficients defined in
Deligne 1989, 1.23.

(b)\emph{\ }Let $k$ be a subfield of $\mathbb{C}{}$, and let $\mathcal{M}%
{}(k;\mathbb{Q}{})$ be the category based on a class $\mathcal{S}{}$
satisfying (\ref{af1}) but using the \textquotedblleft
motivated\textquotedblright\ classes as correspondences (Andr\'{e} 1996,
4.2). In this case, our definition of $\mathcal{M}{}(k;\mathbb{Z}{})$ is
equivalent to that of Andr\'{e} (ibid. 8.1).
\end{remark}

\subsection{Artin motives.}

Fix a separable algebraic closure $\bar{k}$ of $k$, and let $\Gamma =\Gal(%
\bar{k}/k)$. The category of Artin isomotives is the full subcategory $%
\mathcal{\mathcal{M}{}}{}^{\text{Artin}}(k;\mathbb{Q}{})$ of $\mathcal{M}%
{}^{+}(k;\mathbb{Q}{})$ of objects of the form $X=hV$ with $V$ a variety of
dimension $0$. For such an $X$, define $\pi _{0}(X)=\pi _{0}(\bar{V})$ ---
it is a finite set with a continuous action of $\Gamma $. The functor $%
H\colon X\mapsto \mathbb{Q}{}^{\pi _{0}(X)}=_{\text{df}}\Map(\pi _{0}(X),%
\mathbb{Q}{})$ defines an equivalence of $\mathcal{M}^{\text{Artin}}(k;%
\mathbb{Q}{})$ with the category $\mathcal{R}(\Gamma ;\mathbb{Q}{})$ of
continuous representations of $\Gamma $ on finite-dimensional $\mathbb{Q}{}$%
-vector spaces (cf. Deligne and Milne 1982, 6.17). Hence, $\mathcal{M}{}^{%
\text{Artin}}(k;\mathbb{Q}{})$ is a Tannakian category.

We define the category of \emph{Artin motives }to be the full subcategory $%
\mathcal{M}{}^{\text{Artin}}(k;\mathbb{Z)}$ of $\mathcal{M}{}^{+}(k;\mathbb{%
Z)}$ of objects $X$ such that $X_{0}$ is an Artin isomotive and (when $p\neq
1$) $F_{X_{p}}$ is bijective.

Recall (Saavedra 1972, VI 3.1.2; Demazure 1972, p.~ 69) that there is a
fully faithful functor
\begin{equation}
\gamma \colon \mathcal{R}{}(\Gamma ;\mathbb{Z}{}_{p})\rightarrow \Crys%
^{+}(k^{\text{pf}})=\mathcal{R}{}^{+}(k;\mathbb{Z}{}_{p})  \label{e35}
\end{equation}%
whose essential image consists of the crystals with bijective $F$. Let $X$
be an Artin motive. When $p\neq 1$, we choose $X_{p}^{\prime }\in \mathcal{R}%
{}(\Gamma ;\mathbb{Z}{}_{p})$ so that $\gamma (X_{p}^{\prime })=X_{p}$, and
we let $\omega _{p}^{\prime }(X_{0})=\mathbb{(}X_{p}^{\prime })_{\mathbb{Q}%
{}}$; otherwise we set $X_{p}^{\prime }=0=\omega _{p}^{\prime }(X_{0})$.
Define $H(X)$ so that the right hand square in the following diagram is
cartesian:

\begin{equation*}
\begin{CD} 0 @>>> H(X)_{\text{tors}} @>>> H(X) @>>> H(X_0)\\ @. @VV{\cong}V
@VVV @VVV \\ 0 @>>> X_{\text{tors}}^{p}\times X_{p\text{tors}}^{\prime }
@>>> X^{p}\times X_{p}^{\prime } @>>> \omega ^{p}(X_{0})\times \omega
_{p}^{\prime }(X_{0}) \end{CD}
\end{equation*}%
The diagram is exact and commutative, and $H(X)$ is a finitely generated $%
\mathbb{Z}{}$-module with a continuous action of $\Gamma $ (discrete
topology on $H(X)$).

\begin{proposition}
\label{cm15}The functor $H\colon \mathcal{M}{}^{\text{Artin}}(k;\mathbb{Z}%
{})\rightarrow \mathcal{R}{}(\Gamma ;\mathbb{Z}{})$ is an equivalence of
categories.
\end{proposition}

\begin{proof}
The functor $H\colon \mathcal{M}{}^{\text{Artin}}(k;\mathbb{Q}{})\rightarrow
\mathcal{R}{}(\Gamma ;\mathbb{Q}{})$ defines an equivalence of categories%
\begin{equation*}
\mathcal{R}{}(\Gamma ;\mathbb{Z}{})\times _{\mathcal{R}{}(\Gamma ;\mathbb{Q}%
{})}\mathcal{M}{}^{\text{Artin}}(k;\mathbb{Q}{})\rightarrow \mathcal{R}%
{}(\Gamma ;\mathbb{Z})\times _{\mathcal{R}{}(\Gamma ;\mathbb{Q}{})}\mathcal{R%
}{}(\Gamma ;\mathbb{Q}{})=\mathcal{R}{}(\Gamma ;\mathbb{Z}{})
\end{equation*}%
and $\gamma $ defines a full faithful functor%
\begin{equation*}
\mathcal{R}{}(\Gamma ;\mathbb{Z}{})\times _{\mathcal{R}{}(\Gamma ;\mathbb{Q}%
{})}\mathcal{M}{}^{\text{Artin}}(k;\mathbb{Q}{})\rightarrow \mathcal{R}(k;%
\mathbb{\hat{Z}})\times _{\mathcal{R}{}(k;\mathbb{A}_{f})}{}\mathcal{M}{}^{%
\text{Artin}}(k;\mathbb{Q}{})
\end{equation*}%
whose essential image is $\mathcal{M}{}^{\text{Artin}}(k;\mathbb{Z})$.
\end{proof}

\subsection{Rationally decomposed one-motives}

Recall that a one-motive is a triple $(G,X,u)$ with $G$ a semi-abelian
variety, $X$ a finitely generated free $\mathbb{Z}{}$-module on which $%
\Gamma $ acts continuously, and $u$ a $\Gamma $-homomorphism $X\rightarrow
G(k^{\text{sep}})$. A one-motive $X\xr{u}G$ is said to be \emph{decomposed}
if $u=0$ and the extension
\begin{equation*}
0\rightarrow T\rightarrow G\rightarrow A\rightarrow 0\quad (T\text{ a torus,
}A\text{ an abelian variety)}
\end{equation*}
splits, and it is \emph{rationally decomposed} if it is isogenous to a
decomposed one-motive. Over a finite field, every one-motive is rationally
decomposed. The one-motives form an additive category $\mathcal{M}{}_{1}(k)$
with the decomposed and rationally decomposed one-motives as additive
subcategories $\mathcal{M}{}_{1}^{\text{d}}(k)\subset $ $\mathcal{M}{}_{1}^{%
\text{rd}}(k)\subset \mathcal{M}{}_{1}(k)$.

Choose a quasi-inverse $F\colon \mathcal{R}{}(\Gamma ;\mathbb{Z}%
{})\rightarrow \mathcal{M}{}^{\text{Artin}}(k;\mathbb{Z}{})$ to the functor $%
H$ of (\ref{cm15})\thinspace .

\begin{proposition}
\label{cm16}For $L=(X\xr{0}T\times A)$ in $\mathcal{M}{}_{1}^{\text{d}}(k)$,
define $hL$ to be the effective motive $F(X)\oplus h_{\mathbb{Z}%
{}}^{1}A\oplus F(X_{\ast }(T)(1))$. Then $h$ extends to a functor $h\colon
\mathcal{M}{}_{1}^{\text{rd}}(k)\rightarrow \mathcal{M}{}^{+}(k;\mathbb{Z}%
{}) $. (Here $h_{\mathbb{Z}{}}^{1}A$ denotes $h^{1}A$ endowed with the $%
\mathbb{Z}{}$-structure provided by the maps $H^{1}(A,\mathbb{Z}%
{}_{l})\rightarrow H^{1}(A,\mathbb{Q}{}_{l})$.)
\end{proposition}

\noindent \textsc{Proof. }Let $L$ be a rationally decomposed one-motive.
There exists a decomposed one-motive $L^{\prime }$ and an isogeny $\alpha
\colon L\rightarrow L^{\prime }$ such that $(L^{\prime },\alpha )$ is
universal, and $\alpha $ induces isomorphisms $\omega ^{p}(L)\cong \omega
^{p}(L^{\prime })$ and $\omega _{p}(L)\cong \omega _{p}(L^{\prime })$. We
define $hL$ to be the isomotive $(hL^{\prime })_{0}$ equipped with the
structures%
\begin{eqnarray*}
\lambda {}{}^{p}L &\rightarrow &\omega ^{p}(L)\cong \omega ^{p}(L^{\prime
})\cong \omega ^{p}((hL^{\prime })_{0}) \\
\lambda {}{}_{p}L &\rightarrow &\omega _{p}(L)\cong \omega _{p}(L^{\prime
})\cong \omega _{p}((hL^{\prime })_{0}).\quad \quad \square
\end{eqnarray*}

\begin{proposition}
\label{cm17}The functor $h\colon \mathcal{M}{}_{1}^{\text{rd}}(k)\rightarrow
\mathcal{M}{}^{+}(k;\mathbb{Z}{})$ is fully faithful (except possibly when $%
k $ has characteristic $2$ and is not algebraic over $\mathbb{F}_{2}$).
\end{proposition}

\begin{proof}
It suffices to prove this in the case that $k$ is finitely generated over
the prime field, in which case it is a consequence of the following theorem
of Tate (1966a), Zarhin (1975), Faltings (Faltings and W\"{u}stholz 1984),
and de Jong (1998): for abelian varieties $A$ and $B$ over $k$, the map $\Hom%
_{k}(A,B)_{\mathbb{Z}{}_{l}}\rightarrow \Hom_{\mathcal{R}{}(k;\mathbb{Z}%
{}_{l})}(T_{l}A,T_{l}B)$ is bijective (except possibly when $k$ has
characteristic $2$ \ldots ).
\end{proof}

\begin{remark}
Propositions \ref{cm16} and \ref{cm17} show that, for any two rationally
decomposed one-motives $L$, $M$, there is an injective homomorphism%
\begin{equation*}
\Ext_{\mathcal{M}{}_{1}(k)}^{1}(L,M)\rightarrow \Ext_{\mathcal{M}{}^{+}(k;%
\mathbb{Z}{})}^{1}(hL,hM)\text{.}
\end{equation*}%
This map will not in general be surjective because an extension of two $W[F]$%
-modules for which there exist maps $V$ with $FV=p=VF$ need not be similarly
endowed. In fact, the formula Milne 1968, Theorem 3, for the order of the $%
\Ext^{1}$ of abelian varieties disagrees with that in Theorem \ref{gc.a}
below.
\end{remark}

\subsection{The category $\mathcal{M}{}(\mathbb{F}_{q};\mathbb{Z}{})$}

For a motive $X=h(V,e,m)$ over $\mathbb{F}_{q}$, the Frobenius endomorphism
of $V/\mathbb{F}{}_{q}$ defines a canonical Frobenius element $\pi _{X}\in %
\End(X)$. For a motive $X$ over $\mathbb{F}{}$, each model $X_{1}/\mathbb{F}%
{}_{p^{n}}$ of $X$ defines a pair $(\pi ,n)$ with $\pi $ equal to $\pi
_{X_{1}}$ regarded as an endomorphism of $X$. Any two pairs $(\pi ,n)$, $%
(\pi ^{\prime },n^{\prime })$ arising in this way are equivalent in the
sense that $\pi ^{n^{\prime }m}=\pi ^{\prime nm}$ for some integer $m>0$.
The \emph{germ of a Frobenius endomorphism of }$X$ is the equivalence class
of pairs containing the pair arising from one (hence, every) model of $X$
over a finite field.

\begin{proposition}
\label{ct12}Assume the Tate conjecture holds for $\mathcal{M}{}(\mathbb{F}{};%
\mathbb{Q}{})$. Then
\begin{equation*}
X\mapsto X_{/\mathbb{F}{}}\colon \mathcal{M}(\mathbb{F}_{p^{n}};\mathbb{Z}%
)\rightarrow \mathcal{M}(\mathbb{F};\mathbb{\mathbb{Z}{}})
\end{equation*}%
defines an equivalence of $\mathcal{M}(\mathbb{F}_{p^{n}};\mathbb{Z})$ with
the category of pairs $(X,\pi )$ consisting of a motive $X$ in $\mathcal{M}(%
\mathbb{F};\mathbb{Z})$ and an endomorphism $\pi $ of $X$ such that $(\pi
,n) $ represents $\pi _{X}$. A similar statement holds for effective motives.
\end{proposition}

\begin{proof}
Easy consequence of the similar statement for isomotives (Milne 1994, 3.5).
\end{proof}

\subsection{Remarks on the definition of $\mathcal{M}{}(k;\mathbb{Z}{})$}

\begin{plain}
\label{ct14}As Andr\'{e} notes (1996, 8.1), in view of the counterexamples
of Atiyah and Hirzebruch (1962) to the original Hodge conjecture, it would
be \textquotedblleft peu judicieux\textquotedblright\ to define a category
of integral motives by using algebraic cycles with integer coefficients in
Grothendieck's construction, in other words, with%
\begin{equation*}
\Hom(\1,h_{\mathbb{Z}{}}V)=Z_{\text{num}}^{\dim V}(V\times V)\text{.}
\end{equation*}%
There is also the difficulty of knowing what to replace the pseudo-abelian
hull with.

With our definitions, a smooth projective variety $V$ defines a complex of
motives, rather than a motive (see an article in preparation), but we can
make the ad hoc definitions: $h_{\mathbb{Z}{}}^{r}V$ is the isomotive $%
h^{r}V $ with the $\mathbb{Z}{}$-structure provided the maps $H^{r}(V,%
\mathbb{Z}{}_{l})\rightarrow H^{r}(V,\mathbb{Q}_{l})$, and $h_{\mathbb{Z}%
{}}V=\oplus h_{\mathbb{Z}{}}^{r}V$. If $V$ has torsion-free cohomology, then%
\begin{equation*}
\Hom(\1,hV)=Z_{\text{num}}^{\dim V}(V\times V)\cap \End(H^{\ast }(V,\mathbb{Z%
}{}_{l}))\text{.}
\end{equation*}%
There is a map from $Z_{\text{num}}^{\dim V}(V\times V)$ to this group, but
it will not always be onto (p\pageref{ac10}).
\end{plain}

\begin{plain}
\label{ct15}J-M. Fontaine has pointed out to us that there is an
incompatibility between our definitions in characteristic $0$ and in
characteristic $p$: in the first case, we use the \'{e}tale cohomology to
define the integral structure at every prime, whereas in characteristic $p$
we use the crystalline cohomology at the prime $p$; but the crystalline
cohomology in characteristic $p$ corresponds to the de Rham cohomology in
characteristic zero, and the integral structures on the de Rham cohomology
and the $p$-adic \'{e}tale cohomology do not agree. We do not know how to
resolve this problem (or even whether it is resolvable).
\end{plain}

\subsection{Variants\label{variants}}

We list some alternative categories $\mathcal{M}{}(k;\mathbb{Q}{})$ to which
our construction applies.

\begin{plain}
\label{ct16}Let $\mathcal{S}{}$ consist of the smooth projective varieties
over $k$ for which the K\"{u}nneth projectors are algebraic. This class
satisfies the condition \ref{af1}, includes all abelian varieties and
surfaces, and includes all smooth projective varieties when $k$ is finite.
Let $\mathcal{M}{}_{l}(k;\mathbb{Q}{})$ and $\mathcal{M}{}(k;\mathbb{Q}{})$
be the categories of motives based on $\mathcal{S}{}$ using respectively the
algebraic classes modulo $l$-adic homomological equivalence and numerical
equivalence as the correspondences. There is a tensor functor%
\begin{equation*}
\mathcal{M}{}_{l}(k)\rightarrow \mathcal{M}{}(k)
\end{equation*}%
which, according to Andr\'{e} and Kahn 2001, has a tensor section $s_{l}$.
Therefore, $\omega _{l}\circ s_{l}$ is a fibre functor on $\mathcal{M}{}(k;%
\mathbb{Q}{})$. Moreover, $s_{l}$ can be chosen so that $(\omega _{l}\circ
s_{l})(hV)=H^{\ast }(V,\mathbb{Q}{}_{l})$. If, as seems likely, $s_{l}$ can
be chosen\footnote{%
According to Andr\'{e} (email 25.03.02) it does not appear possible to prove
this by the methods of their paper.} so that the action of $Z_{\text{num}%
}^{\dim V}(V\times V)\subset \End(hV)$ on $(\omega _{l}\circ s_{l})(hV)$
preserves the image of $H^{\ast }(V,\mathbb{Z}{}_{l})$ in $H^{\ast }(V,%
\mathbb{Q}{}_{l})$ for all $V$ with torsion-free cohomology, then the $%
\omega _{l}\circ s_{l}$ for $l\neq p$ are the $l$-components of an $\mathbb{A%
}_{f}^{p}$-valued fibre functor on $\mathcal{M}{}(k;\mathbb{Q}{})$, and we
can construct $\mathcal{M}{}^{+}(k;\mathbb{Z})$ and $\mathcal{M}{}(k;\mathbb{%
Z}{})$ as above.
\end{plain}

\begin{plain}
\label{ct17}Let $\mathcal{S}{}$ consist of all smooth projective varieties
over $k$. For any $V$ in $\mathcal{S}{}$, the K\"{u}nneth projectors are
almost algebraic (in the sense of Tate 1994, p76). Therefore, the discussion
in \ref{ct16} applies to the category $\mathcal{M}{}(k;\mathbb{Q}{})$
defined using numerical equivalence classes of almost algebraic classes
rather than algebraic classes.
\end{plain}

\begin{plain}
\label{cm2}Recall that a Lefschetz class modulo $\sim $ on a variety $V$ is
an element of the $\mathbb{Q}{}$-subalgebra of $\oplus Z_{\sim }^{r}(V)$
generated by divisor classes. Let $\mathcal{S}{}$ be the smallest class of
varieties over $k$ satisfying (\ref{af1}) and including the abelian
varieties. Then it is possible to define a category $\mathcal{M}{}(k;\mathbb{%
Q}{})$ of isomotives based on $\mathcal{S}{}$ using the Lefschetz classes
modulo numerical equivalence as correspondences; moreover, $\mathcal{M}{}(k;%
\mathbb{Q}{})$ is a semisimple Tannakian category over $\mathbb{Q}{}$ and
there exist canonical functors $\omega _{l}$ (Milne 1999, \S 1).
\end{plain}

\begin{plain}
\label{cm2a}For any field $k$ of characteristic zero, define $\mathcal{M}%
{}(k;\mathbb{Q}{})$ to be the category of motives based on some class $%
\mathcal{S}{}$ satisfying (\ref{af1}) and using the absolute Hodge classes
as correspondences (Deligne and Milne 1982, \S 6).
\end{plain}

\begin{plain}
\label{cm3}Let $\CM(\mathbb{Q}^{\text{al}})$ be the category of motives
based on the abelian varieties of CM-type over $\mathbb{Q}{}^{\text{al}}$
using the absolute Hodge classes as the correspondences, and let $\mathcal{M}%
{}(\mathbb{F}{};\mathbb{Q}{})$ be the quotient of $\CM(\mathbb{Q}^{\text{al}%
})$ corresponding to a fibre functor on $\CM(\mathbb{Q}^{\text{al}})^{P}$
which is correct at every prime (see Milne 2002; here $P$ is the Weil-number
protorus). As above, we can define $\mathcal{M}{}^{+}(\mathbb{F}{};\mathbb{Z}%
{})\mathcal{\ }$and $\mathcal{M}{}(\mathbb{F}{};\mathbb{Z})\,$. ${}$Almost
by definition, the Tate conjecture holds for $\mathcal{M}{}(\mathbb{F}{};%
\mathbb{Q}{})$.
\end{plain}

\section{Tannakian Properties}

Throughout this section, $R$ is a commutative noetherian ring.

\subsection{Tannakian categories over rings}

In Saavedra 1972, Tannakian categories are defined only over a base field.
We adapt his definition (ibid. III 1.1.1) of an ind-Tannakian category to
obtain the notion of a Tannakian category over an arbitrary noetherian ring.

Recall that a \emph{tensor category }$\mathcal{C}{}$ is a category together
with a bifunctor $\otimes \colon \mathcal{C}{}\times \mathcal{C}%
{}\rightarrow \mathcal{C}{}$ and compatible associativity and commutativity
constraints for which there exists an identity object $\1$. When $\End(\1)=R$%
, we call $(\mathcal{C}{},\otimes )$ a tensor category over $R$ or a tensor $%
R$-category. Let $R^{\prime }$ be an $R$-algebra. An $R^{\prime }\emph{%
-valued}$ \emph{fibre functor }is an $R$-linear exact tensor functor $%
\mathcal{C}{}\rightarrow \Mod(R^{\prime })$.

\subsubsection{Some linear algebra}

Let $\mathcal{C}{}$ be an $R$-linear category, and let $R^{\prime }$ be an $%
R $-algebra. When $\mathcal{C}{}$ has arbitrary direct limits, Saavedra
(1972, I 1.5) constructs an $R^{\prime }$-linear category $\mathcal{C}%
{}_{(R^{\prime })}$ by \textquotedblleft extension of the base
ring\textquotedblright . We adapt his construction to the case that $%
R^{\prime }$ is a finite $R$-algebra and $\mathcal{C}{}$ has finite direct
limits (for example, $\mathcal{C}{}$ is abelian).

Thus, let $\mathcal{C}{}$ be an abelian $R$-linear category. Let $R^{\prime
} $ be a finite $R$-algebra, and define $\mathcal{C}_{(R^{\prime })}$ to be
the category whose objects are pairs $(X,i_{X})$ comprising an object $X$ of
$\mathcal{C}$ and a homomorphism $i_{X}\colon R^{\prime }\rightarrow \End%
_{R}(X)$ of $R$-algebras (Saavedra 1972, II 1.5).

\begin{plain}
\label{ct1}For any finitely generated $R$-module $M$, there is a functor
\begin{equation*}
X\mapsto M\otimes _{R}X\colon \mathcal{C}{}\rightarrow \mathcal{C}{}
\end{equation*}%
and an isomorphism (natural in $X$ and $Y$)
\begin{equation*}
\Hom_{\mathcal{C}{}}(M\otimes _{R}X,Y)\rightarrow \Hom_{R}(M,\Hom_{\mathcal{C%
}{}}(X,Y))
\end{equation*}%
(Saavedra 1972, II 1.5.1.1).
\end{plain}

\begin{plain}
\label{ct2}The category $\mathcal{C}{}_{(R^{\prime })}$ is an abelian $%
R^{\prime }$-linear category. A sequence%
\begin{equation*}
(X^{\prime },i_{X^{\prime }})\rightarrow (X,i_{X})\rightarrow (X^{\prime
\prime },i_{X^{\prime \prime }})
\end{equation*}%
in $\mathcal{C}_{(R^{\prime })}$ is exact if and only if%
\begin{equation*}
X^{\prime }\rightarrow X\rightarrow X^{\prime \prime }
\end{equation*}%
is exact in $\mathcal{C}{}$. The functor $j^{R^{\prime }/R}\colon \mathcal{C}%
\rightarrow \mathcal{C}_{(R^{\prime })}$ sending $X$ to $R^{\prime }\otimes
_{R}X$ (with its canonical $R^{\prime }$-structure) is left adjoint to the
functor $j_{R/R}\colon \mathcal{C}_{(R^{\prime })}\rightarrow \mathcal{C}{}$
sending $(X,i_{X})$ to $X$ (Saavedra 1972, II 1.5.2). For any $X,Y\in \ob(%
\mathcal{C}{})$,
\begin{equation}
\Hom_{\mathcal{C}{}_{(R^{\prime })}}(j^{R^{\prime }/R}X,j^{R^{\prime
}/R}Y)\cong \Hom_{\mathcal{C}{}}(X,Y)\otimes _{R}R^{\prime },  \label{e23}
\end{equation}%
and so $j^{R^{\prime }/R}$ realizes $\mathcal{C}{}_{R^{\prime }}$ as a full
subcategory of $\mathcal{C}_{(R^{\prime })}$.
\end{plain}

\begin{plain}
\label{ct3}When $R^{\prime }$ is a flat $R$-algebra, the functor $%
j^{R^{\prime }/R}\colon \mathcal{C}{}\rightarrow \mathcal{C}{}_{(R^{\prime
})}$ is exact.
\end{plain}

\begin{plain}
\label{ct4}Let $\mathcal{C}{}$ be an abelian tensor category over $R$ such
that $\otimes $ is right exact. The category $\mathcal{C}_{(R^{\prime })}$
has a canonical structure of a tensor category over $R^{\prime }$ for which $%
j^{R^{\prime }/R}$ and $j_{R^{\prime }/R}$ are tensor functors; moreover $%
\otimes $ on $\mathcal{C}{}_{(R^{\prime })}$ is right exact (Saavedra 1972,
II 1.5.4). Let $\mathcal{D}{}$ be an abelian tensor category over $R^{\prime
}$ for which $\otimes $ is right exact. The map $u\mapsto u\circ
j^{R^{\prime }/R}$ defines an equivalence of the category of right exact $%
R^{\prime }$-linear tensor functors $\mathcal{C}{}_{(R^{\prime
})}\rightarrow \mathcal{D}{}$ with the category of right exact $R$-linear
tensor functors $\mathcal{C}{}\rightarrow \mathcal{D}{}$. When $R^{\prime }$
is a flat $R$-algebra, $u\mapsto u\circ j^{R^{\prime }/R}$ carries exact
functors to exact functors. (Ibid. II 1.5.3, 1.5.4.) (But it is not known
that $u\circ j^{R^{\prime }/R}$ exact implies $u$ exact, even when $%
R^{\prime }$ is faithfully flat over $R$, except in the case that $k$ and $%
k^{\prime }$ are fields, when the proof requires the main theorem of Deligne
1990).
\end{plain}

\subsubsection{Definitions}

We say that an abelian tensor $R$-category $\mathcal{C}{}$ is \emph{finitely
generated} if there exists a finite set of objects of $\mathcal{C}{}$
contained in no proper strictly full abelian tensor $R$-subcategory of $%
\mathcal{C}{}$.

\begin{definition}
\label{ct5}(a) A \emph{neutral Tannakian category over }$R$ is a tensor
category over $R$ that is tensor equivalent to $\Rep(G;R)$ for some flat
affine group scheme $G$ over $R$.

(b) An \emph{algebraic} \emph{Tannakian category over }$R$ is a finitely
generated abelian tensor category $\mathcal{C}{}$ over $R$ such that $%
\mathcal{C}_{(R^{\prime })}$ becomes a neutral Tannakian category for some
faithfully flat finite $R$-algebra $R^{\prime }$.

(c) A \emph{Tannakian category over} $R$ is an abelian tensor category $%
\mathcal{C}{}$ over $R$ that is a filtered union of algebraic Tannakian
categories.
\end{definition}

Note that, in a Tannakian category, the Homs${}$ are finitely generated $R$%
-modules, because
\begin{equation*}
\Hom_{\mathcal{C}{}}(X,Y)\otimes _{R}R^{\prime }\overset{(\ref{e23})}{\cong }%
\Hom_{\mathcal{C}{}_{(R^{\prime })}}(j^{R^{\prime }/R}X,j^{R^{\prime
}/R}Y)\subset \Hom_{R^{\prime }}(\omega (X),\omega (Y))\text{.}
\end{equation*}

\subsection{A criterion to be a Tannakian category over a Dedekind domain}

Let $\mathcal{C}{}$ be an abelian tensor category over a noetherian ring $R$%
. When internal Homs exist in $\mathcal{C}{}$ (Saavedra 1972, I 3), we let $%
X^{\vee }=\underline{\Hom}(X,\1)$. There is a canonical morphism $%
X\rightarrow X^{\vee \vee }$, and $X$ is said to be \emph{reflexive }when
this morphism is an isomorphism.

For example, when $G$ is an affine group scheme over $R$, an object $X$ of $%
\Rep(G;R)$ is reflexive if and only if its underlying $R$-module is
projective. Serre (1968, 2.2) proves that, when $R$ is a Dedekind domain,
every object in $\Rep(G;R)$ is a quotient of a representation of $G$ on a
projective $R$-module. Thus every object in $\Rep(G;R)$ is a quotient of a
reflexive object. This property can be used to characterize Tannakian
categories over Dedekind domains.

\begin{proposition}
\label{ct6}Let $\mathcal{C}{}$ be an abelian tensor category over a Dedekind
domain $R$ such that every object of $\mathcal{C}{}$ is a quotient of a
reflexive object. If there exists an $R^{\prime }$-valued fibre functor on $%
\mathcal{C}{}_{(R^{\prime })}$ for some Dedekind domain $R^{\prime }$ finite
and flat over $R$, then $\mathcal{C}{}$ is Tannakian.
\end{proposition}

\begin{proof}
Let $\omega \colon \mathcal{C}{}_{(R^{\prime })}\rightarrow \Mod(R^{\prime
}) $ be such a fibre functor. To prove that $\mathcal{C}{}$ is Tannakian, we
shall show that the functor $\underline{\Aut}^{\otimes }(\omega )$ of $%
R^{\prime }$-algebras is representable by a flat affine group scheme $G$
over $R^{\prime }$, and that $\omega $ defines an equivalence of $R^{\prime
} $-tensor categories $\mathcal{C}_{(R^{\prime })}\rightarrow \Rep%
(G;R^{\prime })$.

When $R^{\prime }=R$, this is proved in Saavedra 1972, II 4.1. In detail:
ibid.~ II 4.1.1 shows that $\underline{\End}^{\otimes }(\omega )$ is
represented by a flat affine monoid scheme $G$ over $R$ and that $\omega $
defines an equivalence of $R$-linear tensor categories $\mathcal{C}%
\rightarrow \Rep(G;R)$; but $\underline{\End}^{\otimes }(\omega )\cong
\underline{\End}^{\otimes }(\omega _{0})$ where $\omega _{0}$ is the
restriction of $\omega $ to the full subcategory of reflexive objects; this
last category is rigid, and so every endomorphism of $\omega _{0}$ is an
automorphism (ibid.~ I 5.2.3), which implies that $G$ is a group scheme.

According to (\ref{ct2}) and (\ref{ct4}), $\mathcal{C}{}_{(R^{\prime })}$ is
an abelian tensor category over $R^{\prime }$. Thus, the general case will
follow from the neutral case once we show that $\mathcal{C}_{(R^{\prime })}$
inherits the property that every object is a quotient of a reflexive object.
Let $(Y,i_{Y})$ be an object of $\mathcal{C}_{(R^{\prime })}$, and let $%
q\colon X\rightarrow Y$ express $Y$ as a quotient of a reflexive object. The
composite%
\begin{equation*}
j^{R^{\prime }/R}X\xr{j^{R'/R}(q)} j^{R^{\prime }/R}Y\rightarrow (Y,i_{Y})
\end{equation*}%
is again surjective, and $j^{R^{\prime }/R}X$ is reflexive because $X$ is.
\end{proof}

\begin{example}
\label{ct6m}For $\ell \neq p$, the category $\mathcal{R}{}(k;\mathbb{Z}%
{}_{\ell })$ is Tannakian over $\mathbb{Z}{}_{\ell }$ (apply \ref{pl10}).
\end{example}

\begin{example}
\label{ct7m}(a) Let $k$ be a finite field, and let $\omega $ be the
forgetful functor $\mathcal{R}{}(k;\mathbb{Z}{}_{p})\rightarrow \Mod(W(k))$.
For $(X,i_{X})$ in $\mathcal{R}{}(k;\mathbb{Z}{}_{p})_{(W(k))}$, there are
two actions of $W(k)$ on $\omega (X)$, namely, that coming from the action $%
i_{X}$ of $W(k)$ on $X$ and that coming from the action of $W(k)$ on $\omega
$, and we define $\omega ^{\prime }(X,i_{X})=W(k)\otimes _{W(k)\otimes _{%
\mathbb{Z}{}_{p}}W(k)}\omega (X)$. Then $\omega ^{\prime }$ is exact.
Indeed, $W(k)\otimes _{\mathbb{Z}{}_{p}}W(k)$ is isomorphic to a product of
copies of $W(k)$ indexed by the elements of $\Gal(k/\mathbb{F}_{p})$, and
the map $W(k)\otimes _{\mathbb{Z}{}_{p}}W(k)\rightarrow W(k)$ corresponds to
the projection to the \textquotedblleft $\id$\textquotedblright -component.
Therefore $\omega ^{\prime }$ is a $W(k)$-valued fibre functor on $\mathcal{R%
}{}(k;\mathbb{Z}{}_{p})_{(W(k))}$, and $\mathcal{R}{}(k;\mathbb{Z}{}_{p})$
is Tannakian (apply \ref{pl13p}(b)).

(b) When $k$ has nonzero characteristic $p$ and is infinite, then $\mathcal{R%
}{}(k;\mathbb{Z}{}_{p})$ is not Tannakian (at least according to our
definition, which may be too strict), because there exist objects $X$ for
which $\End(X)$ is not a finitely generated $\mathbb{Z}{}_{p}$-module.
\end{example}

\subsection{Existence of fibre functors}

We refer to Deligne 1989, \S 6, for the notion of the fundamental group of a
Tannakian category $\mathcal{T}{}$: it is an affine group scheme $\pi (%
\mathcal{T}{})$ in $\Ind\mathcal{T}{}$ acting on the objects of $\mathcal{T}%
{}$; each fibre functor $\omega $ carries $\pi (\mathcal{T}{})$ to $%
\underline{\Aut}^{\otimes }(\omega )$ and the action of $\pi (\mathcal{T}{})$
on an object $X$ to the action of $\underline{\Aut}^{\otimes }(\omega )$ to $%
\omega (X)$. When $\pi (\mathcal{T}{})$ is commutative, it lies in $\Ind%
\mathcal{T}{}^{0}$ where $\mathcal{T}{}^{0}$ is the full subcategory of
trivial objects. Since $\Hom(\1,-)\colon \mathcal{T}{}^{0}\rightarrow \Vct%
(k) $ is an equivalence of categories, in this case $\pi (\mathcal{T}{})$
can be identified with an affine group scheme in the usual sense.

\begin{proposition}
\label{ct8}Let $\mathcal{M}$ be a Tannakian category over $\mathbb{Q}{}$
whose fundamental group $T$ is an algebraic group of multiplicative type,
and let $\omega ^{p}$ be a $\mathbb{Q}{}$-linear tensor functor $\mathcal{M}%
\rightarrow \Mod(\mathbb{A}_{f}^{p})$ such that $\omega _{\ell }=_{\text{df}%
}\omega ^{p}\otimes _{\mathbb{A}_{f}^{p}}\mathbb{Q}{}_{\ell }$ is a $\mathbb{%
Q}{}_{\ell }$-valued fibre functor for all $\ell \neq p$. \noindent Then
there exists a finite field extension $L$ of $\mathbb{Q}{}$, a fibre functor
$\omega \colon \mathcal{M}\rightarrow $ $\Vct(L)$, and an isomorphism $%
\omega \otimes _{L}\mathbb{A}{}_{f,L}^{p}\rightarrow \omega ^{p}\otimes _{%
\mathbb{A}{}_{f}^{p}}\mathbb{A}{}_{f,L}^{p}$ of tensor functors $\mathcal{M}%
\rightarrow \Mod(\mathbb{A}_{f,L}^{p})$.
\end{proposition}

\begin{proof}
The hypotheses imply that $\omega ^{p}(X)$ is a free $\mathbb{A}_{f}^{p}$%
-module for all $X$ (see the proof of \ref{af4}).

Assume first that $T$ is split over $\mathbb{Q}$ and that there exists a
fibre functor $\omega \colon \mathcal{M}\rightarrow \Vct(\mathbb{Q}{})$. We
show in this case that there is an isomorphism $\omega \otimes _{\mathbb{Q}%
{}}\mathbb{A}_{f}^{p}\rightarrow \omega ^{p}$. We use $\omega $ to identify $%
\mathcal{M}$ with $\Rep(T;\mathbb{Q})$. Let $\Xi $ be a basis for the group $%
X^{\ast }(T)$ of characters of $T$, and let $X$ be the representation $%
X=\bigoplus_{\lambda \in \Xi }X_{\lambda }$ where $X_{\lambda }$ is a
one-dimensional representation of $T$ with character $\lambda $. Choose a
graded lattice $\Lambda $ in $X$, i.e., a $\mathbb{Z}{}$-lattice such that $%
\Lambda =$ $\bigoplus_{\lambda \in \Xi }\Lambda \mathbb{\cap }{}X_{\lambda }$%
. Likewise, choose a graded lattice $\Lambda _{f}$ in $\omega ^{p}(X)$. Then
$\Lambda _{\ell }=_{\text{df}}\Lambda _{f}\otimes _{\mathbb{A}_{f}^{p}}%
\mathbb{Z}{}_{\ell }\mathbb{\ }$is a graded lattice in $\omega _{\ell }(X)$.
As $T$ is a split torus, $H^{1}(\mathbb{Q}{}_{\ell },T)=0$, and so the
theory of Tannakian categories shows that there exists an isomorphism of
fibre functors $\alpha _{\ell }\colon \omega \otimes _{\mathbb{Q}{}}\mathbb{Q%
}_{\ell }\rightarrow \omega _{\ell }$ which is uniquely determined up to an
element of $T(\mathbb{Q}_{\ell })$. We may choose $\alpha _{\ell }$ to map $%
\Lambda \otimes _{\mathbb{\mathbb{Z}{}}}\mathbb{Z}{}_{\ell }$ onto $\Lambda
_{\ell }$. Now the family $(\alpha _{\ell })$ defines an isomorphism $\omega
\otimes _{\mathbb{Q}{}}\mathbb{A}_{f}^{p}\rightarrow \omega ^{p}$.

In the general case, there will exist a finite field extension $L$ of $%
\mathbb{Q}{}$ such that $T$ splits over $L$ and such that there exists a
fibre functor $\omega \colon \mathcal{M}\rightarrow \Vct(L)$. Then $\omega $
defines a fibre functor $\mathcal{M}_{(L)}\rightarrow \Vct(L)$ (Deligne and
Milne 1982, 3.11), and the previous argument applied to $\mathcal{M}_{(L)}$
proves the general case of the proposition.
\end{proof}

\subsection{The category $\mathcal{M}{}(k;\mathbb{Z}{})$}

Fix a field $k$, and let $\mathcal{M}{}(k;\mathbb{Q}{})$ and $\mathcal{M}%
{}(k;\mathbb{Z}{})$ be as in \S 5.

\begin{lemma}
\label{ct8f}The category $\mathcal{M}{}(k;\mathbb{Z}{})$ has internal Homs.
\end{lemma}

\begin{proof}
Since each of the categories $\mathcal{R}{}(k;\mathbb{\hat{Z})}{}$, $%
\mathcal{R}{}(k;\mathbb{A}_{f}^{p})$, $\mathcal{M}{}(k;\mathbb{Q}{})$ has
internal Homs, it is straightforward to check this.
\end{proof}

\begin{lemma}
\label{ct8m}Every motive in $\mathcal{M}{}(k;\mathbb{Z})$ is a quotient of a
reflexive motive.
\end{lemma}

\begin{proof}
The proof is similar to that of Proposition \ref{pl13p}: a motive is
reflexive if and only if it is torsion-free; every motive is a direct sum of
a torsion motive and a torsion-free motive; every torsion motive is a
quotient of $\1(r)^{s}$ for some $r$ and $s$.
\end{proof}

\begin{theorem}
\label{ct9}The category $\mathcal{M}{}(\mathbb{F}{}_{q};\mathbb{Z}{})$ is
Tannakian.
\end{theorem}

\begin{proof}
We apply Proposition \ref{ct6}. It remains to prove the existence of fibre
functors. We replace $\mathcal{S}$ with a subset generated by a finite set
of varieties, and prove that $\mathcal{M}{}(k;\mathbb{Z}{};\mathcal{S}%
{})_{(R^{\prime })}$ admits an $R^{\prime }$-valued fibre functor for some
Dedekind domain $R^{\prime }$ finite and flat over $R$. For this, we can
apply Proposition \ref{ct8} because the fundamental group of $\mathcal{M}%
{}=_{\text{df}}\mathcal{M}(k;\mathbb{Q};\mathcal{S}{})$ is of multiplicative
type: it is reductive because $\mathcal{M}{}$ is semisimple, and it is
commutative because its extension to $\mathbb{Q}{}_{\ell }$ is contained in
the algebraic group fixing the Tate classes, which is commutative.

According to Proposition \ref{ct8}, there exists a finite field extension $L$
of $\mathbb{Q}{}$ splitting $T$, a fibre functor $\omega _{0}\colon \mathcal{%
M}\rightarrow \Vct(L)$, and an isomorphism $\xi ^{p}\colon L\otimes _{%
\mathbb{Q}{}}\omega ^{p}\rightarrow \mathbb{A}_{f}^{p}\otimes _{\mathbb{Q}%
{}}\omega _{0}$. After possibly extending $L$, we may suppose that there is
also an isomorphism $\xi _{p}\colon L\otimes _{\mathbb{Q}{}}\omega
_{p}\rightarrow \mathbb{Q}{}_{p}\otimes _{\mathbb{Q}{}}\omega _{0}$ where $%
\omega _{p}$ is $X\mapsto X_{p}$ (regarded as a $B(\mathbb{F}{}_{q})$%
-module). Let $R$ be the integral closure of $\mathbb{Z}{}$ in $L$.

For $X_{0}$ in $\mathcal{M}_{(L)}$, there are two actions of $L$ on $\omega
_{0}(X_{0})$, namely, that coming from the action of $L$ on $\omega _{0}$
and that coming from the action of $L$ on $X_{0}$. Define
\begin{equation*}
\omega _{1}(X_{0})=L\otimes _{L\otimes _{\mathbb{Q}{}}L}\omega _{0}(X_{0})%
\text{.}
\end{equation*}%
Then $\omega _{1}$ is an $L$-valued fibre functor on $\mathcal{M}{}_{(L)}$
(Deligne and Milne 1982, 3.11).

For $(X,i_{X})$ in $\mathcal{M}_{(R)}$, set
\begin{equation*}
\omega _{2}(X)=(X^{p},\omega ^{\prime }(X_{p}))
\end{equation*}%
where $\omega ^{\prime }$ is the functor defined (\ref{ct7m}a). Then $\omega
_{2}$ is an $\hat{R}$-valued fibre functor on $\mathcal{M}{}_{(R)}$, and $%
\xi $ defines an isomorphism $L\otimes _{\mathbb{Q}{}}\omega _{2}\rightarrow
\mathbb{A}_{f}\otimes _{\mathbb{Q}{}}\omega _{1}$.

Thus, we have defined an exact tensor functor%
\begin{equation*}
X\mapsto (\omega _{2}(X),\omega _{1}(X_{0}),\xi )\colon \mathcal{M}%
{}_{(R)}\rightarrow \Modf(\hat{R})\times _{\Modf(\mathbb{A}_{f}\otimes R)}%
\Modf(L)\text{.}
\end{equation*}%
Its image lies in the full subcategory of triples whose $\Modf(\hat{R})$%
-component has finite torsion, which can be identified with the category $%
\Modf(R)$. It is therefore an $R$-valued fibre functor on $\mathcal{M}{}(%
\mathbb{F}{}_{q};\mathbb{Z}{})_{(R)}$.
\end{proof}

\begin{remark}
\label{ct10}If $\mathcal{M}{}(k;\mathbb{Q}{})$ is as in \ref{cm2} or \ref%
{cm3} and $k$ is finite, or as in \ref{cm2a}, then $\mathcal{M}{}(k;\mathbb{Z%
}{})$ is Tannakian. In the first two cases, this can be proved as in Theorem %
\ref{ct9}, and in the third case, there is a Betti fibre functor on $%
\mathcal{M}{}(k;\mathbb{Z}{})$.
\end{remark}

\subsection{Stacks of categories of motives}

We refer the reader to Saavedra 1972, I 4.5, for the notion of \emph{tensor
fibred category} $p\colon \mathcal{M}{}\rightarrow \mathcal{E}{}$. Briefly,
it is a fibred category $p\colon \mathcal{M}{}\rightarrow \mathcal{E}{}$
endowed with an $\mathcal{E}{}$-bifunctor $\otimes \colon \mathcal{M}%
{}\times _{\mathcal{E}{}}\mathcal{M}{}\rightarrow \mathcal{M}{}$ and
compatible associativity and commutativity constraints for which there
exists a cartesian section to $p$ giving identity objects in each fibre;
thus, for each $E$ in $\mathcal{E}{}$, the fibre $\mathcal{M}{}_{E}$ is a
tensor category. When $\mathcal{E}{}$ is endowed with a topology (in the
sense of Grothendieck) and $p\colon \mathcal{M}{}\rightarrow \mathcal{E}$ is
a stack, we call it a \emph{stack of tensor categories}. When the fibres are
abelian categories and the inverse image functors are exact, we call $%
p\colon \mathcal{M}{}\rightarrow \mathcal{E}$ a \emph{stack of abelian
tensor categories}. Finally, when $\End(\1)=R$ for each identity object $\1$
in each fibre, we call $p\colon \mathcal{M}{}\rightarrow \mathcal{E}$ a
\emph{stack of abelian tensor }$R$-\emph{categories}. In particular, then
each fibre is an abelian tensor $R$-category.

For a scheme $S=\tcoprod \Spec k_{i}$, finite and \'{e}tale over $\mathbb{F}%
{}_{q}$, let $\mathcal{M}{}(S;\mathbb{Z}{})=\prod \mathcal{M}{}(k_{i};%
\mathbb{Z}{})$.

\begin{theorem}
\label{ct13}Let $\Et_{\mathbb{F}{}_{q}}$ be the category of all schemes
finite and \'{e}tale over $\mathbb{F}{}_{q}$, and endow $\Et_{\mathbb{F}%
{}_{q}}$ with the \'{e}tale topology. For $k=\mathbb{F}{}_{q}$ or $\mathbb{F}%
{}$, let $\mathcal{M}{}(k;\mathbb{Q}{})$ be as at the start of \S 5 and
assume the Tate conjecture holds for $\mathcal{M}{}(k;\mathbb{Q}{})$. The
categories $\mathcal{M}{}(-;\mathbb{Z}{})$ form (in a natural way) the
fibres of a stack of abelian tensor $\mathbb{Z}{}$-categories over $\Et_{%
\mathbb{F}{}_{q}}$.
\end{theorem}

\begin{proof}
Immediate consequence of Proposition \ref{ct12}.
\end{proof}

\section{Some Spectral Sequences}

\subsection{Abstract spectral sequences}

We prove an abstract version of Tate's spectral sequence (Milne 1986b, I
0.3). We fix a field $k$, and let $\Et_{k}$ be the category of all schemes
finite and \'{e}tale over $k$. Also, we fix a ring $R$, and let $\mathcal{M}%
{}\rightarrow \Et_{k}$ be a stack of noetherian abelian tensor $R$%
-categories with internal Homs over $\Et_{k}$.

\begin{plain}
\label{ha0}When $X$ and $Y$ are objects of an abelian category $\mathcal{A}%
{} $, $\Ext^{r}(X,Y)$ denotes the Yoneda extension group (Mitchell 1965,
VII). It agrees with the group defined in terms of the derived category
(Verdier 1996, III 3.2.12). When $\mathcal{A}{}$ is noetherian, the Yoneda
extension group also agrees with that defined using injective resolutions in
the ind-category $\Ind(\mathcal{A}{})$ (Oort 1964, p.~ 235; see also Huber
1993).
\end{plain}

\begin{plain}
\label{ha0m}An object of a tensor category is \emph{trivial }if it is a
quotient of a finite direct sum of copies of $\1$, and an object $Y$ of an
abelian tensor category is \emph{flat }if the functor $X\mapsto X\otimes Y$
is exact. When $X$ and $Y$ are objects of an abelian tensor category over $R$%
, we write
\begin{equation*}
\Tor_{i}^{R}(X,Y)=0\text{ for all }i>0
\end{equation*}%
to mean $Y$ has a resolution
\begin{equation*}
\cdots \rightarrow Y^{1}\rightarrow Y^{0}\rightarrow Y\rightarrow 0
\end{equation*}%
by flat objects which remains exact when it has been tensored by $X$.
\end{plain}

\begin{plain}
\label{ha1}Let $k^{\prime }$ be a finite Galois extension of $k$ with Galois
group $\Gamma $, and write $X\mapsto X^{\prime }$ for an inverse image
functor $\mathcal{M}{}(k)\rightarrow \mathcal{M}{}(k^{\prime })$. Because $%
\mathcal{M}{}\rightarrow \Et_{k}$ is a stack, the functor
\begin{equation*}
X\mapsto N(X)\overset{\text{df}}{=}\Hom(\1,X^{\prime })
\end{equation*}%
defines an equivalence from the full subcategory of $\mathcal{\mathcal{M}{}}%
(k)$ of objects becoming trivial in $\mathcal{M}{}(k^{\prime })$ to the
category $\Modf(R[\Gamma ])$ with quasi-inverse $N\mapsto N\otimes \1$.
\end{plain}

\begin{plain}
\label{ha2}Let $k^{\prime }$ be a finite \'{e}tale $k$-algebra. The theory
of \textquotedblleft Weil restriction of scalars\textquotedblright\ (Weil
1982, 1.3) applies in the present situation and gives an exact functor $%
\mathcal{M}{}(k^{\prime })\rightarrow \mathcal{M}{}(k)$ that is both left
and right adjoint to the inverse image functor $\mathcal{M}{}(k)\rightarrow
\mathcal{M}{}(k^{\prime })$.
\end{plain}

\begin{proposition}
\label{ha3}Let $k^{\prime }$ be a finite Galois extension of $k$ with Galois
group $\Gamma $. Write $T^{\prime }$ for the inverse image in $\mathcal{M}%
{}(k^{\prime })$ of an object $T$ of $\mathcal{M}{}(k)$. For all $X,Y,Z$ in $%
\mathcal{M}{}(k)$ such that $X^{\prime }$ is trivial and $\Tor%
_{i}^{R}(X,Y)=0 $ for all $i>0$, there is a spectral sequence%
\begin{equation*}
\Ext_{R[\Gamma ]}^{r}(N(X),\Ext_{\mathcal{M}{}(k^{\prime })}^{s}(Y^{\prime
},Z^{\prime }))\Longrightarrow \Ext_{\mathcal{M}{}(k)}^{r+s}(X\otimes Y,Z).
\end{equation*}
\end{proposition}

Consider the functors (for fixed $X,Y$)%
\begin{equation*}
\begin{array}{lll}
\alpha \colon \mathcal{M}{}(k)\rightarrow \Mod(R[\Gamma ]), & \quad & \alpha
(Z)=\Hom_{\mathcal{M}{}(k^{\prime })}(Y^{\prime },Z^{\prime }) \\
\beta \colon \Mod(R[\Gamma ])\rightarrow \Ab & \quad & \beta (M)=\Hom%
_{R[\Gamma ]}(N(X),M) \\
\gamma \colon \mathcal{M}{}(k)\rightarrow \Ab & \quad & \gamma (Z)=\Hom_{%
\mathcal{M}{}(k)}(X\otimes Y,Z).%
\end{array}%
\end{equation*}

To prove the proposition, we shall show:

\begin{enumerate}
\item $\beta \circ \alpha \cong \gamma $;

\item the functor $\alpha $ takes injective objects in $\Ind(\mathcal{M}%
{}(k))$ to $\beta $-acyclic $R[\Gamma ]$-modules;

\item the functor $Z\mapsto Z^{\prime }\colon \mathcal{M}{}(k)\rightarrow
\mathcal{M}{}(k^{\prime })$ takes injective objects in $\Ind(\mathcal{M}%
{}(k))$ to injective objects in $\Ind(\mathcal{M}{}(k^{\prime })).$
\end{enumerate}

Statements (a) and (b) imply that there is a spectral sequence%
\begin{equation*}
(R^{r}\beta \circ R^{s}\alpha )(Z)\Longrightarrow (R^{r+s}\gamma )(Z);
\end{equation*}%
hence,%
\begin{equation*}
\Ext_{R[\Gamma ]}^{r}(N(X),(R^{s}\alpha )(Z))\Longrightarrow \Ext_{\mathcal{M%
}{}(k)}^{r+s}(X\otimes Y,Z).
\end{equation*}%
Statement (c), together with the exactness of $Z\mapsto Z^{\prime }$, shows
that $(R^{s}\alpha )(Z)=\Ext_{\mathcal{M}{}(k^{\prime })}^{s}(Y^{\prime
},Z^{\prime })$.

\smallskip \noindent \textsc{Proof of }(c). The functor $Y\mapsto Y^{\prime
} $ has an exact left adjoint (\ref{ha2}).\hspace*{\fill}$\square $

\smallskip \noindent \textsc{Proof of }(a). We have to prove%
\begin{equation*}
\Hom_{R[\Gamma ]}(N(X),\Hom_{\mathcal{M}{}(k^{\prime })}(Y^{\prime
},Z^{\prime }))\cong \Hom_{\mathcal{M}{}(k)}(X\otimes Y,Z)
\end{equation*}%
functorially in $Z$. From the definition of $\underline{\Hom}$ (Saavedra
1972, I 3.1.1)%
\begin{equation*}
\Hom_{\mathcal{M}{}(k^{\prime })}(X^{\prime }\otimes Y^{\prime },Z^{\prime
})\cong \Hom_{\mathcal{M}{}(k^{\prime })}(X^{\prime },\underline{\Hom}_{%
\mathcal{M}{}(k^{\prime })}(Y^{\prime },Z^{\prime }))\text{.}
\end{equation*}%
But $X^{\prime }\cong N(X)\otimes \1$, and so
\begin{eqnarray*}
\Hom_{\mathcal{M}{}(k^{\prime })}(X^{\prime }\otimes Y^{\prime },Z^{\prime
}) &\overset{\ref{ct1}}{\cong }&\Hom_{R}(N(X),\Hom(\1,\underline{\Hom}_{%
\mathcal{M}{}(k^{\prime })}(Y^{\prime },Z^{\prime }))) \\
&\cong &\Hom_{R}(N(X),\Hom_{\mathcal{M}{}(k^{\prime })}(Y^{\prime
},Z^{\prime }))\text{.}
\end{eqnarray*}%
On taking $\Gamma $-invariants we obtain (a).\hspace*{\fill}$\square $

\smallskip \noindent \textsc{Proof of }(b). We write $\Hom_{k}$ and $\Hom%
_{k^{\prime }}$ for $\Hom_{\Ind(\mathcal{M}{}(k))}$ and $\Hom_{\Ind(\mathcal{%
M}{}(k^{\prime }))}$ respectively. We have to show that if $I$ is injective
in $\Ind(\mathcal{M}{}(k))$, then
\begin{equation*}
\Ext_{R[\Gamma ]}^{r}(N(X),\Hom_{k^{\prime }}(Y^{\prime },I^{\prime }))=0%
\text{ for }r>0.
\end{equation*}

\smallskip \noindent We first show that, if $Y$ is flat and $I$ is
injective, then $\Hom_{k^{\prime }}(Y^{\prime },I^{\prime })$ is an
injective $R[\Gamma ]$-module. For this, we must show that the functor%
\begin{equation*}
\Hom_{R[\Gamma ]}(-,\Hom_{k^{\prime }}(Y^{\prime },I^{\prime }))\colon \Mod%
(R[\Gamma ])\rightarrow \Ab
\end{equation*}%
is exact, at least on finitely generated $R[\Gamma ]$-modules, but (a) shows
that the composite of this functor with the equivalence $X\mapsto N(X)$ is
the composite of the exact functors $-\otimes Y$ and $\Hom_{k}(-,I)$.

We now prove (b). Let
\begin{equation*}
Y^{\bullet }\rightarrow Y
\end{equation*}%
be a flat resolution of $Y$ such that $X\otimes (Y^{\bullet }\rightarrow Y)$
is exact, and let $I$ be an injective object of $\Ind(\mathcal{M}{}(k))$.
Then, by (c), $I^{\prime }$ is injective, and so%
\begin{equation*}
\Hom_{k^{\prime }}(Y^{\prime },I^{\prime })\rightarrow \Hom_{k^{\prime
}}(Y^{\bullet \prime },I^{\prime })
\end{equation*}%
is an exact complex; in fact, it is an injective resolution of $\Hom%
_{k^{\prime }}(Y^{\prime },I^{\prime })$. We shall use it to compute the
groups $\Ext_{R[\Gamma ]}^{r}(X,\Hom_{k^{\prime }}(Y^{\prime },I^{\prime }))$%
. Statement (a) shows that the complex%
\begin{equation*}
\Hom_{R[\Gamma ]}(N(X),\Hom_{k^{\prime }}(Y^{\prime },I^{\prime
}))\rightarrow \Hom_{R[\Gamma ]}(N(X),\Hom_{k^{\prime }}(Y^{\bullet \prime
},I^{\prime }))
\end{equation*}%
is isomorphic to the complex%
\begin{equation*}
\Hom_{k}(X\otimes Y,I)\rightarrow \Hom_{k}(X\otimes Y^{\bullet },I)
\end{equation*}%
which is exact because $X$ $\otimes (Y\rightarrow Y^{\bullet })$ is exact
and $I$ is injective.\hspace*{\fill}$\square $

\begin{corollary}
\label{ha5}Let $\bar{k}$ be a Galois extension of $k$ (possibly infinite),
and let $\mathcal{M}{}(\bar{k})=\varinjlim_{k^{\prime }\subset \bar{k}}\mathcal{M}%
{}(k^{\prime })$ where $k^{\prime }$ runs over the finite Galois extensions
of $k$. Denote the inverse image functor $\mathcal{M}{}(k)\rightarrow
\mathcal{M}{}(\bar{k})$ by $X\mapsto \bar{X}$, and let $\Gamma =\Gal(\bar{k}%
/k)$. Let $Y$ be an object of $\mathcal{M}{}(k)$ that admits a resolution by
flat objects. For all $Z$ in $\mathcal{M}{}(k)$, there is a spectral sequence%
\begin{equation*}
H^{r}(\Gamma ,\Ext_{\mathcal{M}{}(\bar{k})}^{s}(\bar{Y},\bar{Z}%
))\Longrightarrow \Ext_{\mathcal{M}{}(k)}^{r+s}(Y,Z).
\end{equation*}
\end{corollary}

\begin{proof}
The condition on $Y$ is that $\Tor_{i}^{R}(\1,Y)=0$ for all $i>0$. When $%
\bar{k}$ is finite over $k$, the two definitions of $\mathcal{M}{}(\bar{k})$
coincide, and we obtain the spectral sequence by taking $M=\1$ in
Proposition \ref{ha3}. When $\bar{k}$ is of infinite degree over $k$, we
pass to the limit over subfields $k^{\prime }\subset \bar{k}$ that are
finite and Galois over $k$.
\end{proof}

\subsection{Spectral sequences at $\ell $}

\begin{proposition}
\label{lc.dm}Let $\Gamma =\Gal(\mathbb{F}{}/\mathbb{F}{}_{q})$. For $M,N$ in
$\mathcal{R}{}(\mathbb{F}{}_{q};\mathbb{Z}{}_{\ell })$, there is a spectral
sequence%
\begin{equation}
H^{r}(\Gamma ,\Ext_{\mathcal{R}{}(\mathbb{F}{};\mathbb{Z}{}_{\ell
})}^{s}(M,N))\Longrightarrow \Ext_{\mathcal{R}{}(\mathbb{F}{}_{q};\mathbb{Z}%
{}_{\ell })}^{r+s}(M,N)\text{.}  \label{e6}
\end{equation}
\end{proposition}

\begin{proof}
Proposition \ref{pl10} shows $M$ has a resolution by flat objects, $i>0$,
and so this is a special case of Corollary \ref{ha5}.
\end{proof}

Denote Exts${}{}$ in $\Mod(\mathbb{Z}{}_{\ell })$ by $\Ext_{\mathbb{Z}%
{}_{\ell }}$.

\begin{proposition}
\label{lc.e}Let $\Gamma =\Gal(\mathbb{F}{}/\mathbb{F}{}_{q})$. For $M,N$ in $%
\mathcal{R}{}(\mathbb{F}{}_{q};\mathbb{Z}{}_{\ell })$, there is a spectral
sequence%
\begin{equation}
H_{\text{cts}}^{r}(\Gamma ,\Ext_{\mathbb{Z}{}_{\ell
}}^{s}(M,N))\Longrightarrow \Ext_{\mathcal{R}{}(\mathbb{F}{}_{q};\mathbb{Z}%
{}_{\ell })}^{r+s}(M,N)\text{.}  \label{e5}
\end{equation}%
Here $H_{\text{cts}}^{r}$ is computed using continuous cochains relative to
the $\ell $-adic topology on $\Ext_{\mathbb{Z}{}_{\ell }}^{s}(M,N)$.
\end{proposition}

\begin{proof}
When $N$ is finite, (\ref{e5}) will become (\ref{e6}) once we have shown
that the natural map%
\begin{equation*}
\Ext_{\mathcal{R}{}(\mathbb{F}{};\mathbb{Z}{}_{\ell })}^{j}(M,N))\rightarrow %
\Ext_{\mathbb{Z}{}_{\ell }}^{j}(M,N)
\end{equation*}%
is an isomorphism. When $M$ is finite, this is obvious (with no condition on
$N$), because it is true for $j=0$ and the forgetful functor $\Ind\mathcal{R}%
{}(\mathbb{F}{};\mathbb{Z}{}_{\ell })\rightarrow \Mod(\mathbb{Z}{}_{\ell })$
is exact and preserves injectives. The $\Ext(M,-)$ sequences of%
\begin{eqnarray}
0 &\rightarrow &M\xr{\ell^n}M\rightarrow M/\ell ^{n}M\rightarrow 0\,\text{%
,\thinspace }\quad \ell ^{n}N=0\text{,}  \label{e34} \\
0 &\rightarrow &M_{\text{tors}}\rightarrow M\rightarrow M/M_{\text{tors}%
}\rightarrow 0  \notag
\end{eqnarray}%
allow us to deduce it, first for all torsion-free $M$, and then for all $M$.

There are spectral sequences%
\begin{eqnarray*}
E_{2}^{ij} &=&\Ext_{\mathcal{R}{}(\mathbb{F}{}_{q};\mathbb{Z}_{\ell
})}^{i}(M,\varprojlim{}^{j}N^{(\ell ^{n})}) \\
E_{2}^{i,j} &=&\varprojlim{}^{i}\Ext_{\mathcal{R}{}(\mathbb{F}_{q}{};\mathbb{Z}%
_{\ell })}^{j}(M,N^{(\ell ^{n})}),
\end{eqnarray*}%
with isomorphic final terms (cf. Jensen 1972, 4.3). The higher inverse
limits vanish because the groups $\Ext_{\mathcal{R}{}(\mathbb{F}{}_{q};%
\mathbb{Z}_{\ell })}^{j}(M,N^{(\ell ^{n})})$ and $N^{(\ell ^{n})}$ are
finite. Moreover, $\varprojlim{}N^{(\ell ^{n})}\cong N$, and so this shows that%
\begin{equation}
\Ext_{\mathcal{R}{}(\mathbb{F}{}_{q};\mathbb{Z}_{\ell })}^{i}(M,N)\cong \varprojlim%
_{n}\Ext_{\mathcal{R}{}(\mathbb{F}{}_{q};\mathbb{Z}_{\ell })}^{i}(M,N^{(\ell
^{n})})\text{.}  \label{e15}
\end{equation}%
Similarly,
\begin{equation*}
\Ext_{\mathbb{Z}{}_{\ell }}^{j}(M,N)\cong
\varprojlim_{n}\Ext_{\mathbb{Z}{}_{\ell }}^{j}(M,N^{(\ell ^{n})}),
\end{equation*}%
and so (Tate 1976, 2.2)%
\begin{equation*}
H_{\text{cts}}^{i}(\Gamma ,\Ext_{\mathbb{Z}{}_{\ell }}^{j}(M,N))\cong \varprojlim%
_{n}H_{\text{cts}}^{i}(\Gamma ,\Ext_{\mathbb{Z}{}_{\ell }}^{j}(M,N^{(\ell
^{n})})).
\end{equation*}%
Therefore, (\ref{e5}) can be obtained by passing to the limit in the
spectral sequences%
\begin{equation*}
H_{\text{cts}}^{r}(\Gamma ,\Ext_{\mathbb{Z}{}_{\ell }}^{s}(M,N^{(\ell
^{n})}))\Longrightarrow \Ext_{\mathcal{R}{}(\mathbb{F}{}_{q};\mathbb{Z}%
_{\ell })}^{r+s}(M,N^{(\ell ^{n})})\text{.}
\end{equation*}
\end{proof}

\subsection{Spectral sequences at $p$}

\begin{lemma}
\label{lc.eb}For torsion-free objects $M,N$ in $\mathcal{R}{}^{+}(\mathbb{F}%
{};\mathbb{Z}{}_{p})$, $\Ext_{\mathcal{R}{}^{+}(\mathbb{F}{};\mathbb{Z}%
{}_{p})}^{2}(M,N)$ is torsion-free.
\end{lemma}

\begin{proof}
From the exact sequences%
\begin{equation*}
0\rightarrow \Ext_{\mathcal{R}{}^{+}}^{1}(M,N)^{(p^{n})}\rightarrow \Ext_{%
\mathcal{R}{}^{+}}^{1}(M,N^{(p^{n})})\rightarrow \Ext_{\mathcal{R}%
{}^{+}}^{2}(M,N)_{p^{n}}\rightarrow 0
\end{equation*}%
we obtain a surjection%
\begin{equation*}
\varinjlim_{n}\Ext_{\mathcal{R}{}^{+}}^{1}(M,N^{(p^{n})})\rightarrow \Ext_{%
\mathcal{R}{}^{+}}^{2}(M,N)_{\text{tors}}\text{.}
\end{equation*}%
Let $\mathcal{R}{}_{p^{n}}^{+}$ denote the full subcategory of $\mathcal{R}%
{}^{+}(\mathbb{F}{};\mathbb{Z}{}_{p})$ of objects killed by $p^{n}$. The
same argument as in the proof pp\pageref{e2zero}--\pageref{e2one} shows that
the map%
\begin{equation*}
s\colon \Ext_{\mathcal{R}{}^{+}}^{1}(M,N^{(p^{n})})\rightarrow \Ext_{%
\mathcal{R}{}_{p^{n}}^{+}}^{1}(M^{(p^{n})},N^{(p^{n})})
\end{equation*}
that replaces each term $E$ in a short exact sequence with $E^{(p^{n})}$ is
an isomorphism. Therefore, it remains to show that
\begin{equation}
\varinjlim_{n}\Ext_{\mathcal{R}{}_{p^{n}}^{+}}^{1}(M^{(p^{n})},N^{(p^{n})})=0.
\label{e31}
\end{equation}%
We claim (i) that the map%
\begin{equation*}
\Ext_{\Crys^{+}(\mathbb{F}{})}^{1}(M\otimes W(\mathbb{F}{}),N\otimes W(%
\mathbb{F}{}))\rightarrow \Ext_{\mathcal{R}%
{}_{p^{n}}^{+}}^{1}(M^{(p^{n})},N^{(p^{n})})
\end{equation*}%
that replaces each term $E$ in a short exact sequence by $E^{(p^{n})}$ is
surjective, and (ii) that $\Ext_{\Crys^{+}(\mathbb{F}{})}^{1}(M,N)$ is a
torsion group. Together, (i) and (ii) imply (\ref{e31})\thinspace , because
(i) gives a surjection
\begin{equation*}
\Ext_{\Crys^{+}(\mathbb{F}{})}^{1}(M\otimes W(\mathbb{F}{}),M\otimes W(%
\mathbb{F}{}))\otimes \mathbb{Q}{}_{p}/\mathbb{Z}{}_{p}\rightarrow \varinjlim_{n}%
\Ext_{\mathcal{R}{}_{p^{n}}^{+}}^{1}(M^{(p^{n})},N^{(p^{n})})
\end{equation*}%
and (ii) implies that the first group is zero.

To prove (i), note that an element of either group splits when regarded as
an extension of $W$-modules, and is therefore determined by the $\sigma $%
-linear map describing the action of the Frobenius. Since every $\sigma $%
-linear map $M^{(p^{n})}\rightarrow N^{(p^{n})}$ lifts to a $\sigma $-linear
map $M\rightarrow N$, (i) is clear. According to a theorem of Manin
(Demazure 1972, p85), the category $\Isoc^{+}(k^{\text{al}})$ is semisimple,
which implies (ii).
\end{proof}

\begin{lemma}
\label{lc.ec}When $k=\mathbb{F}{}$ or $\mathbb{F}{}_{q}$, every object in $%
\mathcal{R}{}^{+}(k;\mathbb{Z}{}_{p})$ is a quotient of a torsion-free
object.
\end{lemma}

\begin{proof}
Suppose that the object $M$ is an extension of objects $M^{\prime }$ and $%
M^{\prime \prime }$ each of which is a quotient of a torsion-free objects,
say, $N^{\prime }\twoheadrightarrow M^{\prime }$, $N^{\prime \prime
}\twoheadrightarrow M^{\prime \prime }$. Pull back the original extension by
$N^{\prime \prime }\rightarrow M^{\prime \prime }$:%
\begin{equation*}
\begin{CD} 0@>>>M'@>>>E@>>>N''@>>>0\\ @.@|@VVV@VVV\\
0@>>>M'@>>>M@>>>M''@>>>0 \end{CD}
\end{equation*}%
If the class of $E$ in $\Ext^{1}(N^{\prime \prime },M^{\prime })$ is the
image of an element of $\Ext^{1}(N^{\prime \prime },N^{\prime })$, then $M$
also is a quotient of an torsion-free object. The obstruction to this is an
element of $\Ext^{2}(N^{\prime \prime },K)$ where $K$ is the kernel of $%
N^{\prime }\rightarrow M^{\prime }$. Certainly, when $k=\mathbb{F}{}_{q}$,
this group will be zero if there exists a resolution%
\begin{equation*}
0\rightarrow A\xr{\cdot\lambda}A\rightarrow N^{\prime \prime }\rightarrow 0
\end{equation*}%
(because in this case Exts in $\mathcal{R}{}_{p}^{+}$ agree with those in $%
\Mod(A)$ --- see \ref{ha0}).

It suffices to prove the lemma in the case $k=\mathbb{F}{}_{q}$, and, for a
given $M$, an argument using restriction of scalars allows us to replace $%
\mathbb{F}{}_{q}$ by a finite extension. After such an extension, a torsion $%
M$ will have a composition series whose quotients are isomorphic to one of $%
(k,\sigma )$ or $(k,0)$. Each of these is a quotient%
\begin{equation*}
\begin{array}{ccccccccc}
0 & \rightarrow & (W,\sigma ) & \xr{p} & (W,\sigma ) & \rightarrow &
(k,\sigma ) & \rightarrow & 0 \\
0 & \rightarrow & (W,p\sigma ) & \xr{p} & (W,p\sigma ) & \rightarrow & (k,0)
& \rightarrow & 0.%
\end{array}%
\end{equation*}%
Since $(W,\sigma )\cong A/A(F-\sigma )$ and $(W,p\sigma )\cong A/A(F-p\sigma
)$ the above remarks prove the lemma for $M$. Let $M$ be an arbitrary
object, and let $M_{1}=M/M_{\text{tors}}$. The obstruction to extending the
statement from $M_{\text{tors}}$ to $M$ is a torsion element of $\Ext%
^{2}(M_{1},K)$, which (see \ref{lc.eb}) becomes zero after a finite
extension of $\mathbb{F}{}_{q}$.
\end{proof}

\begin{proposition}
\label{lc.ed}Let $\Gamma =\Gal(\mathbb{F}{}/\mathbb{F}{}_{q})$. For $M,N$ in
$\mathcal{R}{}^{+}(\mathbb{F}{}_{q};\mathbb{Z}{}_{p})$, there is a spectral
sequence%
\begin{equation}
H^{r}(\Gamma ,\Ext_{\mathcal{R}{}^{+}(\mathbb{F}{};\mathbb{Z}{}_{p})}^{s}(%
\bar{M},\bar{N}))\Longrightarrow \Ext_{\mathcal{R}^{+}(\mathbb{F}{}_{q};%
\mathbb{Z}{}_{p})}^{r+s}(M,N)\text{.}  \label{e24}
\end{equation}
\end{proposition}

\begin{proof}
Lemma \ref{lc.ec} shows that $M$ has a flat resolution. We can not apply
Corollary \ref{ha5} directly, because $\mathcal{R}{}^{+}(\mathbb{F}%
{}_{q^{\prime }};\mathbb{Z}{}_{p})$ does not have internal Homs${}{}$.
However, this is only used in the proof of (a) in \ref{ha3}, namely, in
showing that%
\begin{equation*}
\Hom_{\mathbb{Z}{}_{p}[\Gamma ]}(P,\Hom_{\mathcal{R}{}^{+}(\mathbb{F}%
{}_{q^{\prime }};\mathbb{Z}{}_{p})}(M^{\prime },N^{\prime }))\cong \Hom_{%
\mathcal{R}{}^{+}(\mathbb{F}{}_{q};\mathbb{Z}{}_{p})}(P\otimes M,N)
\end{equation*}%
where $\Gamma =\Gal(\mathbb{F}{}_{q^{\prime }}/\mathbb{F}{}_{q})$ and $P$ is
an object of $\mathcal{R}^{+}{}(\mathbb{F}{}_{q};\mathbb{Z}_{p})$ that
becomes isomorphic to a direct sum of copies of $(W,\sigma )$ in $\mathcal{R}%
{}^{+}(\mathbb{F}{}_{q^{\prime }};\mathbb{Z}{}_{p})$. This can be shown
directly.
\end{proof}

\begin{proposition}
\label{lc.em}Let $\Gamma =\Gal(\mathbb{F}{}/\mathbb{F}{}_{q})$. For $M,N$ in
$\mathcal{R}{}^{+}(\mathbb{F}{}_{q};\mathbb{Z}{}_{p})$, there is a spectral
sequence%
\begin{equation}
H_{\text{cts}}^{r}(\Gamma ,\Ext_{\Crys^{+}(\mathbb{F}{})}^{s}(M_{W(\mathbb{F}%
{})},N_{W(\mathbb{F}{})}))\Longrightarrow \Ext_{\mathcal{R}^{+}(\mathbb{F}%
{}_{q};\mathbb{Z}{}_{p})}^{r+s}(M,N)\text{.}  \label{e17}
\end{equation}%
Here $H_{\text{cts}}^{r}$ is computed using continuous cochains relative to
the $p$-adic topology on $\Ext_{W(\mathbb{F}{})}^{s}(M,N)$.
\end{proposition}

\begin{proof}
When $N$ (or $M$) is finite, (\ref{e17}) becomes the (\ref{e24}) --- see the
proof of Proposition \ref{lc.e}. The remainder of the proof of Proposition %
\ref{lc.e} applies to this case, except that $N^{(p^{n})}$ and $\Ext_{\Crys%
^{+}(\mathbb{F}{})}^{j}(M,N^{(p^{n})})$ need no longer be finite. However, $%
N^{(p^{n})}$ is of finite length over $W$, and so we can apply
Jensen 1972, 7.2, to see that $\varprojlim{}^{i}N^{(p^{n})}=0$ for
$i>0$. Moreover, there
exists an $N$ such that $p^{N}\cdot \Ext_{\smallskip \noindent \Crys^{+}(%
\mathbb{F}{})}^{j}(M,N^{(p^{n})})$ is finite for all $n$, and so the inverse
system $\Ext_{\noindent \Crys^{+}(\mathbb{F}{})}^{j}(M,N^{(p^{n})})$
satisfies the Mittag-Leffler condition$.$
\end{proof}

\section{Homological Algebra in the Category of Motives}

In this section, we assume that the $l$-adic Tate conjectures hold for $%
\mathcal{M}{}^{+}(k;\mathbb{Q}{})$ (see p\pageref{e27}).

\subsection{Homological algebra in $\mathcal{M}^{+}(k{};\mathbb{\mathbb{Z}})$%
}

\begin{theorem}
\label{ha6}Let $X$ and $Y$ be effective motives in $\mathcal{M}^{+}(k;%
\mathbb{Z})$.

\begin{enumerate}
\item The map $(\alpha ,\beta )\mapsto \alpha -\beta $,
\begin{equation*}
\Hom_{\mathcal{M}{}^{+}(k;\mathbb{Q)}}(X_{0},Y_{0})\times \Hom_{\mathcal{R}%
^{+}(k;\mathbb{\hat{Z}})}(X_{f},Y_{f})\rightarrow \Hom_{\mathcal{R}^{+}(k;%
\mathbb{A}_{f})}(\omega _{f}(X_{0}),\omega _{f}(Y_{0})),
\end{equation*}%
is surjective with kernel $\Hom_{\mathcal{M}{}^{+}(k;\mathbb{Z}{})}(X,Y)$.

\item For $i>0$, $\Ext_{\mathcal{M}{}^{+}(k;\mathbb{Z})}^{i}(X,Y)$ is
torsion and the sequence
\begin{equation}
0\rightarrow \Ext_{\mathcal{M}{}^{+}(k;\mathbb{Z}{})}^{i}(X,Y)_{\mathbb{Z}%
_{l}}\rightarrow \Ext_{\mathcal{R}(k;\mathbb{Z}_{l})}^{i}(X_{l},Y_{l})%
\rightarrow T_{l}\Ext_{\mathcal{M}{}^{+}(k;\mathbb{Z}{})}^{i+1}(X,Y)%
\rightarrow 0  \label{e19}
\end{equation}%
is exact (all $l$); hence,
\begin{equation*}
\Ext_{\mathcal{M}{}^{+}(k;\mathbb{\mathbb{Z}{}})}^{i}(X,Y)\cong \Ext_{%
\mathcal{R}^{+}(k;\mathbb{\hat{Z}}{})}^{i}(X_{f},Y_{f})_{\text{tors}}.
\end{equation*}

\item If $k$ is separably closed, then $\Ext_{\mathcal{M}{}^{+}(k;\mathbb{Z}%
{})}^{i}(X,Y)=0$ for all $i>2$, and for all $i>1$ when $X$ and $Y$ are
torsion-free.
\end{enumerate}
\end{theorem}

\begin{proof}
For any $\mathbb{Z}{}$-module $M$ of finite rank, the sequence%
\begin{equation*}
0\rightarrow M\rightarrow M_{\mathbb{\hat{Z}}{}}\times M_{\mathbb{Q}%
{}}\rightarrow M_{\mathbb{A}_{f}}\rightarrow 0
\end{equation*}%
is exact. We obtain (a) by taking $M=\Hom(X,Y)$ in this sequence and
applying (\ref{cm10}, \ref{cm11}, \ref{cm13}).

The remainder of the proof will require several steps.

\medskip\noindent$\Ext_{\mathcal{M}{}^{+}{}(k;\mathbb{Z}{})}^{i}(X,Y)$
\textsc{is torsion for }$i>0$\textsc{: }It suffices to prove that $\Ext_{%
\mathcal{M}{}^{+}{}(k;\mathbb{Z}{})}^{1}(X,Y)$ is torsion. Moreover, we may
assume that $X$ and $Y$ are torsion-free. Consider an exact sequence
\begin{equation*}
0\rightarrow Y\rightarrow E\overset{\pi }{\rightarrow }X\rightarrow 0
\end{equation*}%
in $\mathcal{M}{}^{+}(k;\mathbb{Z}{})$. Because $\mathcal{M}{}^{+}(k;\mathbb{%
Q}{})$ is semisimple, there exists a morphism $\alpha _{0}\colon
X_{0}\rightarrow E_{0}$ such that $\pi _{0}\circ \alpha _{0}=\id_{X_{0}}$.
For some $m\in \mathbb{Z}{}$, $m\alpha _{0}\in \Hom_{\mathcal{M}{}^{+}(k;%
\mathbb{Z})}(X,E)$ (see \ref{cm9}). Now the equation $\pi \circ (m\alpha
_{0})=m\id_{X}$ implies that the class of the sequence is killed by $m$.

\smallskip Write $f^{i}(X,Y)$ for the natural map%
\begin{equation*}
\Ext_{\mathcal{M}{}^{+}(k;\mathbb{\mathbb{Z}{}})}^{i}(X,Y)_{\mathbb{Z}%
{}_{l}}\rightarrow \Ext_{\mathcal{R}^{+}(k;\mathbb{Z}_{l})}^{i}(X_{l},Y_{l})%
\text{.}
\end{equation*}%
Recall (\ref{cm13}) that $f^{0}(X,Y)$ is an isomorphism.

\medskip \noindent \textsc{The maps }$f^{i}(X,Y)$ \textsc{are isomorphisms
when }$X$ \textsc{and }$Y$ \textsc{are finite. }The functor $X\mapsto X_{f}$
defines an equivalence of the category of torsion effective motives with the
category of torsion objects in $\mathcal{R}^{+}{}(k;\mathbb{\hat{Z})}$, and
so it suffices to prove that the Exts${}$ in the torsion subcategories
coincide with the Exts${}$ in the full categories\label{p1}. We do this for $%
\mathcal{M}{}^{+}(k;\mathbb{Z}{})$ since $\mathcal{R}{}^{+}(k;\mathbb{Z}%
{}_{l})$ is similar. Because $\mathcal{M}{}^{+}(k;\mathbb{Z}{})$ is
noetherian, $\Ind(\mathcal{M}{}^{+}(k;\mathbb{Z}{}))$ has enough injectives,
and the Exts${}$ computed using injective resolutions coincide with the
Yoneda Exts${}$ (\ref{ha0}). Because $n\colon \1\rightarrow \1$ is
injective, every injective object in $\Ind(\mathcal{M}{}^{+}(k;\mathbb{Z}%
{})){}$ is divisible, and so an injective resolution $Y\rightarrow
I^{\bullet }$ in $\Ind(\mathcal{M}{}^{+}(k;\mathbb{Z}{})){}$ gives an
injective resolution $Y\rightarrow I_{\text{tors}}^{\bullet }$ in the
subcategory of torsion objects. As $\Hom(X,I_{\text{tors}}^{\bullet })\cong %
\Hom(X,I^{\bullet })$, this proves that the Exts${}$ in the two categories
agree.

\medskip \noindent \textsc{The maps }$f^{i}(X,Y)$ \textsc{are isomorphisms
when }$Y$\textsc{\ is finite.} It suffices to prove this when $X$ is
torsion-free and $Y$ is $l$-torsion, say $l^{n}Y=0$. The $\Ext$ sequences of
\begin{equation*}
0\rightarrow X\xr{l^n}X\rightarrow X/l^{n}X\rightarrow 0
\end{equation*}%
give a diagram%
\begin{equation*}
\begin{CD} 0 @>>> \Ext^{i-1}_{\mathcal{M}^+}(X,Y) @>>>
\Ext_{\mathcal{M}^+}^{i}(X/l^{n}X,Y) @>>> \Ext_{\mathcal{M}^+}^{i}(X,Y) @>>>
0\\ @. @VVf^{i-1}(X,Y)V @VV{\cong}V @VVf^i(X,Y)V @.\\ 0 @>>>
\Ext^{i-1}_{\mathcal{R}^+}(X_{l},Y_{l}) @>>>
\Ext_{\mathcal{R}^+}^{i}(X_{l}/l^{n}X_{l},Y_{l}) @>>>
\Ext_{\mathcal{R}^+}^{i}(X_{l},Y_{l}) @>>> 0 \end{CD}
\end{equation*}%
from which it follows that $f^{i-1}(X,Y)$ is injective and $f^{i}(X,Y)$ is
surjective. Since this holds for all $i$, $f^{i}(X,Y)$ is an isomorphism.

\medskip \noindent \textsc{The maps }$f^{i}(X,Y)$ \textsc{are isomorphisms
when }$X$\textsc{\ is finite.} Similar to the preceding case.

\medskip \noindent \textsc{The sequence (\ref{e19}) is exact. }Because $%
f^{i}(X,Y)$ is an isomorphism when $Y$ is finite, $f^{i}(X,Y^{\bullet })$
(map of hyper-Exts${}{}$) is an isomorphism when $Y^{\bullet }$ is a bounded
complex with finite cohomology, for example, when $Y^{\bullet
}=Y^{\{l^{n}\}}=_{\text{df}}Y\xr{l^n}Y$:%
\begin{equation*}
\Ext_{\mathcal{M}^{+}{}}^{i}(X,Y^{\{l^{n}\}})\cong \Ext_{\mathcal{R}%
{}^{+}}^{i}(X_{l},Y_{l}^{\{l^{n}\}})\text{.}
\end{equation*}%
There are exact sequences%
\begin{equation}
0\rightarrow \Ext_{\mathcal{M}^{+}{}}^{i}(X,Y)^{(l^{n})}\rightarrow \Ext_{%
\mathcal{M}{}^{+}}^{i}(X,Y^{\{l^{n}\}})\rightarrow \Ext_{\mathcal{M}%
{}^{+}}^{i}(X,Y)_{l^{n}}\rightarrow 0\text{,}  \label{e20}
\end{equation}%
and%
\begin{equation*}
\varprojlim{}\Ext_{\mathcal{M}{}^{+}}^{i}(X,Y^{\{l^{n}\}})\cong \varprojlim{}\Ext_{%
\mathcal{R}{}^{+}}^{i}(X_{l},Y_{l}^{\{l^{n}\}})\cong \varprojlim{}\Ext_{\mathcal{R}%
^{+}{}}^{i}(X_{l},Y_{l}^{(l^{n})})\cong \Ext_{\mathcal{R}%
^{+}{}}^{i}(X_{l},Y_{l})
\end{equation*}%
(for the third isomorphism, see the proofs of \ref{lc.e} and \ref{lc.em}).
Thus, on passing to the inverse limit in (\ref{e20}), we obtain (\ref{e19}).

\smallskip This completes the proof of (b). We now assume that $k=k^{\text{%
sep}}$.

\medskip \noindent $\Ext_{\mathcal{M}{}^{+}}^{2}(X,Y)=0$\emph{\ }\textsc{if }%
$X$\textsc{\ and }$Y$\textsc{\ are torsion-free.} \label{e2zero}For an
integer $m$, let $\mathcal{M}{}_{m}^{+}$ (resp. $\mathcal{R}{}_{m}^{+})$
denote the subcategory of $\mathcal{M}{}^{+}(k;\mathbb{Z}{})$ (resp. $%
\mathcal{R}^{+}(k{};\mathbb{\hat{Z}})$) of objects killed by $m$. On passing
to the direct limit in the sequence
\begin{equation*}
0\rightarrow \Ext_{\mathcal{M}{}^{+}}^{1}(X,Y)^{(m)}\rightarrow \Ext_{%
\mathcal{M}{}^{+}}^{1}(X,Y^{(m)})\rightarrow \Ext_{\mathcal{M}%
{}^{+}}^{2}(X,Y)_{m}\rightarrow 0
\end{equation*}%
and using that $\Ext_{\mathcal{M}{}^{+}}^{1}(X,Y)$ and $\Ext_{\mathcal{M}%
{}^{+}}^{2}(X,Y)$ are torsion, we find that
\begin{equation*}
\varinjlim_{m}\Ext_{\mathcal{M}{}^{+}}^{1}(X,Y^{(m)})\cong \Ext_{\mathcal{M}%
{}^{+}}^{2}(X,Y).
\end{equation*}%
Consider the diagram%
\begin{equation*}
\begin{diagram} \Ext_{\mathcal{M}^+}^{1}(X^{(m)},Y^{(m)}) & &\rTo^{\pi }&
&\Ext_{\mathcal{M}^+}^{1}(X,Y^{(m)}) \\ &\luTo^i & &\ldTo^s& \\
&&\Ext_{\mathcal{M}^+_{m}}^{1}(X^{(m)},Y^{(m)}) &&\\ \end{diagram}
\end{equation*}%
in which $\pi $ is the surjection induced by $X\rightarrow X^{(m)}$, $i$ is
the obvious injection, and $s$ applied to a short exact sequence replaces
each motive $E$ in the sequence with $E^{(m)}$. One checks easily that $%
s\circ (\pi \circ i)=\id$ and $(\pi \circ i)\circ s=\id$, and so $s$ is an
isomorphism. Thus%
\begin{equation*}
\Ext_{\mathcal{M}{}^{+}}^{2}(X,Y)\cong \varinjlim_{m}\Ext_{\mathcal{M}%
{}_{m}^{+}}^{1}(X^{(m)},Y^{(m)})\cong \varinjlim_{m}\Ext_{\mathcal{\mathcal{R}{}}%
{}_{m}^{+}}^{1}(X_{f}^{(m)},Y_{f}^{(m)})\text{.}
\end{equation*}%
For $l\neq p$, $\Ext_{\mathcal{\mathcal{R}{}}%
{}_{l^{n}}^{+}}^{1}(X_{f}^{(l^{n})},Y_{f}^{(l^{n})})$ is obviously zero, and
so it remains to show that
\begin{equation}
\varinjlim_{p^{n}}\Ext_{\mathcal{\mathcal{R}{}}%
{}_{p^{n}}^{+}}^{1}(X_{p}^{(p^{n})},Y_{p}^{(p^{n})})=0\text{.}  \label{e25}
\end{equation}%
But this is shown in the proof of (\ref{lc.eb}) (for $k=\mathbb{F}{}$, but
the argument works generally).

\medskip \noindent $\Ext_{\mathcal{M}{}^{+}}^{i}(X,Y)=0$\textsc{\ for all }$%
i>2$\textsc{.} This follows from (b) and the fact that $\Ext_{\mathcal{R}%
{}^{+}(k;\mathbb{Z}{}_{l})}^{i}(M,N)=0$ for $i>2$ (for $l\neq p$, this is
obvious, and for $l=p$ it can be proved by the arguments used to prove the $%
p $-case of Theorem \ref{lc.a} below).\label{e2one}
\end{proof}

\begin{remark}
\label{ha7}The group $\Ext_{\mathcal{R}{}^{+}(\mathbb{F}{};\mathbb{\mathbb{Z}%
}_{\ell })}^{1}(M,N)$ need not be torsion. Consider, for example, the
extension
\begin{equation*}
0\rightarrow \mathbb{Z}_{\ell }\rightarrow \mathbb{Z}_{\ell }\oplus \mathbb{Z%
}_{\ell }\rightarrow \mathbb{Z}_{\ell }\rightarrow 0
\end{equation*}%
\noindent on which the Frobenius element in $\Gamma =\Gal(\mathbb{F}{}/%
\mathbb{F}{}_{p})$ acts on the middle term as the matrix $\left(
\begin{smallmatrix}
1 & a \\
0 & 1%
\end{smallmatrix}%
\right) $ with $a$ a nonzero element of $\mathbb{Z}{}_{\ell }$. Multiplying
the class of the extension by $m$ corresponds to multiplying $a$ by $m$, and
so this represents a nontorsion element of $\Ext_{\mathcal{R}{}^{+}(\mathbb{F%
}{}_{p};\mathbb{\mathbb{Z}}_{\ell })}^{1}(M,N)$, which remains nontorsion in
$\Ext_{\mathcal{R}{}^{+}(\mathbb{F}{};\mathbb{\mathbb{Z}}_{\ell })}^{1}(M,N)$%
. Thus, in general,
\begin{equation*}
\Ext_{\mathcal{M}{}^{+}(\mathbb{F}{};\mathbb{\mathbb{Z}{}})}^{1}(X,Y)(\ell
)\not\approx \Ext_{\mathcal{R}(\mathbb{F}{};\mathbb{Z}{}_{\ell
})}^{1}(X_{\ell },Y_{\ell }).
\end{equation*}%
Similarly, the group $\Ext_{\mathcal{R}{}^{+}(\mathbb{F}{};\mathbb{Z}%
{}_{p})}^{1}(M,N)$ need not be torsion (because the category $\mathcal{R}%
{}^{+}(\mathbb{F}{};\mathbb{Z}{}_{p})$, in contrast to $\Isoc^{+}(\mathbb{F}%
{})$, is not semisimple).
\end{remark}

\subsection{Homological algebra in $\mathcal{M}^{+}{}(\mathbb{F}{}_{q};%
\mathbb{Z}{})$}

\begin{lemma}
\label{ha9a}Every effective motive over $\mathbb{F}{}_{q}$ or $\mathbb{F}{}$
is a quotient of a torsion-free effective motive.
\end{lemma}

\begin{proof}
Straightforward adaptation of the proof of Lemma \ref{lc.ec}.
\end{proof}

\begin{theorem}
\label{ha9}Let $\Gamma =\Gal(\mathbb{F}{}/\mathbb{F}{}_{q})$. For all
motives $X$ and $Y$ over $\mathbb{F}_{q}$, there is a spectral sequence:%
\begin{equation*}
H^{r}(\Gamma ,\Ext_{\mathcal{M}^{+}(\mathbb{F}{};\mathbb{Z}{})}^{s}(X_{/%
\mathbb{F}{}},Y_{/\mathbb{F}{}}))\implies \Ext_{\mathcal{M}^{+}(\mathbb{F}%
{}_{q};\mathbb{Z}{})}^{r+s}(X,Y)\text{.}
\end{equation*}
\end{theorem}

\begin{proof}
Lemma \ref{ha9a} shows that $X$ has a flat resolution. We can not apply
Corollary \ref{ha5} directly because $\mathcal{M}{}^{+}(\mathbb{F}{}_{q};%
\mathbb{Z}{})$ does not have internal Homs, but its proof can be adapted as
in the proof of \ref{lc.ed}.
\end{proof}

\subsection{Homological algebra in $\mathcal{M}{}(k{};\mathbb{Z}{})$}

The results of the last two subsections hold mutatis mutandis for motives.

\subsection{Examples}

\begin{example}
\label{ha8}We apply (\ref{ha6}) to compute $\Ext_{\mathcal{M}^{+}(\mathbb{F}%
{};\mathbb{Z}{})}^{1}(\mathbb{Z}{},\mathbb{Z}{}(r))$ (we are writing $%
\mathbb{Z}{}$ for $\1$ and $\mathbb{Z}{}(r)$ for $\mathbb{T}^{\otimes r}$).

For $r=0$, $\Hom_{\mathcal{M}^{+}(\mathbb{F};\mathbb{Z}{})}(\mathbb{Z}{},%
\mathbb{Z}{}(0))\cong \mathbb{Z}{}$ and $\Ext_{\mathcal{M}^{+}(\mathbb{F}{};%
\mathbb{Z}{})}^{i}(\mathbb{Z}{},\mathbb{Z}{}(0))=0$ for $i=1,2$..

For $r\neq 0$, $\Hom_{\mathcal{M}^{+}(\mathbb{F}{};\mathbb{Z}{})}(\mathbb{Z}%
{},\mathbb{Z}{}(r))=0$, and so
\begin{equation*}
\Hom_{\mathcal{M}^{+}(\mathbb{F}{};\mathbb{Z}{})}(\mathbb{Z}{},\mathbb{Z}%
{}(r)^{(m)})\overset{\cong }{\rightarrow }\Ext_{\mathcal{M}^{+}(\mathbb{F}{};%
\mathbb{Z}{})}^{1}(\mathbb{Z}{},\mathbb{Z}{}(r))_{m}.
\end{equation*}%
For $m=p^{n}$, the left hand side is zero, and for $m=\ell ^{n}$ it is $(%
\mathbb{Z}{}^{(\ell ^{n})}{})(r)$. \noindent Therefore, on passing to the
direct limit over $m$, we find that
\begin{equation*}
\Ext_{\mathcal{M}^{+}(\mathbb{F}{};\mathbb{Z}{})}^{1}(\mathbb{Z}{},\mathbb{Z}%
{}(r))\cong \left( \bigoplus_{\ell \neq p}\mathbb{Q}_{\ell }/\mathbb{Z}%
{}_{\ell }\right) (r)
\end{equation*}%
where $\mathbb{(Q}_{\ell }/\mathbb{Z}{}_{\ell })(r)$ is $\mathbb{Q}{}_{\ell
}/\mathbb{Z}{}_{\ell }$ with the Frobenius element of $\Gal(\mathbb{F}{}/%
\mathbb{F}{}_{q})$ acting as $q^{r}$. Moreover, $\Ext_{\mathcal{M}^{+}(%
\mathbb{F},\mathbb{Z})}^{i}(\mathbb{Z}{},\mathbb{Z}{}(r))=0$ for $i\neq 0,2$%
. (Without the $+$, these statement holds only modulo $p$-torsion.)
\end{example}

\begin{example}
\label{ha10}We compute $\Ext_{\mathcal{M}^{+}(\mathbb{F}{}_{q},\mathbb{Z}%
{})}^{1}(\mathbb{Z},\mathbb{Z}{}(r))$ for $r>0$. \noindent\ The spectral
sequence in (\ref{ha9}) shows that $\Ext_{\mathcal{M}^{+}(\mathbb{F}{}_{q},%
\mathbb{Z}{})}^{1}(\mathbb{Z},\mathbb{Z}{}(r))$ and $\Ext_{\mathcal{M}^{+}(%
\mathbb{F}{}_{q},\mathbb{Z}{})}^{2}(\mathbb{Z},\mathbb{Z}{}(r))$ are the
kernel and cokernel respectively of
\begin{equation*}
\oplus _{\ell \neq p}\mathbb{Q}_{\ell }/\mathbb{\mathbb{Z}}_{\ell }%
\xr{x\mapsto (q^{r}-1)x}\oplus _{\ell \neq p}\mathbb{Q}_{\ell }/\mathbb{%
\mathbb{Z}}_{\ell }.
\end{equation*}%
Therefore $\Ext_{\mathcal{M}^{+}(\mathbb{F}{}_{q};\mathbb{Z}{})}^{1}(\mathbb{%
Z}{},\mathbb{Z}{}(r))$ is cyclic of order $q^{r}-1$, and $\Ext_{\mathcal{M}%
^{+}(\mathbb{F}{}_{q};\mathbb{Z}{})}^{i}(\mathbb{Z}{},\mathbb{Z}{}(r))=0$
for $i\neq 0,2$.
\end{example}

\begin{example}
\label{ha10m}Let $A$ be an abelian variety over $\mathbb{F}{}$, and let $%
h_{1}A$ be the isomotive $h_{1}(A)_{0}$ endowed with the $\mathbb{Z}{}$%
-structure provided by the maps $H_{1}(A,\mathbb{Z}{}_{l})\rightarrow
H_{1}(A,\mathbb{Q}{}_{l})$. Since $\Hom(\mathbb{Z}{},h_{1}A)=0$,\qquad
\begin{equation*}
\Hom_{\mathcal{M}{}^{+}(\mathbb{F}{};\mathbb{Z}{})}(\mathbb{Z}%
{},(h_{1}A)^{(m)})\cong \Ext_{\mathcal{M}{}^{+}(\mathbb{F}{};\mathbb{Z}%
{}){}}^{1}(\mathbb{Z}{},h_{1}A)_{m}
\end{equation*}%
for all $m$, and so%
\begin{equation*}
\Hom_{\mathcal{M}{}^{+}(\mathbb{F}{};\mathbb{Z}{}){}}(\mathbb{Z}%
{},h_{1}A\otimes \mathbb{Q}{}/\mathbb{Z}{})\cong \Ext_{\mathcal{M}{}^{+}(%
\mathbb{F}{};\mathbb{Z}{}){}}^{1}(\mathbb{Z}{},h_{1}A)\text{.}
\end{equation*}%
But%
\begin{equation*}
\Hom_{\mathcal{M}{}^{+}(\mathbb{F}{};\mathbb{Z}{}){}}(\mathbb{Z}%
{},h_{1}A\otimes \mathbb{Q}{}/\mathbb{Z}{})\cong \Hom_{\mathcal{\mathcal{R}{}%
}^{+}(\mathbb{F}{};\mathbb{\hat{Z}}{}){}}(\mathbb{Z}{},(h_{1}A)_{f\mathbb{Q}%
{}}/(h_{1}A)_{f})\cong A(\mathbb{F}{})\text{.}
\end{equation*}%
Thus,%
\begin{equation*}
\Ext_{\mathcal{M}{}^{+}(\mathbb{F}{};\mathbb{Z}{}){}}^{1}(\mathbb{Z}%
{},h_{1}A)\cong A(\mathbb{F}{}).
\end{equation*}%
When $A$ is defined over $\mathbb{F}{}_{q}$, $\Gal(\mathbb{F}{}/\mathbb{F}%
_{q})$ acts, and, on taking invariants, we find that%
\begin{equation*}
\Ext_{\mathcal{M}{}^{+}(\mathbb{F}{}_{q};\mathbb{Z}{})}^{1}(\mathbb{Z}%
{},h_{1}A)\cong A(\mathbb{F}{}_{q}).
\end{equation*}

A similar argument using that $\Hom_{\mathcal{M}{}^{+}(\mathbb{F}{};\mathbb{Z%
}{})}((h_{1}A)(-1),\mathbb{Z}{})=0$ shows that%
\begin{equation*}
\Hom_{\mathcal{M}{}^{+}(\mathbb{F}{};\mathbb{Z}{}){}}((h_{1}A)(-1),\mathbb{Q}%
{}/\mathbb{Z}{})\cong \Ext_{\mathcal{M}{}^{+}(\mathbb{F}{};\mathbb{Z}%
{}){}}^{1}((h_{1}A)(-1),\mathbb{Z}{})\text{.}
\end{equation*}%
But
\begin{equation*}
\Hom_{\mathcal{M}{}^{+}(\mathbb{F}{};\mathbb{Z}{}){}}((h_{1}A)(-1),\mathbb{Q}%
{}/\mathbb{Z}{})\cong \Hom_{\mathcal{R}^{+}(\mathbb{F}{};\mathbb{\hat{Z}}%
{}){}}((h_{1}A)_{f}(-1),\mathbb{Q}{}/\mathbb{Z}{})\cong A^{\vee }(\mathbb{F}%
{})
\end{equation*}%
where $A^{\vee }$ is the dual abelian variety. For an abelian variety $A$
over $\mathbb{F}{}_{q}$, this becomes%
\begin{equation*}
\Ext_{\mathcal{M}{}^{+}(\mathbb{F}{}_{q};\mathbb{Z}{})}^{1}((h_{1}A)(-1),%
\mathbb{Z}{})\cong A^{\vee }(\mathbb{F}{}_{q}).
\end{equation*}
\end{example}

\begin{remark}
\label{ha11}It is generally conjectured, that for a smooth variety $V$ over
a field $F$,%
\begin{equation}
H_{\text{mot}}^{i}(V,\mathbb{Q}{}(r))\cong (K_{2r-i}V)_{\mathbb{Q}{}}^{(r)}
\label{e12}
\end{equation}%
\noindent\ \noindent and
\begin{equation}
H_{\text{mot}}^{i}(V,\mathbb{Q}{}(r))\cong \Ext_{\mathcal{MM}(F;\mathbb{Q}%
{})}^{i}(\1,(hV)(r));  \label{e13}
\end{equation}%
\noindent see, for example, Deligne 1994, 3.7, 3.3. Here $\Ext_{\mathcal{M%
\mathcal{M}{}}(F;\mathbb{Q}{})}^{i}(\1,(hV)(r))$ is computed in the
conjectural category of mixed isomotives over $F$ and $(K_{2r-i}V)_{\mathbb{Q%
}{}}^{(r)}$ is the subspace of $(K_{2r-1}V)_{\mathbb{Q}{}}{}$ on which the
Adams operators $\psi ^{k}$ act as $k^{r}$.

In particular, for $F$ itself, it is conjectured that
\begin{equation}
\Ext_{\mathcal{M}{}\mathcal{M}{}(F;\mathbb{Q}{})}^{1}(\mathbb{Q}{},\mathbb{Q}%
{}(r))\cong (K_{2r-1}F)_{\mathbb{Q}{}}^{(r)}  \label{e14}
\end{equation}%
For $F=\mathbb{F}{}_{q}$ or $\mathbb{F}{}$, the statement (\ref{e14}) is
true because the $K$-groups are torsion and, under any definition, the
category of mixed isomotives will be semisimple, so both sides are $0$.

With $\mathbb{Z}{}$-coefficients, statement (\ref{e14}) is expected to be
false in general, and (\ref{e12}) is definitely false\footnote{%
For a projective smooth variety over a field $k$,
\begin{equation*}
H_{\text{mot}}^{3}(V,\mathbb{Z}{}(1))=H^{2}(V,\mathbb{G}_{m}),
\end{equation*}%
which is a nonzero torsion group in general whereas $K_{-1}V=0$. See also
Bloch and Esnault 1996, p.~ 304, which explains why a closely related
conjecture should not be true integrally.}. However, Quillen (1972) shows
that $K_{2r-1}\mathbb{F}{}\approx \oplus _{\ell \neq p}\mathbb{Q}_{\ell }/%
\mathbb{\mathbb{Z}}_{\ell }$ with the Frobenius automorphism acting as $%
p^{r} $ and $K_{2r-1}\mathbb{F}{}_{q}=(K_{2r-1}\mathbb{F}{})^{\Gal(\mathbb{F}%
{}/\mathbb{F}{}_{q})}$, and Hiller (1981) and Kratzer (1980) show that the
Adams operator $\psi ^{k}$ acts on $K_{2r-1}\mathbb{F}{}$ as multiplication
by $k^{r}$, and so an integral version of (\ref{e14}) \emph{does }hold when $%
F=\mathbb{F}{}_{q}$ or $\mathbb{F}{}$:
\begin{equation}
\Ext_{\mathcal{M}{}^{+}(F;\mathbb{Z}{})}^{1}(\mathbb{Z}{},\mathbb{Z}%
{}(r))\approx (K_{2r-1}F)^{(r)}.  \label{e26}
\end{equation}

The isomorphism (\ref{e26}) can be made canonical: Soul\'{e} (1979, IV.2,
p.~ 284) proves that the Chern class map $c_{i,1}\colon K_{2i-1}(\mathbb{F}%
{}_{q})_{\mathbb{\mathbb{Z}{}}_{l}}\rightarrow H^{1}(\mathbb{F}{}_{q},%
\mathbb{Z}{}_{l}(i))$ is an isomorphism if $i<l$; later Dyer and Friedlander
(1985) defined a Chern character $ch_{i,1}\colon K_{2i-1}(\mathbb{F}{}_{q})_{%
\mathbb{\mathbb{Z}{}}_{l}}\rightarrow H^{1}(\mathbb{F}{}_{q},\mathbb{Z}%
{}_{l}(i))$ that is related the Chern class map by $%
(-1)^{i-1}(i-1)!ch_{i,1}=c_{i,1}$ (Soul\'{e} 1999, 3.2.2, p.~ 267).
\end{remark}

\section{Extensions and Zeta Functions: Local Case}

\label{localext}We use the notations: $\Gamma =\Gal(\mathbb{F}{}/\mathbb{F}%
{}_{q})$, $q=p^{a}$, $W=W(\mathbb{F}{}_{q})$, $B=B(\mathbb{F}_{q})$, $%
A=W[F,\sigma ]$, $\mathcal{R}{}_{\ell }=\mathcal{R}{}(\mathbb{F}{}_{q};%
\mathbb{Z}{}_{\ell })=_{\text{df}}\mathcal{R}{}(\Gamma ;\mathbb{Z}{}_{\ell
}) $, $\mathcal{\bar{R}}_{\ell }=\Mod(\mathbb{Z}{}_{\ell }){}$, $\mathcal{R}%
{}_{p}=\mathcal{R}^{+}{}(\mathbb{F}{}_{q},\mathbb{Z}{}_{p})=_{\text{df}}\Crys%
^{+}(\mathbb{F}{}_{q}),$ $\mathcal{\bar{R}}_{p}=\Crys^{+}(\mathbb{F}{})$.${}$

\subsection{Statement of the theorem}

Let $M$ be an object of $\mathcal{R}{}_{l}$. The Frobenius endomorphism $\pi
_{M}$ of $M$ is the map induced by the element $x\mapsto x^{q}$ of $\Gamma $
when $l\neq p$ and by $F^{a}$ when $l=p$. We denote the characteristic and
minimum polynomials of $\pi _{M}$ acting on $M_{\mathbb{Q}{}_{l}}$ by $%
m_{M}(t)$ and $P_{M}(t)$ respectively. They are monic polynomials with
coefficients in $\mathbb{Z}{}_{l}$ (see Demazure 1972, p89, for the case $%
l=p $). For $M$ in $\mathcal{R}{}_{p}$, let $r(M)$ denote the rank of $M$
and $s(M)$ the sum of the slopes of $M$; thus if
\begin{equation*}
P_{M}(T)=t^{h}+\cdots +c
\end{equation*}%
then $r(M)=h$ and $s(M)=\ord_{p}(c)/a$ (cf. Demazure 1972, p90). For $M$ in $%
\mathcal{R}{}_{\ell }$ we set $s(M)=0$.

Let $M,N$ be objects of $\mathcal{R}{}_{l}$. From the spectral sequences (%
\ref{lc.e}, \ref{lc.em}) we obtain a map%
\begin{equation*}
\Hom_{\mathcal{\bar{R}}{}_{l}}(\bar{M},\bar{N})_{\Gamma }\rightarrow \Ext_{%
\mathcal{R}{}_{l}}^{1}(M,N)\text{.}
\end{equation*}%
Define\label{p2}%
\begin{equation*}
f=f(M,N)\colon \Hom_{\mathcal{R}_{l}}(M,N)\rightarrow \Ext_{\mathcal{R}%
_{l}}^{1}(M,N)
\end{equation*}%
to be the composite of this with the obvious map $\Hom_{\mathcal{R}%
{}_{l}}(M,N)\rightarrow \Hom_{\mathcal{\bar{R}}{}_{l}}(\bar{M},\bar{N}%
)_{\Gamma }$.

\begin{theorem}
\label{lc.a}If $m_{M}(T)$ and $m_{N}(T)$ have no multiple root in common,
then $z(f)$ is defined and satisfies%
\begin{equation}
z(f)\cdot \lbrack \Ext_{\mathcal{R}{}_{l}}^{2}(M,N)]=\left\vert q^{s(M)\cdot
r(N)}\cdot \prod_{a_{i}\neq b_{j}}\left( 1-\frac{b_{j}}{a_{i}}\right)
\right\vert _{l}\text{.}  \label{e7}
\end{equation}%
where $(a_{i})_{1\leq i\leq r(M)}$ and $(b_{j})_{1\leq j\leq r(N)}$ are the
roots of $P_{M}(t)$ and $P_{N}(t)$ respectively.
\end{theorem}

\begin{remark}
\label{lc.c}

\begin{enumerate}
\item When $l\neq p$, the term $q^{s(M)\cdot r(N)}=1$.

\item In the course of the proof, we shall show that $\Ext_{\mathcal{R}%
{}_{l}}^{i}(M,N)=0$ for $i>2$ (without any condition on $M$ and $N$) and
that $\Ext_{\mathcal{R}{}_{l}}^{2}(M,N)$ is finite.

\item When one of $M$ or $N$ is finite, the theorem simply states that the
groups $\Ext_{\mathcal{R}{}_{l}}^{i}(M,N)$ are finite and the alternating
product of their orders is $1$.

\item Let $P_{M^{\vee }\otimes N}(t)=\prod \left( 1-\frac{b_{j}}{a_{i}}%
t\right) $, and let $\rho (M,N)$ be the rank of $\Hom_{\mathcal{R}%
{}_{l}}(M,N)$. Then $\rho (M,N)$ is the number of pairs $(i,j)$ with $%
a_{i}=b_{j}$, and (\ref{e7}) becomes the equation%
\begin{equation*}
z(f)\cdot \lbrack \Ext^{2}(M,N)]=\left\vert q^{s(M)\cdot r(N)}\cdot
\lim_{t\rightarrow 1}\frac{P_{M^{\vee }\otimes N}(t)}{(1-t)^{\rho (M,N)}}%
\right\vert _{l}.
\end{equation*}

\item When the groups $\Ker(f)$, $\Coker(f)$, $\Ext_{\mathcal{R}%
{}_{l}}^{2}(M,N)$, $\Ext_{\mathcal{R}{}_{l}}^{3}(M,N)$, \ldots\ are finite
and almost all zero, then we say that $\chi ^{\times }(M,N)$ is defined, and
we set it equal to the alternating product of the orders of the groups.
Granted (b), Theorem \ref{lc.a} is equivalent to the statement: if $m_{M}(T)$
and $m_{N}(T)$ have no multiple root in common, then $\chi ^{\times }(M,N)$
is defined and equals%
\begin{equation}
\left\vert q^{s(M)\cdot r(N)}\cdot \prod_{a_{i}\neq b_{j}}\left( 1-\frac{%
b_{j}}{a_{i}}\right) \right\vert _{l}  \label{e4}
\end{equation}%
\qquad

\item Let
\begin{equation*}
0\rightarrow N^{\prime }\rightarrow N\rightarrow N^{\prime \prime
}\rightarrow 0
\end{equation*}%
be an exact sequence in $\mathcal{R}{}_{l}$. If $\chi ^{\times }(M,N^{\prime
})$ and $\chi ^{\times }(M,N^{\prime \prime })$ are defined, then so also is
$\chi ^{\times }(M,N)$, and
\begin{equation*}
\chi ^{\times }(M,N)=\chi ^{\times }(M,N^{\prime })\cdot \chi ^{\times
}(M,N^{\prime \prime })\text{.}
\end{equation*}%
Since a similar statement holds for exact sequences in $M$ and the
expression (\ref{e4}) is multiplicative, we see that it suffices to prove
Theorem \ref{lc.a} for $M$ and $N$ running through a set of generators for
the Grothendieck group $K(\mathcal{R}{}_{l})$ of $\mathcal{R}_{l}$.

Even more is true: both $\chi ^{\times }$ and the expression in (\ref{e4})
take values in $\mathbb{Q}{}_{+}^{\times }$, the group of positive rational
numbers, which is torsion-free. Therefore, in order to prove Theorem \ref%
{lc.a}, it suffices to verify it for $M$ and $N$ running through a set of
generators for $K(\mathcal{R}{}_{l})\otimes \mathbb{Q}{}$.
\end{enumerate}
\end{remark}

\subsection{Proof of Theorem \protect\ref{lc.a} in the case $l\neq p$}

From Proposition \ref{lc.e}, we obtain an exact sequence and isomorphisms%
\begin{eqnarray*}
\Hom_{\mathcal{R}{}_{\ell }}(M,N) &\cong &\Hom_{\mathbb{Z}{}_{\ell
}}(M,N)^{\Gamma } \\
0\rightarrow \Hom_{\mathbb{Z}{}_{\ell }}(M,N)_{\Gamma } &\xr{c}&\Ext_{%
\mathcal{R}{}_{\ell }}^{1}(M,N)\rightarrow \Ext_{\mathbb{Z}{}_{\ell
}}^{1}(M,N)^{\Gamma }\rightarrow 0 \\
\Ext_{\mathbb{Z}{}_{\ell }}^{1}(M,N)_{\Gamma } &\cong &\Ext_{\mathcal{R}%
{}_{l}}^{2}(M,N) \\
\Ext_{\mathcal{R}{}_{l}}^{i}(M,N) &=&0\text{ for }i>2\text{.}
\end{eqnarray*}%
Let $f_{0}$ be the map
\begin{equation*}
\Hom_{\mathbb{Z}{}_{\ell }}(M,N)^{\Gamma }\rightarrow \Hom_{\mathbb{Z}%
{}_{\ell }}(M,N)_{\Gamma }
\end{equation*}%
induced by the identity map. Then $f=c\circ f_{0}$, and so (see \ref{ez.d})%
\begin{equation}
z(f)=z(f_{0})/[\Ext_{\mathbb{Z}{}_{\ell }}^{1}(M,N)^{\Gamma }]\text{.}
\label{e1}
\end{equation}%
Let $\gamma $ be the generator $x\mapsto x^{q}$ of $\Gamma $. The family of
eigenvalues of $\gamma $ on $\Hom_{\mathbb{Z}{}_{\ell }}(M,N)$ is $(\frac{%
b_{i}}{a_{j}})_{i,j}$. Because of our hypothesis on $m_{M}(t)$ and $m_{N}(t)$%
, (\ref{ez.f}) applies, and shows that
\begin{equation*}
z(f_{0})=\left\vert \tprod\nolimits_{a_{i}\neq b_{j}}\left( 1-\frac{b_{i}}{%
a_{j}}\right) \right\vert _{\ell }.
\end{equation*}%
Because $\Ext_{\mathbb{Z}{}_{\ell }}^{1}(M,N)$ is finite,
\begin{equation*}
\lbrack \Ext_{\mathbb{Z}{}_{\ell }}^{1}(M,N)^{\Gamma }]=[\Ext_{\mathbb{Z}%
{}_{\ell }}^{1}(M,N)_{\Gamma }]
\end{equation*}%
and so%
\begin{equation*}
\lbrack \Ext_{\mathcal{R}{}_{l}}^{2}(M,N)]=[\Ext_{\mathbb{Z}{}_{\ell
}}^{1}(M,N)^{\Gamma }]\text{.}
\end{equation*}%
On substituting for $z(f_{0})$ and $[\Ext_{\mathbb{Z}{}_{\ell
}}^{1}(M,N)^{\Gamma }]$ in (\ref{e1}), we obtain (\ref{e7})$.$

\subsection{Proof of Theorem \protect\ref{lc.a} in the case $l=p$}

An $A$-module $M$ will be said to be $\emph{special}$ if it is cyclic, $%
M=A/A\cdot \lambda $, and $\lambda $ lies in the centre $\mathbb{Z}%
_{p}[F^{a}]$ of $A$, in which case $\lambda (F)=m_{M}(F^{a})$ and $%
P_{M}(t)=m_{M}(t)^{a}$. We let $k$ denote the $A$-module with underlying $W$%
-module $\mathbb{F}{}_{q}$ and with $F$ acting as $0$.

Let $\mathcal{R}{}_{p,f}$ be the full subcategory of $\mathcal{R}{}_{p}$ of
finite objects. Then Exts computed in $\mathcal{R}{}_{p,f}$ coincide with
those in $\mathcal{R}{}_{p}$ (see the argument p\pageref{p1}) and with those
in Pro-$\mathcal{R}{}_{p,f}$ ($\mathcal{R}{}_{p,f}$ is Artinian and so we
can apply the opposite of \ref{ha0}). There is a canonical functor $\mathcal{%
R}{}_{p}\rightarrow $Pro-$\mathcal{R}{}_{p,f}$ and by using the usual exact
sequences (\ref{e34}) one deduces that Exts in $\mathcal{R}{}_{p}$ coincide
with those in Pro-$\mathcal{R}{}_{p,f}$. This allows us to use resolutions
by free finitely generated $A$-modules to compute Exts in $\mathcal{R}{}_{p}$%
.

\subsubsection{The case $M$ or $N$ is finite with $F$ acting invertibly}

Because of (\ref{e35}), this case is essentially the same as the non-$p$
case. From the spectral sequence (\ref{lc.em}), we obtain exact sequences%
\begin{equation*}
0\rightarrow \Ext_{\mathcal{\bar{R}}_{p}{}}^{i-1}(M_{W(\mathbb{F}{})},N_{W(%
\mathbb{F}{})})_{\Gamma }\rightarrow \Ext_{\mathcal{R}{}_{p}}^{i}(M,N)%
\rightarrow \Ext_{\mathcal{\bar{R}}_{p}{}}^{i}(M_{W(\mathbb{F}{})},N_{W(%
\mathbb{F}{})})^{\Gamma }\rightarrow 0\text{.}
\end{equation*}%
In this case, the groups $\Ext_{\mathcal{\bar{R}}_{p}{}}^{i}(M_{W(\mathbb{F}%
{})},N_{W(\mathbb{F}{})})$ are finite, and it follows immediately that the
alternating product of the orders of the $\Ext_{\mathcal{R}%
{}_{p}{}}^{i}(M,N) $ is $1$.

\subsubsection{The case of $M=k$ and $N$ finite}

We have the following Koszul-type resolution of $k$

\begin{equation}
0\rightarrow A\xrightarrow{t}A\oplus A\xrightarrow{s}A\rightarrow
k\rightarrow 0  \label{e2}
\end{equation}%
where $t(\alpha )=(\alpha F,p\alpha )$ and $s(x,y)=px-yF$. For any $N$, the
groups $\Ext^{i}(k,N)$ are the cohomology groups of the complex%
\begin{equation}
\Hom(A,N)\xr{\circ s}\Hom(A\oplus A,N)\xr{\circ t}\Hom(A,N)\rightarrow
0\rightarrow \cdots  \label{e3}
\end{equation}%
Therefore, $\Ext^{i}(k,N)=0$ for $i>2$, and if $N$ is finite, then
\begin{equation*}
\Ext^{0}(k,N)\cdot \Ext^{2}(k,N)=\Ext^{1}(k,N)
\end{equation*}%
because the groups in (\ref{e3}) are finite.

\subsubsection{The case $M$ cyclic and $N=k$}

The resolution,%
\begin{equation*}
0\rightarrow A\xr{\cdot\lambda}A\rightarrow M\rightarrow 0,
\end{equation*}%
gives an exact sequence
\begin{equation*}
0\rightarrow \Hom(M,k)\rightarrow \Hom(A,k)\rightarrow \Hom(A,k)\rightarrow %
\Ext^{1}(M,k)\rightarrow 0
\end{equation*}%
and equalities $\Ext^{i}(M,k)=0$ ($i>1$). As $\Hom(A,k)\cong k$ is finite,
this implies that%
\begin{equation*}
\lbrack \Ext^{0}(M,k)]=[\Ext^{1}(M,k)]\text{.}
\end{equation*}

\subsubsection{The case $M=k$ and $N$ cyclic}

From (\ref{e2}) and (\ref{e3})$,$ we find that $\Ext^{i}(k,A)=0$ unless $i=2$
in which case $\Ext^{2}(k,A)=k$. Now the $\Ext(k,-)$-sequence of%
\begin{equation*}
0\rightarrow A\xr{\cdot\lambda}A\rightarrow M\rightarrow 0
\end{equation*}%
gives an exact sequence%
\begin{equation*}
0\rightarrow \Ext^{1}(k,N)\rightarrow \Ext^{2}(k,A)\rightarrow \Ext%
^{2}(k,A)\rightarrow \Ext^{2}(k,N)\rightarrow 0
\end{equation*}%
and equalities $\Ext^{i}(k,N)=0$ ($i\neq 1,2$). As $\Ext^{2}(k,A)$ is
finite, the sequence gives that%
\begin{equation*}
\lbrack \Ext^{1}(k,N)]=[\Ext^{2}(k,N)].
\end{equation*}

\subsubsection{The case $M$, $N$ special, $m_{M}$ and $m_{N}$ relatively
prime}

Let%
\begin{eqnarray*}
M &=&A/A\lambda _{1},\quad \lambda _{1}=m_{M}(F^{a}) \\
N &=&A/A\lambda _{2},\quad \lambda _{2}=m_{N}(F^{a})
\end{eqnarray*}%
with $m_{M}(t)$ and $m_{N}(t)$ relatively prime. From the resolution
\begin{equation*}
0\rightarrow A\xr{\cdot\lambda_1}A\rightarrow M\rightarrow 0
\end{equation*}%
we obtain an exact sequence%
\begin{equation*}
0\rightarrow \Hom(M,N)\rightarrow N\xr{\lambda_1\cdot}N\rightarrow \Ext%
^{1}(M,N)\rightarrow 0\rightarrow \ldots \text{.}
\end{equation*}%
As $\cdot \lambda _{1}$ is injective, we see that $\Ext^{i}(M,N)=0$ for $%
i\neq 1$, and so%
\begin{equation*}
z(f)\cdot \lbrack \Ext^{2}(M,N)]=z(\lambda _{1}\cdot )\overset{\ref{ez.c}}{=}%
\left\vert \det (\lambda _{1}\cdot )\right\vert _{p}^{a}.
\end{equation*}%
Here $\det (\lambda _{1}\cdot )$ can be computed as the determinant of $%
m_{M}(F^{a})$ acting on the $B$-vector space $N_{\mathbb{Q}{}}$. But $%
P_{M}(t)=m_{M}(t)^{a}$, and so%
\begin{equation*}
\det (m_{M}(F^{a})|N_{\mathbb{Q}{}})^{a}=\prod\nolimits_{%
\begin{smallmatrix}
1\leq i\leq r(M) \\
1\leq j\leq r(N)%
\end{smallmatrix}%
}(a_{i}-b_{j}).
\end{equation*}%
Thus,%
\begin{equation*}
z(f)\cdot \lbrack \Ext^{2}(M,N)]=\left\vert \prod (a_{i}-b_{j})\right\vert
_{p}=\left\vert q^{s(M)\cdot r(N)}\prod \left( 1-\frac{b_{j}}{a_{i}}\right)
\right\vert _{p}
\end{equation*}%
as required.

\subsubsection{The case $M$ and $N$ are special and equal}

We first give an explicit description of the map
\begin{equation*}
f\colon \Hom(M,N)\rightarrow \Ext^{1}(M,N)
\end{equation*}%
in the case that $M=A/A\cdot \lambda _{1}$ and $N=A/A\cdot \lambda _{2}$ are
special but not necessarily equal. We use $^{-}$ to denote $-\otimes _{W(%
\mathbb{F}{}_{q})}W(\mathbb{F}{})$, and we let $\sigma _{k}$ be the
canonical generator of $\Gal(\mathbb{F}{}/\mathbb{F}{}_{q}$). The $\Ext(-,N)$
sequence of
\begin{equation*}
0\rightarrow A\overset{\cdot \lambda _{1}}{\rightarrow }A\rightarrow
M\rightarrow 0
\end{equation*}%
is%
\begin{equation*}
0\rightarrow \Hom(M,N)\rightarrow N\overset{\lambda _{1}\cdot }{\rightarrow }%
N\rightarrow \Ext^{1}(M,N)\rightarrow 0
\end{equation*}%
and the $\Ext^{r}(-,\bar{N})$ sequence of%
\begin{equation*}
0\rightarrow \bar{A}\overset{\cdot \lambda _{1}}{\rightarrow }\bar{A}%
\rightarrow \bar{M}\rightarrow 0
\end{equation*}%
is%
\begin{equation*}
0\rightarrow \Hom(\bar{M},\bar{N})\rightarrow \bar{N}\overset{\lambda
_{1}\cdot }{\rightarrow }\bar{N}\rightarrow \Ext^{1}(\bar{M},\bar{N}%
)\rightarrow 0\text{.}
\end{equation*}%
The map $f$ can be described as follows: let $u\in \Hom(M,N)$ and regard $u$
as an element of $N$ such that $\lambda _{1}u=0$; then $u$ can be written $%
u=(\sigma _{k}-1)v$ for some $v\in \bar{N}$; now%
\begin{equation*}
(\sigma _{k}-1)(\lambda _{1}v)=\lambda _{1}(\sigma _{k}-1)v=\lambda _{1}u=0%
\text{,}
\end{equation*}%
and so $\lambda _{1}v\in \bar{N}^{\Gamma }=N$; the image of $\lambda _{1}v$
under $N\rightarrow \Ext^{1}(M,N)$ is $f(u)$.

Now take $N=M$, so that multiplication by $\lambda =\lambda _{1}=\lambda
_{2} $ is zero on $M$. Thus%
\begin{equation*}
\Hom(M,N)=A/A\lambda =\Ext^{1}(M,N)
\end{equation*}%
and $f$ is an endomorphism of $A/A/\lambda $. Since $A/A\lambda $ is
torsion-free, $z(f)$ is defined if and only if
\begin{equation*}
f_{\mathbb{Q}{}}\colon A_{\mathbb{Q}{}}/A_{\mathbb{Q}{}}\lambda \rightarrow
A_{\mathbb{Q}{}}/A_{\mathbb{Q}{}}\lambda
\end{equation*}%
has nonzero determinant (\ref{ez.c}), in which case%
\begin{equation*}
z(f)=|\det (f)|_{p}^{a}\text{.}
\end{equation*}

Let $u\in A_{\mathbb{Q}{}}/A_{\mathbb{Q}{}}\lambda $ and choose $v\in \bar{A}%
_{\mathbb{Q}{}}/\bar{A}_{\mathbb{Q}{}}\lambda $ such that $u=\sigma _{k}v-v$%
. Then $\sigma _{k}^{i}v=iu+v$ for all $i$. Let%
\begin{equation*}
\lambda (F)=F^{ma}+b_{ma-a}F^{ma-a}+\cdots +b_{0}=m(F^{a})\text{.}
\end{equation*}%
Then%
\begin{eqnarray*}
f(u) &=&\lambda (F)v \\
&=&muF^{ma}+(m-1)b_{ma-a}uF^{ma-a}+\cdots +0\quad \quad \text{(as }v\lambda
=0\text{)} \\
&=&uF^{a}\frac{d}{dF^{a}}(m(F^{a})) \\
&=&F^{a}\frac{d}{dF^{a}}(m(F^{a}))u\text{.}
\end{eqnarray*}%
From (\ref{ez.c}) we find that%
\begin{equation*}
z(f)=\left\vert \det \left( \frac{d}{dF^{a}}(m(F^{a}))\right) \right\vert
_{p}^{a}\cdot \left\vert \det (F^{a})\right\vert _{p}^{a}.
\end{equation*}%
But%
\begin{equation*}
|\det (F^{a})|_{p}^{a}=|q^{a\cdot s(N)}|_{p}
\end{equation*}%
and%
\begin{equation*}
\left\vert \det \left( \frac{d}{dF^{a}}(m(F^{a}))\right) \right\vert
_{p}^{a}=\left\vert \prod_{a_{i}\neq b_{j}}(a_{i}-b_{j})\right\vert
_{p}=\left\vert q^{s(N)(r(N)-a)}\prod_{a_{i}\neq b_{j}}\left( 1-\frac{b_{j}}{%
a_{i}}\right) \right\vert _{p}.
\end{equation*}%
Thus,%
\begin{equation*}
z(f)=\left\vert q^{r(N)\cdot s(N)}\prod_{a_{i}\neq b_{j}}\left( 1-\frac{b_{j}%
}{a_{i}}\right) \right\vert _{p}.
\end{equation*}%
As $\Ext^{2}(M,N)=0$ this proves Theorem \ref{lc.a} in this case.

\subsubsection{ The general case}

An isogeny (denoted $\sim $) of objects in $\mathcal{R}{}_{p}$ is a
homomorphism whose kernel and cokernel have finite length over $W$. We say
an object in $\mathcal{R}{}_{p}$ is indecomposable if it is not isogenous to
a direct sum of two nonzero objects.

\begin{lemma}
\label{lc.f}Every object of $\mathcal{R}{}_{p}$ is isogenous to a direct sum
of indecomposable objects, and the decomposition is unique up to isogeny.
\end{lemma}

\begin{proof}
Apply the Krull-Schmidt theorem in the category of $A_{\mathbb{Q}{}}$%
-modules.
\end{proof}

\begin{lemma}
\label{lc.g}If $M$ in $\mathcal{R}{}_{p}$ is indecomposable, then $M_{%
\mathbb{Q}{}}\approx A_{\mathbb{Q}{}}/A_{\mathbb{Q}{}}\lambda $ for some $%
\lambda \in A_{\mathbb{Q}{}}$; moreover, the minimum polynomial $m_{M}(t)$
of the Frobenius endomorphism of $M$ is a power of a $\mathbb{Z}{}_{p}$%
-irreducible polynomial, and there exists an integer $e$ such that $\oplus
^{e}M\sim A/A\cdot m_{M}(t)$.
\end{lemma}

\begin{proof}
The ring $A_{\mathbb{Q}{}}$ is the skew polynomial ring $B[F,\sigma ]$. It
is therefore a principal ideal domain (Jacobson 1943, III 1, p30), and every
finitely generated module over such a ring is isomorphic a direct sum of
cyclic modules (ibid. III, Theorem 19, p44). This proves the first part.

For the second, we use that the two-sided ideals in $A_{\mathbb{Q}{}}$ are
precisely those generated by a polynomial in the centre $\mathbb{Q}%
{}_{p}[F^{a}]$ of $A_{\mathbb{Q}{}}$ (ibid. III 5, p38). It follows that $A_{%
\mathbb{Q}{}}\cdot m(F^{a})$ is the largest two-sided ideal contained in $A_{%
\mathbb{Q}{}}\cdot \lambda $, i.e., it is the bound of the ideal $A_{\mathbb{%
Q}{}}\cdot \lambda $ in the sense of ibid. III 6, p38. Because $A_{\mathbb{Q}%
{}}/A_{\mathbb{Q}{}}\lambda $ is indecomposable, this implies that $m(t)$ is
a power of a $\mathbb{Z}{}_{p}$-irreducible polynomial (ibid. III Theorem
13, p40). Finally, when $A_{\mathbb{Q}{}}/A_{\mathbb{Q}{}}\cdot m(F^{a})$ is
decomposed into a direct sum of indecomposable modules $A_{\mathbb{Q}{}}/A_{%
\mathbb{Q}{}}\lambda ^{\prime }$, each of the ideals $A_{\mathbb{Q}%
{}}\lambda ^{\prime }$ has $A_{\mathbb{Q}{}}\cdot m(F^{a})$ as its bound,
and is therefore isomorphic to $M$ (ibid. III Theorem 20, p45).
\end{proof}

\begin{proposition}
\label{lc.h}Objects of the following types generate $K(\mathcal{R}{}_{p})_{%
\mathbb{Q}{}}$:

\begin{enumerate}
\item $k=(\mathbb{F}{}_{q},0)$;

\item $(M,F_{M})$ with $M$ finite and $F_{M}$ invertible;

\item $M$ special and $m_{M}(t)$ is a power of an irreducible polynomial.
\end{enumerate}
\end{proposition}

\begin{proof}
For every finite module $M$, there exists an $n$ such that $M=\im
(F^{n})\oplus \Ker(F^{n})$ (Fitting's Lemma). It follows that such an $M$
has a composition series with quotients of type (a) or (b). The general
statement now follows from Lemma \ref{lc.g}.
\end{proof}

Therefore (see \ref{lc.c}(f)), it suffices to prove Theorem \ref{lc.a} when
each of $M$ and $N$ is of type (a), (b), or (c). Moreover, because of the
condition on the minimum polynomials, in case (c) we can assume that either $%
m_{M}$ and $m_{N}$ have no common factor or they are irreducible and equal.
This we have done.

\subsection{Restatement of the theorem for complexes}

We sketch a generalization of Theorem \ref{lc.a} to complexes. Let $%
M^{\bullet }$ and $N^{\bullet }$ be bounded complexes of $A$-modules such
that $H^{i}(M^{\bullet })$ lies in $\mathcal{R}{}_{p}$ for all $i$. Define%
\begin{equation*}
P_{M^{\bullet \vee }\otimes N^{\bullet }}(t)=\prod P_{H^{i}(M^{\bullet
})^{\vee }\otimes H^{j}(N^{\bullet })}(t)^{(-1)^{i+j}}.
\end{equation*}%
Let $f^{i}\colon \Ext^{i}(M^{\bullet },N^{\bullet })\rightarrow \Ext%
^{i+1}(M^{\bullet },N^{\bullet })$ be the map rendering

\begin{equation*}
\begin{CD} \Ext^{i}(M^{\bullet },N^{\bullet }) @>f^{i}>>
\Ext^{i+1}(M^{\bullet },N^{\bullet }) \\ @VVV @VVV\\
\Ext^{i}(\bar{M}^{\bullet },\bar{N}^{\bullet })^{\Gamma } @>\cup \theta
_{p}>> \Ext^{i}(\bar{M}^{\bullet },\bar{N}^{\bullet })_{\Gamma }\end{CD}
\end{equation*}%
commutative. Here the vertical maps arise from a spectral sequence (cf. \ref%
{lc.em}) and the lower map is cup-product with the canonical generator $%
\theta _{p}\in H^{1}(\Gamma ,W)$ (equal to the map induced by the identity
map on $\Ext^{i}(\bar{M}^{\bullet },\bar{N}^{\bullet })$). Because $\theta
_{p}^{2}=0$,
\begin{equation*}
\Ext^{\bullet }(M^{\bullet },N^{\bullet })\colon \quad \cdots \rightarrow %
\Ext^{i}(M^{\bullet },N^{\bullet })\overset{f^{i}}{\rightarrow }\Ext%
^{i+1}(M^{\bullet },N^{\bullet })\rightarrow \cdots
\end{equation*}%
is a complex. Define%
\begin{equation*}
\chi ^{\times }(M^{\bullet },N^{\bullet })=\prod [H^{i}(\Ext^{\bullet
}(M^{\bullet },N^{\bullet }))]^{(-1)^{i}}
\end{equation*}%
when these numbers are finite.

\begin{theorem}
\label{lc.i}Let $M^{\bullet }$ and $N^{\bullet }$ be bounded complexes of $A$%
-modules such that $H^{i}(M^{\bullet })$ and $H^{i}(N^{\bullet })$ are
semisimple and finitely generated over $W$ for all $i$.

\begin{enumerate}
\item The groups $\Ext^{i}(M^{\bullet },N^{\bullet })$ are finitely
generated $\mathbb{Z}{}_{p}$-modules for all $i$.

\item The alternating sum $\tsum (-1)r_{i}$ of the ranks of the $\Ext%
^{i}(M^{\bullet },N^{\bullet })$ is zero.

\item The order of zero of $P(M^{\bullet \vee }\otimes N^{\bullet },t)$ at $%
t=1$ is equal to the secondary Euler characteristic $\tsum (-1)^{i+1}\cdot
i\cdot r_{i}.$

\item The groups $H^{i}(\Ext^{\bullet }(M^{\bullet },N^{\bullet }))$ are
finite, and%
\begin{equation*}
\left\vert q^{\chi (M^{\bullet },N^{\bullet })}\cdot \lim_{s\rightarrow 1}%
\frac{P_{M^{\bullet \vee }\otimes N^{\bullet }}(t)}{(1-t)^{\rho (M^{\bullet
},N^{\bullet })}}\right\vert _{l}=\chi ^{\times }(M^{\bullet },N^{\bullet })%
\text{ }
\end{equation*}%
where $\chi (M^{\bullet },N^{\bullet })=\tsum (-1)^{i+j}s(H^{i}(M^{\bullet
}))\cdot r(H^{j}(N^{\bullet }))$.
\end{enumerate}
\end{theorem}

\begin{proof}
Consider first the case that $M^{\bullet }$ and $N^{\bullet }$ are modules $%
M $ and $N$ regarded as complexes concentrated in degree $0$.

(a) Because $M$ and $N$ are semisimple, $m(M,t)$ and $m(N,t)$ have only
simply roots, and therefore the hypothesis of Theorem \ref{lc.a} hold.
Hence, the groups $\Ext^{i}(M,N)$ are finitely generated $\mathbb{Z}{}_{p}$%
-modules for all $i$, finite for $i=2$, and zero for $i>2$.

(b) Moreover, $f\colon \Hom(M,N)\rightarrow \Ext^{1}(M,N)$ has finite kernel
and cokernel, and so $\Hom(M,N)$ and $\Ext^{1}(M,N)$ have the same rank.

(c) In fact, one sees easily that the rank of $\Hom(M,N)$ is the number of
pairs $(i,j)$ such that $a_{i}=b_{j}$ (see \ref{lc.c}d).

(d) The map $f$ of Theorem \ref{lc.a} coincides with $f^{0}$, and so (d) is
a restatement of the theorem.

The general case case now follows from the usual induction argument on the
number of nonzero terms in the complexes $M^{\bullet }$ and $N^{\bullet }$
(cf. \ref{lc.c}f).
\end{proof}

There is a similar Theorem for $l\neq p$, whose formulation and proof we
leave to the reader.

\subsection{Why $\mathcal{R}^{+}(\mathbb{F}_{q};\mathbb{Z}_{p})$ rather than
$\mathcal{R}(\mathbb{F}_{q};\mathbb{Z}_{p})$?}

Write $\mathcal{R}{}^{+}$ and $\mathcal{R}{}$ for $\mathcal{R}^{+}(\mathbb{F}%
_{q};\mathbb{Z}_{p})$ and $\mathcal{R}(\mathbb{F}_{q};\mathbb{Z}_{p})$
respectively. Even the torsion subgroup of $\Ext_{\mathcal{R}}^{1}(M,N)$ can
be infinite. Suppose, for example, that $M$ and $N$ are torsion-free, and
consider the diagram%
\begin{equation*}
\begin{CD} 0 @>>> \Hom_{\mathcal{R}^{+}}(M,N)^{(p^{n})} @>>>
\Hom_{\mathcal{R}^{+}}(M,N^{(p^{n})}) @>>>
\Ext_{\mathcal{R}^{+}}^{1}(M,N)_{p^{n}} @>>> 0 \\ @. @VV{\cong}V @VVV@VVV \\
0 @>>> \Hom_{\mathcal{R}}(M,N)^{(p^{n})} @>>>
\Hom_{\mathcal{R}}(M,N^{(p^{n})}) @>>> \Ext_{\mathcal{R}}^{1}(M,N)_{p^{n}}
@>>> 0. \end{CD}
\end{equation*}%
Since $[\Hom_{\mathcal{R}{}}(M,N^{(p^{n})})]=(q^{n})^{r(M)r(N)}$ is
unbounded, whereas $[\Hom_{\mathcal{R}{}^{+}}(M,N^{(p^{n})})]$ is bounded if
$\Hom_{R{}^{+}}(M,N)=0$, we see that $\Ext_{\mathcal{R}}^{1}(M,N)(p)$ is
frequently infinite. From a different perspective, it follows from Theorem
9.1 that $\Ext_{\mathcal{R}{}^{+}}^{1}(M,N)$ grows when $M$ and $N$ are
replaced by $M(1)$ and $N(1)$ because $s(M(1))>s(M)$, but the functor $%
M\mapsto M(1)$ is an equivalence on the category $\mathcal{R}{}(\mathbb{F}%
{}_{q};\mathbb{Z}{}_{p})$.

\subsection{Appendix: Lemmas on abelian groups}

Recall that, for a homomorphism $f\colon M\rightarrow N$ of abelian groups,
we define $z(f)=\frac{[\Ker(f)]}{[\Coker(f)]}$ when both the top and bottom
are finite. There are the following elementary statements (Tate 1966b, \S %
5). We let $W=W(\mathbb{\mathbb{F}{}}_{q})$, $q=p^{a}$.

\begin{E}
\label{ez.c}Let $M$ and $N$ be finitely generated $\mathbb{Z}{}$-modules
(resp. $W$-modules) of the same rank and let $(x_{i})$ and $(y_{i})$ be
bases for $M$ and $N$ respectively modulo torsion. Suppose $f(x_{i})\equiv
\sum a_{ij}y_{j}$ modulo torsion. Then $z(f)$ is defined if and only if $%
\det (a_{ij})\neq 0$, in which case%
\begin{equation*}
z(f)=\frac{[M_{\text{tors}}]}{|\det (a_{ij})|[N_{\text{tors}}]}
\end{equation*}%
(resp.%
\begin{equation*}
z(f)=\frac{[M_{\text{tors}}]|\det (a_{ij})|_{p}^{a}}{[N_{\text{tors}}]}\text{%
).}
\end{equation*}
\end{E}

Note that, when $M=N$, $\det (a_{ij})$ is the determinant of $f_{\mathbb{Q}%
{}}\colon M_{\mathbb{Q}{}}\rightarrow M_{\mathbb{Q}{}}$.

\begin{E}
\label{ez.d}Consider maps $f\colon M\rightarrow N$ and $g\colon N\rightarrow
P$. If any two of the three numbers $z(f)$, $z(g)$, $z(g\circ f)$ are
defined, then so is the third, and%
\begin{equation*}
z(g\circ f)=z(g)\cdot z(f).
\end{equation*}
\end{E}

\begin{E}
\label{ez.f}Let $\Gamma =\mathbb{\hat{Z}}{}$, and let $\gamma $ be its
canonical topological generator. Let $M$ be a $\Gamma $ -module that is
finitely generated as a $\mathbb{Z}{}$-module, and let $f\colon M^{\Gamma
}\rightarrow M_{\Gamma }$ be the map induced by the identity map on $M$.
Then $z(f)$ is defined if and only if the minimum polynomial of $\gamma $ on
$M_{\mathbb{Q}{}}$ does not have $1$ as a multiple root, in which case%
\begin{equation*}
z(f)\left| \prod_{a_{i}\neq 1}(1-a_{i})\right| =1
\end{equation*}%
where $a_{1},a_{2},\ldots $ is the family of eigenvalues of $\gamma $ acting
on $M_{\mathbb{Q}{}}$.
\end{E}

\section{Extensions and Zeta Functions: Global Case}

In this section, $\mathcal{M}{}(\mathbb{F};\mathbb{Q}{})$ will be one of the
following categories:

\begin{enumerate}
\item the category $\mathcal{M}_{\text{num}}(\mathbb{F};\mathbb{Q};\mathcal{S%
})$ where $\mathcal{S}{}$ consists of all smooth projective varieties over $%
\mathbb{F}{}$ and we assume that the Tate conjecture holds for all $V\in
\mathcal{S}{};$

\item the category $\mathcal{M}{}(\mathbb{F}{};\mathbb{Q}{})$ defined in (%
\ref{cm3}).
\end{enumerate}

\noindent In fact, if the category in (a) exists and the Hodge conjecture
holds for complex abelian varieties of CM-type, then these categories are
essentially the same: the reduction functor $\CM(\mathbb{Q}^{\text{al}%
})\rightarrow \mathcal{M}_{\text{num}}(\mathbb{F};\mathbb{Q};\mathcal{S})$
defines a fibre functor $\omega $ on $\CM(\mathbb{Q}^{\text{al}})^{P}$, and
the quotient of $\CM(\mathbb{Q}^{\text{al}})$ corresponding to $\omega $ is
canonically equivalent with $\mathcal{M}_{\text{num}}(\mathbb{F};\mathbb{Q};%
\mathcal{S})$.

\subsection{The main theorem for extensions}

\begin{theorem}
\label{gc.a}For all effective motives $X$ and $Y$ in $\mathcal{M}{}^{+}(%
\mathbb{F}_{q};\mathbb{Z}{})$, the groups $\Ext^{1}(X,Y)$ and $\Ext%
^{2}(X,Y)_{\text{cotors}}$ are finite, and%
\begin{equation*}
q^{\chi (X,Y)}\cdot \zeta (X^{\vee }\otimes Y,s)\sim \pm \frac{\lbrack \Ext%
^{1}(X,Y)]\cdot D(X,Y)}{[\Hom(X,Y)_{\text{tors}}]\cdot \lbrack \Ext%
^{2}(X,Y)_{\text{cotors}}]}\cdot (1-q^{-s})^{\rho (X,Y)}\text{ }\quad \text{%
as }s\rightarrow 0
\end{equation*}%
where $\chi (X,Y)=r(X_{p})s(Y_{p})$ (notations as on p\pageref{localext}), $%
D=D(X,Y)$ is the discriminant of the pairing%
\begin{equation}
\Hom(Y,X)\times \Hom(X,Y)\xr{\circ}\End(Y)\xr{\text{trace}}\mathbb{Z}{},
\label{e16}
\end{equation}%
and $\rho (X,Y)$ is the rank of $\Hom(X,Y)$.
\end{theorem}

\begin{remark}
\label{gc.c}Let $X$ and $Y$ be effective motives over $\mathbb{F}{}_{q}$.

\begin{enumerate}
\item Because $\mathcal{M}{}^{+}(\mathbb{F}{}_{q};\mathbb{Q}{})$ is a
semisimple category, $\End(X_{0}\times Y_{0})$ is a semisimple ring. As it
is also finite dimensional over $\mathbb{Q}{}$, the trace pairing%
\begin{equation*}
\End(X_{0}\times Y_{0})\times \End(X_{0}\times Y_{0})\rightarrow \mathbb{Q}{}
\end{equation*}%
is nondegenerate: $D(X,Y)\neq 0$.

\item Below we prove that $[\Ext^{2}(X,Y)_{\text{cotors}}]=[\Hom(Y,X)_{\text{%
tors}}]$ . This makes the formula more symmetric --- note that $%
D(X,Y)=D(Y,X) $.

\item When we replace $X$ by $X(r)$ or $Y$ by $Y(r)$ ($r\in \mathbb{N}{}$)
in Theorem \ref{gc.a}, we obtain a description of the behaviour of $\zeta
(X^{\vee }\otimes Y,s)$ near $-r$ or $+r$.
\end{enumerate}
\end{remark}

\subsubsection{Proof of Theorem \protect\ref{gc.a} with $z(\protect%
\varepsilon )$ for $D(X,Y)$}

Recall (\ref{ha6}) that for all $l$, there is an exact sequence%
\begin{equation*}
0\rightarrow \Ext_{\mathcal{M}{}^{+}(\mathbb{F}{}_{q};\mathbb{Z}%
{})}^{1}(X,Y)_{\mathbb{Z}{}_{l}}\xr{c_l}\Ext_{\mathcal{R}{}^{+}(\mathbb{F}%
{}_{q};\mathbb{Z}_{l}{})}^{1}(X_{l},Y_{l})\rightarrow T_{l}\Ext_{\mathcal{M}%
{}^{+}(\mathbb{F}{}_{q};\mathbb{Z}{}_{l})}^{2}(X,Y)\rightarrow 0\text{.}
\end{equation*}%
Define
\begin{equation*}
\varepsilon =(\varepsilon _{l})_{l}\colon \Hom_{\mathcal{M}{}^{+}(\mathbb{F}%
{}_{q};\mathbb{Z}{})}(X,Y)_{\mathbb{\hat{Z}}{}}\rightarrow T\Ext_{\mathcal{M}%
{}^{+}(\mathbb{F}{}_{q};\mathbb{Z}{})}^{2}(X,Y)
\end{equation*}%
to be the map making%
\begin{equation*}
\begin{CD}
\Hom_{\mathcal{M}^{+}(\mathbb{F}_{q};\mathbb{Z})}(X,Y)_{\mathbb{Z}_{l}}
@>\varepsilon_{l}>>
T_{l}\Ext_{\mathcal{M}^{+}(\mathbb{F}_{q};\mathbb{Z})}^{2}(X,Y) \\
@VV{\cong}V @AA{c_{l}}A \\
\Hom_{\mathcal{R}^{+}(\mathbb{F}_{q};\mathbb{Z}_{l})}(X_{l},Y_{l}) @>f_{l}>>
\Ext_{\mathcal{R}^{+}(\mathbb{F}_{q};\mathbb{Z}_{l})}^{1}(X_l,Y_l)\end{CD}
\end{equation*}%
commute. Here $f_{l}$ $=f(X_{l},Y_{l})$ is the map defined on p\pageref{p2}.
According to \ref{ez.d} and \ref{lc.a}%
\begin{equation}
z(\varepsilon _{l})=z(c_{l})\cdot z(f_{l})=\frac{[\Ext^{1}(X,Y)(l)]}{[\Ext%
^{2}(X_{l},Y_{l})]}\cdot \left\vert q^{s(X_{l})\cdot r(Y_{l})}\cdot
\prod_{a_{i}\neq b_{j}}\left( 1-\frac{b_{j}}{a_{i}}\right) \right\vert _{l}%
\text{.}  \label{e28}
\end{equation}%
The $a_{i},b_{j}$ are algebraic numbers independent of $l$. Because $\Ext%
^{i}(X_{l},Y_{l})=0$ for $i\geq 3$, (\ref{ha6}b) shows that $\Ext^{i}(X,Y)=0$
for $i\geq 3$, and it follows (from \ref{ha6}b again) that%
\begin{equation*}
\Ext^{2}(X,Y)_{\text{cotors}}(l)\cong \Ext^{2}(X_{l},Y_{l})\text{.}
\end{equation*}%
On multiplying the equations (\ref{e28}) for the different $l$ and applying
the product formula, we find that%
\begin{equation*}
q^{\chi (X,Y)}\cdot \prod_{a_{i}\neq b_{j}}\left( 1-\frac{b_{j}}{a_{i}}%
\right) \cdot z(\varepsilon )=\pm \frac{\lbrack \Ext^{1}(X,Y)]}{[\Ext%
^{2}(X,Y)_{\text{cotors}}]}.
\end{equation*}%
This will become Theorem \ref{gc.a} once we have shown that%
\begin{equation*}
z(\varepsilon )=\pm \frac{\lbrack \Hom(X,Y)_{\text{tors}}]}{D(X,Y)}.
\end{equation*}

\subsubsection{Comparison of $z(\protect\varepsilon )$ with $D(X,Y)$}

Let $N$ be a finitely generated torsion-free discrete $\Gamma $-module, and
let $N^{\vee }=\Hom(N,\mathbb{Z})$. The pairing%
\begin{equation*}
N^{\Gamma }\times H^{2}(\Gamma ,N^{\vee })\rightarrow H^{2}(\Gamma ,\mathbb{Z%
}{})\cong \mathbb{Q}{}/\mathbb{Z}{}
\end{equation*}%
realizes the compact group $(N^{\Gamma })_{\mathbb{\hat{Z}}}$ as the
Pontryagin dual of the discrete group $H^{2}(\Gamma ,N^{\vee })$ (e.g.,
Milne 1986b, I 1.10). There is therefore a canonical isomorphism%
\begin{equation*}
\Hom(N^{\Gamma },\mathbb{\hat{Z}})\rightarrow \Hom(\mathbb{Q}{}{}/\mathbb{Z}%
{},H^{2}(\Gamma ,N^{\vee }))\overset{\text{df}}{=}TH^{2}(\Gamma ,N^{\vee })%
\text{.}
\end{equation*}%
Suppose we are also given a nondegenerate pairing
\begin{equation*}
\psi \colon N\times M\rightarrow \mathbb{Z}{}
\end{equation*}%
where $M$ is a discrete $\Gamma $-module such that $M_{\text{tors}}$ is
killed by some integer and $M/M_{\text{tors}}$ is finitely generated. This
defines a map $M\rightarrow N^{\vee }$ with torsion kernel and cokernel,
which induces isomorphisms $H^{2}(\Gamma ,M)\rightarrow H^{2}(\Gamma
,N^{\vee })$ and $TH^{2}(\Gamma ,M)\rightarrow TH^{2}(\Gamma ,N^{\vee })$.
Hence, there is a canonical isomorphism%
\begin{equation*}
\Hom(N^{\Gamma },\mathbb{\hat{Z}}{})\rightarrow TH^{2}(\Gamma ,M).
\end{equation*}%
When combined with the map%
\begin{equation*}
M^{\Gamma }\rightarrow \Hom(N^{\Gamma },\mathbb{Z}{}\mathbb{)}{}
\end{equation*}%
defined by $\psi $, this gives a homomorphism%
\begin{equation*}
\varepsilon (\psi )\colon (M^{\Gamma })_{\mathbb{\hat{Z}}{}}\rightarrow
TH^{2}(\Gamma ,M)\text{.}
\end{equation*}%
Note that $z(\varepsilon (\psi ))$ is defined if and only if $M_{\text{tors}%
}^{\Gamma }$ is finite and the restriction $\psi ^{\Gamma }$ to a pairing $%
N^{\Gamma }\times M^{\Gamma }\rightarrow \mathbb{Z}{}$ is nondegenerate, in
which case%
\begin{equation*}
z(\varepsilon (\psi ))=\pm \frac{\lbrack M_{\text{tors}}^{\Gamma }]}{\det
(\psi ^{\Gamma })}\text{.}
\end{equation*}%
When we apply this remark to the pairing
\begin{equation*}
\Hom(\bar{Y},\bar{X})/\{\text{torsion}\}\times \Hom(\bar{X},\bar{Y}%
)\rightarrow \mathbb{Q}{},
\end{equation*}
we obtain a homomorphism%
\begin{equation*}
\varepsilon \colon \Hom(X,Y)_{\mathbb{\hat{Z}}{}}\rightarrow TH^{2}(\Gamma ,%
\Hom(X,Y))\text{,}
\end{equation*}%
and%
\begin{equation*}
z(\varepsilon )=\pm \frac{\lbrack \Hom(X,Y)_{\text{tors}}]}{D(X,Y)}\text{.}
\end{equation*}%
Since $\varepsilon $ coincides with the map defined in the preceding
subsubsection, this completes the proof of Theorem \ref{gc.a}.

\subsection{Comparison of $\Hom(Y,X)_{\text{tors}}$ with $\Ext^{2}(X,Y)_{%
\text{cotors}}$}

In this subsection, we prove that%
\begin{equation*}
\lbrack \Hom(Y,X)_{\text{tors}}]=[\Ext^{2}(X,Y)_{\text{cotors}}].
\end{equation*}%
It follows from Theorem \ref{ha6} that%
\begin{equation*}
\Hom(Y,X)(l)\cong \Hom(Y_{l},X_{l})_{\text{tors}},\quad \Ext^{2}(X,Y)_{\text{%
cotors}}(l)\cong \Ext^{2}(X_{l},Y_{l})_{\text{cotors}},
\end{equation*}%
and so it suffices to show that, for $M$ and $N$ in $\mathcal{R}{}^{+}(%
\mathbb{F}_{q};\mathbb{Z}_{l})$,%
\begin{equation*}
\lbrack \Hom(N,M)_{\text{tors}}]=[\Ext^{2}(M,N)_{\text{cotors}}]\text{.}
\end{equation*}%
Because $\Ext^{3}=0$ in $\mathcal{R}{}^{+}(\mathbb{F}_{q};\mathbb{Z}_{l})$, $%
\Ext^{2}(M/M_{\text{tors}},N)$ is divisible and%
\begin{equation*}
\Ext^{2}(M,N)_{\text{cotors}}\cong \Ext^{2}(M_{\text{tors}},N)\text{.}
\end{equation*}%
As $\Hom(N,M_{\text{tors}})\cong \Hom(N,M)_{\text{tors}}$, this shows that
it suffices to prove that%
\begin{equation}
\lbrack \Hom(N,M)]=[\Ext^{2}(M,N)]  \label{e36}
\end{equation}%
when $M$ is torsion. In fact, we shall show that%
\begin{equation}
\Hom(N,M)\cong \Ext^{2}(M,N)^{\ast }  \label{e37}
\end{equation}%
when $N$ is torsion-free and $M$ is torsion. Here $\ast $ denotes $\Hom(-,%
\mathbb{Q}{}/\mathbb{Z}{})$. This implies (\ref{e36}) because an arbitrary $N
$ has a resolution by torsion-free objects (1.4b, 7.10).

Define%
\begin{equation*}
\Ext^{i}(M,N\otimes \mathbb{Q}{}/\mathbb{Z}{})=\varinjlim_{n}\Ext%
^{i}(M,N^{(l^{n})}).
\end{equation*}%
Then%
\begin{equation*}
\Ext^{1}(M,N\otimes \mathbb{Q}{}/\mathbb{Z}{})\cong \Ext^{2}(M,N)
\end{equation*}%
and so we have to show that%
\begin{equation}
\Hom(N,M)\cong \Ext^{1}(M,N\otimes \mathbb{Q}{}/\mathbb{Z}{})^{\ast }
\label{e38}
\end{equation}%
($M$ torsion, $N$ torsion-free).

\subsubsection{Proof of (\protect\ref{e38}) in the case $l\neq p$}

We let $\bar{M}$ and $\bar{N}$ denote $M$ and $N$ regarded as $\mathbb{Z}%
{}_{l}$-modules. There is a perfect pairing of finite $\mathbb{Z}{}_{l}$%
-modules%
\begin{equation*}
\Hom(\bar{N},\bar{M})\times \Hom(\bar{M},\bar{N}\otimes \mathbb{Q}{}/\mathbb{%
Z}{})\overset{\circ }{\rightarrow }\Hom(\bar{N},\bar{N}\otimes \mathbb{Q}{}/%
\mathbb{Z}{})\cong \End(\bar{N})\otimes \mathbb{Q}{}/\mathbb{Z}{}%
\xr{\Tr\otimes 1} \mathbb{Q}{}/\mathbb{Z}{}\text{.}
\end{equation*}%
This gives a perfect pairing of $\Gamma $-cohomology groups%
\begin{equation*}
\Hom(\bar{N},\bar{M})^{\Gamma }\times \Hom(\bar{M},\bar{N}\otimes \mathbb{Q}%
{}/\mathbb{Z}{})_{\Gamma }\rightarrow \mathbb{Q}{}/\mathbb{Z}{},
\end{equation*}%
which can be identified with (\ref{e38}).

\subsubsection{Proof of (\protect\ref{e38}) in the case $l=p$}

If $F$ acts as an isomorphism on $M$, this can be proved as in the case $%
l\neq p$. Thus, we may assume that $F$ act nilpotently on $M$, which
therefore has a composition series whose quotients are isomorphic to $k=_{%
\text{df}}(\mathbb{F}{}_{q},0)$.

Any extension of $N^{(p^{n})}$ by $N^{(p^{n})}$ splits as an extension of $W$%
-modules, and therefore is described by a $\sigma $-linear endomorphism $%
\alpha $ of $N^{(p^{n})}$. The trace of $\alpha $ acting on $N^{(p^{n})}$ is
an element of $W/p^{n}W$, and the trace of this element is a well-defined
element of $\mathbb{Z}{}_{p}/p^{n}\mathbb{Z}{}_{p}$. Thus, we have a
homomorphism%
\begin{equation*}
\Ext^{1}(N^{(p^{n})},N^{(p^{n})})\rightarrow \mathbb{Z}{}_{p}/p^{n}\mathbb{Z}%
{}_{p}\cong p^{-n}\mathbb{Z}{}_{p}/\mathbb{Z}{}_{p}\text{.}
\end{equation*}%
On passing to the limit, we obtain a homomorphism%
\begin{equation*}
t\colon \Ext^{1}(N,N\otimes \mathbb{Q}{}/\mathbb{Z}{})\rightarrow \mathbb{Q}%
{}/\mathbb{Z}{}.
\end{equation*}%
We claim that the pairings%
\begin{equation}
\Ext^{i}(N,M)\times \Ext^{1-i}(M,N\otimes \mathbb{Q}{}/\mathbb{Z}%
{})\rightarrow \Ext^{1}(N,N\otimes \mathbb{Q}{}/\mathbb{Z}{})\overset{t}{%
\rightarrow }\mathbb{Q}{}/\mathbb{Z}{}  \label{e39}
\end{equation}%
are perfect for $i\geq 0$. When $N$ is cyclic, $N=A/A\cdot \lambda $, and $%
M=k$, this can be proved directly: either both groups are zero, or the
pairing is the trace pairing%
\begin{equation*}
\mathbb{F}{}_{q}\times \mathbb{F}{}_{q}\rightarrow \mathbb{F}{}_{p}\text{.}
\end{equation*}%
Consider an exact sequence%
\begin{equation*}
0\rightarrow N^{\prime }\rightarrow N\rightarrow k\rightarrow 0\text{.}
\end{equation*}%
If the pairings (\ref{e39}) are perfect for $N$, $k$, and $i\geq 0$, then
the diagram%
%TCIMACRO{\TeXButton{footnotesize}{\footnotesize} }%
%BeginExpansion
\footnotesize
%EndExpansion
\begin{equation*}
\begin{CD} 0 @>>> \Hom(k,k) @>>> \Hom(N,k) @>>>
\Hom(N^{\prime },k) @>>> \cdots\\ @.@VV{\cong}V @VV{\cong}V
 @VVV @VV{\cong}V \\
 0 @>>> \Ext^{2}(k,k)^{\ast } @>>>
 \Ext^{1}(k,N\otimes \mathbb{Q}{}/\mathbb{Z}{})^{\ast } @>>>
 \Ext^{1}(k,N^{\prime }\otimes \mathbb{Q}{}/\mathbb{Z}{})^{\ast } @>>>
\cdots\end{CD}
\end{equation*}%
%TCIMACRO{\TeXButton{normalsize}{\normalsize}}%
%BeginExpansion
\normalsize%
%EndExpansion
shows that they are perfect for $N^{\prime }$, $k$, and $i\geq 0$ (note that
(24) shows that $\Ext^{i}(k,k)$ equals $k,k\oplus k,k$ respectively for $%
i=0,1,2$). A similar argument applies when $N/N^{\prime }$ is a finite
module on which $F$ acts invertibly. Since an arbitrary torsion-free $N$ can
be realized as a submodule of a direct sum of cyclic modules (9.3), this
proves that the pairings (\ref{e39}) are perfect when $M=k$. An induction
argument on the length of $M$ extends this to an arbitrary torsion $M$. The
case $i=0$ completes the proof of (\ref{e38}).

\subsection{The main theorem for Weil extensions}

Let $\Gamma _{0}$ be the subgroup of $\Gamma $ generated by the Frobenius
element (so $\Gamma _{0}\cong \mathbb{Z}$). Inspired by Lichtenbaum 2002, we
define \textquotedblleft Weil\textquotedblright\ extension groups satisfying
an analogue of the spectral sequence (\ref{e29}).

As we noted in \ref{ct12}, to give a motive over $\mathbb{F}{}_{q}$ amounts
to giving a motive $X$ over $\mathbb{F}{}$ together with an endomorphism $%
\pi $ of $X$ that represents its Frobenius germ. The $\pi $ can be
considered a descent datum on $X$, and the condition on it a continuity
requirement. We now define the category of \emph{Weil motives }over $\mathbb{%
F}{}_{q}$, $\mathcal{M}(\mathbb{F}{}_{q};\mathbb{Z}{})_{W}$ to be the
category of pairs $(X,\pi )$ consisting of a motive $X$ over $\mathbb{F}{}$
and an endomorphism $\pi $ (\textquotedblleft noncontinuous descent
datum\textquotedblright ). Equivalently, $\mathcal{M}{}(\mathbb{F}{}_{q};%
\mathbb{Z}{})_{W}$ is the category of motives over $\mathbb{F}{}$ with a $%
\mathbb{Z}{}$-action. It contains $\mathcal{M}{}(\mathbb{F}{}_{q};\mathbb{Z}%
{})$ as a full subcategory.

Similarly, we define the category $\mathcal{M}^{+}(\mathbb{F}{}_{q};\mathbb{Z%
}{})_{W}$ of \emph{effective Weil motives}. For effective motives $X$ and $Y$
over $\mathbb{F}{}_{q}$, let%
\begin{equation*}
\Ext^{i}(X,Y)_{0}=\Ext_{\mathcal{M}{}^{+}(\mathbb{F}{}_{q};\mathbb{Z}%
{})_{W}}^{i}(X,Y)\text{.}
\end{equation*}

\begin{proposition}
\label{gc.b}For effective motives $X,Y$ over $\mathbb{F}{}_{q}$, there is a
spectral sequence%
\begin{equation*}
H^{i}(\Gamma _{0},\Ext_{\mathcal{M}{}^{+}(\mathbb{F}{};\mathbb{Z}{})}^{j}(%
\bar{X},\bar{Y}))\Longrightarrow \Ext^{i+j}(X,Y)_{0}
\end{equation*}
\end{proposition}

\begin{proof}
The proof of Theorem \ref{ha9} can be adapted to the present situation.
\end{proof}

\begin{lemma}
\label{gc.d}The groups $\Ext^{i}(X,Y)_{0}$ are finitely generated for all $i$%
, torsion for $i\geq 2$, and zero for $i\geq 3$.
\end{lemma}

\begin{proof}
The groups $\Ext_{\mathcal{M}{}^{+}(\mathbb{F}{};\mathbb{Z}{})}^{j}(\bar{X},%
\bar{Y})$ have composition series whose quotients are finitely generated $%
\mathbb{Z}{}$-modules or $\mathbb{F}{}$ with its natural $\Gamma _{0}$%
-action. Moreover, they are torsion for $i\geq 1$ and zero for $i\geq 2$.
Thus the statement follows from \ref{gc.b}.
\end{proof}

Let $f\colon \Hom(M,N)_{0}\rightarrow \Ext^{1}(M,N)_{0}$ be the map rendering

\begin{equation*}
\begin{CD} \Hom(M,N)_{0} @>f>> \Ext^{1}(M,N)_{0} \\ @VVV@AAA \\
\Hom_{\mathcal{M}^{+}(\mathbb{F};\mathbb{Z})}(\bar{M},\bar{N})^{\Gamma _{0}}
@>>>
\Hom_{\mathcal{M}{}^{+}(\mathbb{F}{};\mathbb{Z}{})}(\bar{M},\bar{N})_{\Gamma
_{0}}\end{CD}
\end{equation*}%
commutative. Here the vertical maps arise from a spectral sequence (\ref%
{gc.b}) and the lower map is induced by the identity map on $\Hom(\bar{M},%
\bar{N})$ (equal to cup-product with the canonical generator of $%
H^{1}(\Gamma _{0},\mathbb{Z}{})\cong \mathbb{Z}{}$).

\begin{theorem}
\label{gc.i}For effective motives $X$ and $Y$ over $\mathbb{F}{}_{q}$, $z(f)$
is defined and%
\begin{equation*}
q^{\chi (X,Y)}\cdot \zeta (X^{\vee }\otimes Y,s)\sim \pm z(f)\cdot \lbrack %
\Ext^{2}(X,Y)_{0}]\cdot (1-q^{-s})^{\rho (X,Y)}\text{ }\quad \text{as }%
s\rightarrow 0\text{.}
\end{equation*}%
(Notations as in \ref{gc.a}.)
\end{theorem}

\begin{proof}
If either $X$ or $Y$ is finite, then the spectral sequence \ref{gc.b}
coincides with \ref{ha9} and so%
\begin{equation*}
\Ext^{i}(X,Y)\cong \Ext^{i}(X,Y)_{0}\text{.}
\end{equation*}%
Now the same argument as in \S 8 proves that, for arbitrary $X,Y$, there are
exact sequences%
\begin{equation*}
0\rightarrow \Ext^{i}(X,Y)_{0}\otimes \mathbb{Z}{}_{l}\rightarrow \Ext%
^{i}(X_{l},Y_{l})\rightarrow T_{l}\Ext^{i+1}(X,Y)_{0}\rightarrow 0.
\end{equation*}%
Lemma \ref{gc.d} shows that $\Ext^{i+1}(X,Y)_{0}$ has no $l$-divisible
elements, and so these sequences reduce to isomorphisms%
\begin{equation*}
\Ext^{i}(X,Y)_{0}\otimes \mathbb{Z}{}_{l}\cong \Ext^{i}(X_{l},Y_{l})\text{.}
\end{equation*}%
The spectral sequence (\ref{gc.b}) is compatible with the spectral sequences
(\ref{lc.e}) and (\ref{lc.em}), and so the diagrams%
\begin{equation*}
\begin{CD} \Hom(X,Y)_{0}\otimes \mathbb{Z}_{l} @>f\otimes 1>>
\Ext^{1}(X,Y)_{0}\otimes \mathbb{Z}_{l} \\ @VV{\cong}V@VV{\cong}V \\
\Hom(X_{l},Y_{l}) @>{f_{l}}>> \Ext^{1}(X_{l},Y_{l})\end{CD}
\end{equation*}%
commute. Here $f_{l}$ is the map in Theorem \ref{lc.a}. Therefore, Theorem %
\ref{gc.i} follows from Theorem \ref{lc.a} and the product formula.
\end{proof}

\begin{remark}
\label{gc.e}The advantage of statement \ref{gc.i} over \ref{gc.a} is that it
has a natural extension to complexes of motives --- see the sequel to this
paper. Compare also Theorem \ref{lc.i}.
\end{remark}

\subsection{Motivic cohomology}

Let $V$ be a smooth projective variety\footnote{%
If the category $\mathcal{M}{}(\mathbb{F}{};\mathbb{Q}{})$ at the start of
this section is taken to be as in \ref{cm3} (case (b)), then $V$ must be
taken to be an abelian variety.} over $\mathbb{F}{}_{q}$. In contrast to the
pure case, $V$ should define, not a mixed motive, but rather a complex which
will not (in general) decompose when there is torsion, the torsion at $p$
being particularly complicated. Nevertheless, for an $r\geq 0$ such that the
truncated de Rham-Witt cohomology groups $H^{i}(V,W\Omega ^{\leq r-1})$ are
finitely generated over $W$, we make the ad hoc definition: $h^{i}V(r)$ is
the isomotive $h_{0}^{i}V(r)$ endowed with the $\mathbb{Z}{}$-structure
provided by the maps $H^{i}(V,\mathbb{Z}{}_{l}(r))\rightarrow H^{i}(V,%
\mathbb{Q}{}_{l}(r))$ (\'{e}tale cohomology for $l\neq p$ and crystalline
cohomology $l=p$). Define the (Weil) motivic cohomology groups by\footnote{%
According to Deligne 1994, 3.2.1, the motivic cohomology groups ($\mathbb{Q}%
{}$-coefficients) should be the final term of a spectral sequence%
\begin{equation*}
E_{2}^{ij}=\Ext^{i}(\1,h^{j}(V))\Longrightarrow H_{\text{mot}}^{i+j}(V,%
\mathbb{Q}{}(r))\text{.}
\end{equation*}%
Over $\mathbb{F}{}_{q}$ it is natural to expect this also with $\mathbb{Z}{}$%
-coefficients and (following Lichtenbaum's ideas) replace the Exts with the
Weil Exts:%
\begin{equation*}
E_{2}^{ij}=\Ext^{i}(\1,h^{j}(V))_{0}\Longrightarrow H_{\text{mot}}^{i+j}(V,%
\mathbb{\mathbb{Z}{}}{}(r))\text{.}
\end{equation*}%
If this spectral sequence degenerates, then we arrive at (\ref{e30}). In the
sequel to this paper, we shall take a less ad hoc approach.}%
\begin{equation}
H_{\text{mot}}^{j}(V,\mathbb{Z}{}(r))=\oplus _{i}\Ext^{i}(\1,h^{j-i}V(r))_{0}%
\text{.}  \label{e30}
\end{equation}%
Let
\begin{equation*}
H_{\text{mot}}^{\bullet }(V,\mathbb{Z}{}(r))=\cdots \rightarrow H_{\text{mot}%
}^{j}(V,\mathbb{Z}{}(r))\xr{f^j}H_{\text{mot}}^{j+1}(V,\mathbb{Z}%
{}(r))\rightarrow \cdots
\end{equation*}%
be the complex with%
\begin{equation*}
\begin{diagram} H_{\text{mot}}^{j}(V,Z(r))&\quad&\Ext^{0}(\1,h^{j}(V)(r))_0
& \oplus & \Ext^{1}(\1,h^{j-1}(V)(r))_0 & \oplus &
\Ext^{2}(\1,h^{j-2}(V)(r))_0\\ \dTo^{f^j}&=& &\rdTo^{f} & & \rdTo^{f} & \\
H_{\text{mot}}^{j+1}(V,Z(r))&\quad&\Ext^{0}(\1,h^{j+1}(V)(r))_0 &\oplus &
\Ext^{1}(\1,h^{j}(V)(r))_0 &\oplus & \Ext^{2}(\1,h^{j-1}(V)(r))_0
\end{diagram}
\end{equation*}%
The maps $f$ at right are as in Theorem \ref{gc.i}.

\begin{theorem}
\label{gc.n}Let $V$ be a smooth projective variety over $\mathbb{F}{}_{q}$,
and let $r\geq 0$ be such that the groups $H^{i}(V,W\Omega ^{\leq r-1})$ are
finitely generated over $W$.

\begin{enumerate}
\item The groups $H_{\text{mot}}^{i}(V,\mathbb{Z}{}(r))$ are finitely
generated abelian groups.

\item The alternating sum $\tsum (-1)^{i}r_{i}$ of the ranks $r_{i}$ of the
groups $H_{\text{mot}}^{i}(V,\mathbb{Z}{}(r))$ is zero.

\item The order of the zero of $\zeta (V,s)$ at $s=r$ is equal to the
secondary Euler characteristic $\rho (V,r)=_{\text{df}}\tsum (-1)^{i}ir_{i}$.

\item The cohomology groups of the complex $(H_{\text{mot}}^{\bullet }(V,%
\mathbb{Z}{}(r)),f)$ are finite, and the alternating product of their orders
$\chi ^{\times }(V,\mathbb{Z}{}(r))$ satisfies%
\begin{equation*}
\zeta (V,s)\sim \pm \chi ^{\times }(V,\mathbb{Z}{}(r))\cdot q^{\chi (V,%
\mathcal{O}{},r)}\cdot (1-q^{r-s})^{\rho }\text{ as }s\rightarrow r
\end{equation*}%
where
\begin{equation*}
\chi (V,\mathcal{O}{},r)=\sum_{%
\begin{smallmatrix}
0\leq i\leq r \\
0\leq j\leq \dim V%
\end{smallmatrix}%
}(-1)^{i+j}(r-i)\cdot \dim _{\mathbb{F}{}_{q}}H^{j}(X,\Omega ^{i}).
\end{equation*}
\end{enumerate}
\end{theorem}

\begin{proof}
Immediate consequence of Theorem \ref{gc.i} and Milne 1986a, 4.1 (the
condition on the de Rham-Witt groups implies that the numbers $d^{i}(r-1)$
in the second reference are zero).
\end{proof}

\begin{remark}
\label{gc.o}When $r=0$, the condition on the de Rham-Witt groups is
satisfied vacuously, and the groups $H_{\text{mot}}^{i}(V,\mathbb{Z}{}(0))$
coincide with the Weil-\'{e}tale groups $H_{W}^{i}(V,\mathbb{Z}{})$ of
Lichtenbaum 2002. Thus, in this case the theorem coincides with the smooth
projective case of ibid., Theorem 8.2.
\end{remark}

\section*{References}

%TCIMACRO{\TeXButton{footnotesize}{\footnotesize}}%
%BeginExpansion
\footnotesize%
%EndExpansion

Andr\'e, Yves, Pour une th\'eorie inconditionnelle des motifs.
Inst. Hautes \'Etudes Sci. Publ. Math. No. 83, (1996), 5--49.

Andr\'e, Yves, and Kahn, Bruno, Construction inconditionelle de
groupes de Galois motiviques, Preprint, October 2001.

%Artin, M., Grothendieck, A., and Verdier J.L., Th\'eorie des topos et cohomologie \'etale des sch\'emas. S\'eminaire de G\'eom\'etrie Alg\'ebrique du Bois-Marie 1963--1964 (SGA 4). Avec la collaboration de N. Bourbaki, P. Deligne et B. Saint-Donat. Lecture Notes in Mathematics, Vol. 269, 270, 305. Springer-Verlag, Berlin-New York, 1972.

Atiyah, M. F., and Hirzebruch, F., Analytic cycles on complex
manifolds.  Topology, 1, 1962, 25--45.

%Atiyah, M. F. and Macdonald, I. G., Introduction to commutative algebra. Addison-Wesley Publishing Co., Reading, Mass.-London-Don Mills, Ont. 1969.

Ballico, E., Ciliberto, C., and Catanese F., Trento examples. In:
Classification of irregular varieties. Minimal models and abelian
varieties. Proceedings of the conference held in Trento, December
17--21, 1990. Edited by E. Ballico, F. Catanese and C. Ciliberto.
Lecture Notes in Mathematics, 1515. Springer-Verlag, Berlin, 1992,
p. 134--139.

%Beilinson, A. A. Notes on absolute Hodge cohomology. Applications of algebraic $K$-theory to algebraic geometry and number theory, Part I, II (Boulder, Colo., 1983), 35--68, Contemp. Math., 55, Amer. Math. Soc., Providence, R.I., 1986.

%Beilinson, A. A., On the derived category of perverse sheaves. $K$-theory, arithmetic and geometry (Moscow, 1984--1986), 27--41, Lecture Notes in Math., 1289, Springer, Berlin-New York, 1987.

%Beilinson, A. A.; Bernstein, J.; Deligne, P., Faisceaux pervers. Analysis and topology on singular spaces, I (Luminy, 1981), 5--171, Ast\'erisque, 100, Soc. Math. France, Paris, 1982.

%Berthelot, P., Alt\'erations de vari\'et\'es alg\'ebriques (d'apr\`es A. J. de Jong), Ast\'erisque No. 241 (1997), Exp.\ No.\ 815, 5, 273--311.

Bloch, Spencer, and Esnault, H\'el\`ene, The coniveau filtration
and non-divisibility for algebraic cycles. Math. Ann. 304 (1996),
no. 2, 303--314.

Brosnan, P. 1999. Steenrod Operations in Chow Theory, preprint,
www.math.uiuc.edu/K-theory/0370.

%Bruhat, F.; Tits, J. Groupes r\'eductifs sur un corps local. II. Sch\'e mas en groupes. Existence d'une donn\'ee radicielle valu\'ee. Inst. Hautes \'Etudes Sci. Publ. Math. No. 60 (1984), 197--376.

%Deligne, Pierre, La conjecture de Weil. I. Inst. Hautes \'Etudes Sci. Publ. Math. No. 43, (1974a), 273--307.

%Deligne, Pierre, Th\'eorie de Hodge. III. Inst. Hautes \'Etudes Sci. Publ. Math. No. 44, (1974b), 5--77.

%Deligne, Pierre, et al., Cohomologie \'etale. (S\'eminaire de G\'eom \'etrie Alg\'ebrique du Bois-Marie SGA 4$\frac{1}{2}$.) Lecture Notes in Mathematics, Vol. 569. Springer-Verlag, Berlin-New York, 1977.

%Deligne, Pierre, La conjecture de Weil. II. Inst. Hautes \'Etudes Sci. Publ. Math. No. 52, (1980), 137--252.

%Deligne, Pierre, (Notes by J.S. Milne), Hodge cycles on abelian varieties, in Hodge cycles, motives, and Shimura varieties. Lecture Notes in Mathematics, 900. Springer-Verlag, Berlin-New York, 1982, 9--100.

Deligne, P., Le groupe fondamental de la droite projective moins
trois points. Galois groups over $\mathbb{Q}{}$ (Berkeley, CA,
1987), 79--297, Math. Sci. Res. Inst. Publ., 16, Springer, New
York-Berlin, 1989.

Deligne, P., Cat\'egories tannakiennes. The Grothendieck
Festschrift, Vol. II, 111--195, Progr. Math., 87, Birkh\"{a}user
Boston, Boston, MA, 1990.

%Deligne, Pierre, D\'ecompositions dans la cat\'egorie d\'eriv\'ee. Motives (Seattle, WA, 1991), 115--128, Proc. Sympos. Pure Math., 55, Part 1, Amer. Math. Soc., Providence, RI, 1994a.

Deligne, Pierre, \`{A} quoi servent les motifs? Motives (Seattle,
WA, 1991), 143--161, Proc. Sympos. Pure Math., 55, Part 1, Amer.
Math. Soc., Providence, RI, 1994.

Deligne, Pierre, and Milne, James, Tannakian Categories. In: Hodge
cycles, motives, and Shimura varieties. Lecture Notes in
Mathematics, 900. Springer-Verlag, Berlin-New York, 1982,
101--228.

Demazure, Michel, Lectures on $p$-divisible groups. Lecture Notes
in Mathematics, 302. Springer-Verlag, Berlin-New York, 1972.

Dwyer, William G., and Friedlander, Eric M., Algebraic and etale
$K$-theory. Trans. Amer. Math. Soc. 292 (1985), no. 1, 247--280.

Ekedahl, Torsten, Diagonal complexes and $F$-gauge structures.
Travaux en Cours. Hermann, Paris, 1986.

%Ekedahl, Torsten, On the adic formalism. The Grothendieck Festschrift, Vol. II, 197--218, Progr. Math., 87, Birkh\"{a}user Boston, Boston, MA, 1990.

%Euclid, -300, $\Sigma \tau o\iota \chi \epsilon \tilde{\iota}\alpha $.

Faltings, Gerd, and W\"{u}stholz, Gisbert (Eds). Rational points.
Papers from the seminar held at the Max-Planck-Institut f\"{u}r
Mathematik, Bonn, 1983/1984. Aspects of Mathematics, E6. Friedr.
Vieweg \& Sohn, Braunschweig; distributed by Heyden \& Son, Inc.,
Philadelphia, Pa., 1984.

%Fan, Yun, and Xu, Bangteng, On the $2$-category of triangulated categories and its $t$-split exact sequences. Bull. Hong Kong Math. Soc. 2 (1998), no. 1, 99--113.

Fontaine, Jean-Marc, and Mazur, Barry Geometric Galois
representations. Elliptic curves, modular forms, \& Fermat's last
theorem (Hong Kong, 1993), 41--78, Ser. Number Theory, I,
Internat. Press, Cambridge, MA, 1995.

Gabber, Ofer, Sur la torsion dans la cohomologie $l$-adique d'une
vari\'et \'e. C. R. Acad. Sci. Paris S\'er. I Math. 297 (1983),
no. 3, 179--182.

%Gabriel, Pierre, Des cat\'egories ab\'eliennes. Bull. Soc. Math. France 90 (1962) 323--448.

%Gabriel, P., and Zisman, M., Calculus of fractions and homotopy theory. Ergebnisse der Mathematik und ihrer Grenzgebiete, Band 35 Springer-Verlag New York, Inc., New York 1967.

%Gelfand, Sergei I., and Manin, Yuri I., Methods of homological algebra. Translated from the 1988 Russian original. Springer-Verlag, Berlin, 1996.

Griffiths, Phillip, and Harris, Joe, On the Noether-Lefschetz
theorem and some remarks on codimension-two cycles. Math. Ann. 271
(1985), no. 1, 31--51.

%Grothendieck, Alexander, Sur quelques points d'alg\`ebre homologique. T \^{o}hoku Math. J. (2) 9 1957 119--221.

%Grothendieck, Alexander, Techniques de construction et th\'eor\`emes d'existence en g\'eom\'etrie alg\'ebrique, III: Pr\'esch\'emas quotients, S\'eminaire Bourbaki, 1960/61, 212.

%Grothendieck, Alexander, Letter to Illusie, dated May 3, 1973 (In: Motives (Seattle, WA, 1991), 296--300, Proc. Sympos. Pure Math., 55, Part 1, Amer. Math. Soc., Providence, RI, 1994).

%Grothendieck, Alexander, et al., Cohomologie $l$-adique et fonctions $L$. (S \'eminaire de G\'eom\'etrie Alg\'ebrique du Bois-Marie 1965-66, SGA 5.), Lecture Notes in Mathematics 589, Springer-Verlag, Berlin-Heidelberg-New York, 1977.

Hanamura, Masaki Mixed motives and algebraic cycles. I. Math. Res.
Lett. 2 (1995), no. 6, 811--821.

Hanamura, Masaki Mixed motives and algebraic cycles. III. Math.
Res. Lett. 6 (1999), no. 1, 61--82.

Hiller, Howard L., $\lambda $-rings and algebraic $K$-theory. J.
Pure Appl. Algebra 20 (1981), no. 3, 241--266.

Huber, Annette, Calculation of derived functors via
Ind-categories. J. Pure Appl. Algebra 90 (1993), no. 1, 39--48.

%Huber, Annette, Mixed motives and their realizations in derived categories, Lecture Notes in Math., 1604, Springer, Berlin-New York, 1995.

%Illusie, L., Finiteness, duality, and K\"{u}nneth theorems in the cohomology of the de Rham Witt complex. Algebraic geometry (Tokyo/Kyoto, 1982), 20--72, Lecture Notes in Math., 1016, Springer, Berlin-New York, 1983.

Jacobson, Nathan, The Theory of Rings. American Mathematical
Society Mathematical Surveys, vol. I. American Mathematical
Society, New York, 1943.

%Jannsen, Uwe, Continuous \'etale cohomology. Math. Ann. 280 (1988), no. 2, 207--245.

Jannsen, Uwe, Motives, numerical equivalence, and semi-simplicity.
Invent. Math. 107 (1992), 447--452.

%Jannsen, Uwe, Mixed motives, motivic cohomology, and Ext-groups. Proceedings of the International Congress of Mathematicians, Vol. 1, 2 (Z\"{u}rich, 1994), 667--679, Birkh\"{a}user, Basel, 1995.

Jensen, C. U. Les foncteurs d\'eriv\'es de $\varprojlim$ et leurs
applications en th\'eorie des modules. Lecture Notes in
Mathematics, Vol. 254. Springer-Verlag, Berlin-New York, 1972.

de Jong, A. J., Smoothness, semi-stability and alterations. Inst.
Hautes \'Etudes Sci. Publ. Math. No. 83, (1996), 51--93.

de Jong, A. J. Homomorphisms of Barsotti-Tate groups and crystals
in positive characteristic. Invent. Math. 134 (1998), no. 2,
301--333, Erratum ibid. 138 (1999), no. 1, 225.

%Katz, Nicholas M. and Messing, William, Some consequences of the Riemann hypothesis for varieties over finite fields. Invent. Math. 23 (1974), 73--77.

Katz, Nicholas M., Review of $l$-adic cohomology. Motives
(Seattle, WA, 1991), 21--30, Proc. Sympos. Pure Math., 55, Part 1,
Amer. Math. Soc., Providence, RI, 1994.

Kleiman, Steven L. The standard conjectures. Motives (Seattle, WA,
1991), 3--20, Proc. Sympos. Pure Math., 55, Part 1, Amer. Math.
Soc., Providence, RI, 1994.

Kratzer, Ch., $\lambda $-structure en $K$-th\'eorie alg\'ebrique.
Comment. Math. Helv. 55 (1980), no. 2, 233--254.

%Laumon, G\'erard, Sur la cat\'egorie d\'eriv\'ee des $\mathcal{D}$ -modules filtr\'es, Algebraic geometry (Tokyo/Kyoto, 1982), 151--237, Lecture Notes in Math., 1016, Springer, Berlin-New York, 1983.

 Lang, Serge, and Tate, John, Principal homogeneous spaces over abelian varieties. Amer. J. Math. 80 1958 659--684.

Levine, Marc, Mixed motives. Mathematical Surveys and Monographs,
57. American Mathematical Society, Providence, RI, 1998.

%Lichtenbaum, Stephen, Zeta functions of varieties over finite fields at $s=1$. Arithmetic and geometry, Vol. I, 173--194, Progr. Math., 35, Birkh\"{a}user Boston, Boston, Mass., 1983.

%Lichtenbaum, Stephen, Motivic complexes. Motives (Seattle, WA, 1991), 303--313, Proc. Sympos. Pure Math., 55, Part 1, Amer. Math. Soc., Providence, RI, 1994.

Lichtenbaum, Stephen, The Weil-\'etale topology. Preprint (2002).

%Mac Lane, Saunders, Categories for the working mathematician. Second edition. Graduate Texts in Mathematics, 5. Springer-Verlag, New York, 1998.

%Matsumura, Hideyuki, Commutative algebra. W. A. Benjamin, Inc., New York 1970.

Milne, J. S., Extensions of abelian varieties defined over a
finite field. Invent. Math. 5 (1968), 63--84.

%Milne, J. S., The homological dimension of commutative group schemes over a perfect field. J. Algebra 16 1970 436--441.

Milne, J. S., \'Etale cohomology. Princeton Mathematical Series,
33. Princeton University Press, Princeton, N.J., 1980.

Milne, J. S., Values of zeta functions of varieties over finite
fields. Amer. J. Math. 108 (1986a), no. 2, 297--360.

Milne, J. S., Arithmetic duality theorems. Perspectives in
Mathematics, 1. Academic Press, Inc., Boston, Mass., 1986b.

%Milne, J.S., Motivic cohomology and values of zeta functions, Compos. Math. 68 (1988), 59-102.

Milne, J. S., Motives over finite fields. Motives (Seattle, WA,
1991), 401--459, Proc. Sympos. Pure Math., 55, Part 1, Amer. Math.
Soc., Providence, RI, 1994.

%Milne, J. S., Lefschetz classes on abelian varieties, Duke Math. J. 96 (1999a), pp. 639-675.

Milne, J. S., Lefschetz motives and the Tate conjecture,
Compositio Math. 117 (1999), pp. 47-81.

%Milne, J. S., The Tate conjecture for certain abelian varieties over finite fields. Preprint, August 1, 1999c, 28pp, arXiv:math.NT/9911218 (to appear in Acta Arith.).

Milne, J. S., Motives for rational Tate classes (preprint 2002, in
preparation).

Mitchell, Barry, Theory of categories. Pure and Applied
Mathematics, Vol. XVII Academic Press, New York-London 1965.

%Mumford, David Geometric invariant theory. Ergebnisse der Mathematik und ihrer Grenzgebiete, Neue Folge, Band 34 Academic Press, Inc., New York; Springer-Verlag, Berlin-New York 1965

%Neeman, Amnon, Some new axioms for triangulated categories. J. Algebra 139 (1991), no. 1, 221--255.

Quillen, Daniel, On the cohomology and $K$-theory of the general
linear groups over a finite field. Ann. of Math. (2) 96 (1972),
552--586.

Oort, F., Yoneda extensions in abelian categories. Math. Ann. 153
1964 227--235.

%Raskind, Wayne, Algebraic $K$-theory, \'etale cohomology and torsion algebraic cycles. Algebraic $K$-theory and algebraic number theory (Honolulu, HI, 1987), 311--341, Contemp. Math., 83, Amer. Math. Soc., Providence, RI, 1989.

Saavedra Rivano, Neantro, Cat\'egories Tannakiennes. Lecture Notes
in Mathematics, Vol. 265. Springer-Verlag, Berlin-New York, 1972.

Schoen, Chad, An integral analog of the Tate conjecture for
one-dimensional cycles on varieties over finite fields. Math. Ann.
311 (1998), no. 3, 493--500.

%Scholl, A. J., Classical motives. Motives (Seattle, WA, 1991), 163--187, Proc. Sympos. Pure Math., 55, Part 1, Amer. Math. Soc., Providence, RI, 1994.

Serre, Jean-Pierre, Groupes de Grothendieck des sch\'emas en
groupes r\'eductifs d\'eploy\'es. Inst. Hautes \'Etudes Sci. Publ.
Math. No. 34 (1968), 37--52.

%Shimura, Goro, Introduction to the arithmetic theory of automorphic functions. Kan\^{o} Memorial Lectures, No. 1. Publications of the Mathematical Society of Japan, No. 11. Iwanami Shoten, Publishers, Tokyo; Princeton University Press, Princeton, N.J., 1971.

Soul\'e, Christophe, $K$-th\'eorie des anneaux d'entiers de corps
de nombres et cohomologie \'etale. Invent. Math. 55 (1979), no. 3,
251--295.

%Soul\'e, Christophe, Groupes de Chow et $K$-th\'eorie de vari\'et\'e s sur un corps fini. Math. Ann. 268 (1984), no. 3, 317--345.

Soul\'e, Christophe, $K$-theory and values of zeta functions.
Algebraic $K$ -theory and its applications (Trieste, 1997),
255--283, World Sci. Publishing, River Edge, NJ, 1999.

Tate, John, Endomorphisms of abelian varieties over finite fields.
Invent. Math. 2 1966a, 134--144.

Tate, John, On the conjectures of Birch and Swinnerton-Dyer and a
geometric analog, S\'eminaire Bourbaki: Vol. 1965/66, Expose 306,
1966b.

Tate, John, Relations between $K_{2}$ and Galois cohomology.
Invent. Math. 36 (1976), 257--274.

Tate, John, Conjectures on algebraic cycles in $\ell $-adic
cohomology. Motives (Seattle, WA, 1991), 71--83, Proc. Sympos.
Pure Math., 55, Part 1, Amer. Math. Soc., Providence, RI, 1994.

Verdier, Jean-Louis, Des cat\'egories d\'eriv\'ees des
cat\'egories ab\'eliennes. Ast\'erisque No. 239 , 1996.

Voevodsky, V., Triangulated categories and motives over a field,
in Cycles, transfers, and motivic homology theories. Annals of
Mathematics Studies, 143. Princeton University Press, Princeton,
NJ, 2000, 188--238.

%Waterhouse, W. C.; Milne, J. S., Abelian varieties over finite fields. 1969 Number Theory Institute (Proc. Sympos. Pure Math., Vol. XX, State Univ. New York, Stony Brook, N.Y., 1969), pp. 53--64. Amer. Math. Soc., Providence, R.I., 1971.

%Weibel, Charles, \'Etale Chern classes at the prime $2$. Algebraic $K$ -theory and algebraic topology (Lake Louise, AB, 1991), 249--286, NATO Adv. Sci. Inst. Ser. C Math. Phys. Sci., 407, Kluwer Acad. Publ., Dordrecht, 1993.

Weil, Andr\'e, Adeles and algebraic groups. With appendices by M.
Demazure and Takashi Ono. Progress in Mathematics, 23.
Birkh\"{a}user, Boston, Mass., 1982.

%Yan, F., Tate property and isogeny estimate for semi-abelian varieties, thesis, ETH Z\"{u}rich, 1994.

%Yao, Dongyuan, On equivalence of derived categories. $K$-Theory 10 (1996), no. 3, 307--322.

%Zarhin, Ju. G., Isogenies of abelian varieties over fields of finite characteristic. (Russian) Mat. Sb. (N.S.) 95(137) (1974a), 461--470, 472.

%Zarhin, Ju. G., A finiteness theorem for isogenies of abelian varieties over function fields of finite characteristic. (Russian) Funkcional. Anal. i Prilov zen. 8 (1974b), no. 4, 31--34.

%Zarhin, Ju. G., A remark on endomorphisms of abelian varieties over function fields of finite characteristic. (Russian) Izv. Akad. Nauk SSSR Ser. Mat. 38 (1974c), 471--474.

Zarhin, Ju. G., Endomorphisms of Abelian varieties over fields of
finite characteristic. (Russian) Izv. Akad. Nauk SSSR Ser. Mat. 39
(1975), no. 2, 272--277, 471.

%Zarhin, Ju. G. Abelian varieties in characteristic $p$. (Russian) Mat. Zametki 19 (1976), no. 3, 393--400.

\end{document}